\documentclass[a4paper,fleqn]{cas-sc}

\usepackage[numbers]{natbib}

\newtheorem{cj}{Conjectures}[section]
\newtheorem{summary}{Summary (secondary bifurcations)}[section]
\def\tsc#1{\csdef{#1}{\textsc{\lowercase{#1}}\xspace}}
\tsc{WGM}
\tsc{QE}

\newtheorem{theorem}{Theorem}

\newdefinition{remk}{Remark}
\newproof{proof}{Proof}

\begin{document}
\let\WriteBookmarks\relax
\def\floatpagepagefraction{1}
\def\textpagefraction{.001}

\shorttitle{{The 1D nonlocal Fisher-KPP equation with a top hat kernel. Part 3. The effect of perturbations in the kernel}}    

\shortauthors{D. J. Needham and J. Billingham}  

\title [mode = title]{{The evolution problem for the 1D nonlocal Fisher-KPP equation with a top hat kernel. Part 3. The effect of perturbations in the kernel}}  

\tnotemark[1] 


%

\author[1]{David John Needham}



\ead{d.j.needham@bham.ac.uk}


\credit{Conceptualization, Methodology, Formal Analysis, Verification, Writing - Original Draft, Writing - Reviewing and Editing}

\affiliation[1]{organization={University of Birmingham},
addressline={School of Mathematics, Watson Building, University of Birmingham, Edgbaston}, 
city={Birmingham},
postcode={B15 2TT}, 
country={UK}}
\affiliation[2]{organization={University of Nottingham},
addressline={School of Mathematical Sciences, University of Nottingham}, 
city={Nottingham},
postcode={NG7 2RD}, 
country={UK}}

\author[2]{John Billingham}



\ead{john.billingham@nottingham.ac.uk}

\ead[url]{}

\credit{Conceptualization, Computation, Verification, Writing - Reviewing and Editing}

\cortext[1]{Corresponding author}


\nonumnote{The evolution problem for the 1D nonlocal Fisher-KPP equation with a top hat kernel. Part 3. The effect of perturbations in the kernel}

\begin{abstract}
In the third part of this series of papers, we address the same Cauchy problem that was considered in part 1 (see \cite{NBLM}), namely the nonlocal Fisher-KPP equation in one spatial dimension, $u_t = D u_{xx} + u(1-\phi_T*u)$, where $\phi_T*u$ is a spatial convolution with the top hat kernel, $\phi_T(y) \equiv H\left(\frac{1}{4}-y^2\right)$, except that now we include a specified perturbation to this kernel, which we denote as $\overline{\phi}:\mathbb{R}\to \mathbb{R}$. Thus the top hat kernel $\phi_T$ is now replaced by the perturbed kernel $\phi:\mathbb{R} \to \mathbb{R}$, where $\phi(x) = \phi_T(x) + \overline{\phi}(x)~~\forall~~x\in \mathbb{R}$. When the magnitude of the kernel perturbation is small in a suitable norm, the situation is shown to be generally a regular perturbation problem when the diffusivity $D$ is formally of O(1) or larger. However when $D$ becomes small, and in particular, of the same order as the magnitude of the perturbation to the kernel, this becomes a strongly singular perturbation problem, with considerable changes in overall structure. This situation is uncovered in detail In terms of its generic interest, the model forms a natural extension to the classical Fisher-KPP model, with the introduction of the simplest possible nonlocal effect into the saturation term. Nonlocal reaction-diffusion models arise naturally in a variety of (frequently biological or ecological) contexts, and as such it is of fundamental interest to examine its properties in detail, and to compare and contrast these with the well known properties of the classical Fisher-KPP model.
\end{abstract}


\begin{highlights}
\item The robustness of the qualitative and quantitative dynamical features of the nonlocal Fisher-KPP equation to a broad class of perturbations of a top hat kernel is shown to be a regular perturbation when the diffusivity coefficient, $D = O(1)$, but a singular perturbation when $D = o(1)$, in relation to the magnitude of the perturbation.
\item Details in both cases are fully unfolded.
\item Two specific and representative forms of perturbation are considered in detail, and in particular when the perturbation is singular. This uncovers a complicated network of secondary bifurcation structures.
\end{highlights}


\maketitle


\section{Introduction}\label{sec_intro}
In this third paper we consider the evolution problem detailed in part 1 of this series of papers \cite{NBLM} (henceforth referred to as (NB)) for the nonlocal Fisher-KPP equation with top hat kernel, but now we augment the top hat kernel by including a specified perturbation to this kernel. Specifically, the top hat kernel, henceforth denoted by $\phi_T:\mathbb{R}\to \mathbb{R}$, is replaced by the perturbed kernel $\phi \equiv (\phi_T + \overline{\phi}) :\mathbb{R}\to \mathbb{R}$, and the associated nonlocal Fisher-KPP equation becomes,
\begin{equation}
    u_t = Du_{xx} + u\left(1 - \int_{x-\frac{1}{2}}^{x+\frac{1}{2}}{u(y,t)dy} - \overline{\phi}\ast u\right)~~\text{with}~~(x,t)\in D_{\infty}\equiv \mathbb{R}\times \mathbb{R}^+. \label{eqn1.1}
\end{equation}
The perturbation to the kernel, $\overline{\phi}: \mathbb{R}\to \mathbb{R}$ is limited so that the perturbed kernel $(\phi_T + \overline{\phi})$ is everywhere non-negative, which we refer to as admissible, and restricted to satisfy the following general conditions, 
\begin{itemize}
    \item[(P1)]~~$\overline{\phi}\in PC^2(\mathbb{R}) \cap C(\mathbb{R}\setminus \{-\frac{1}{2}, \frac{1}{2}\}) \cap L^{\infty}(\mathbb{R}) \cap L^1(\mathbb{R})$
    \item[(P2)]~~$\text{There exists}~R \in \mathbb{N}$ such that
        $\int_{-\infty}^{\infty}{|y|^R|\overline{\phi}|(y)}dy < \infty$
    \item[(P3)]~~$\overline{\phi}(y) = \overline{\phi}(-y)~~\forall~~y\in \mathbb{R}$
    \item[(P4)]~~$\overline{\phi}(y)\to 0$ as $y\to \infty$
    \item[(P5)]~~$\int_{-\infty}^{\infty}{\overline{\phi}(y)}dy = 0$ 
\end{itemize}
We denote the set of admissible kernel perturbations which satisfy the conditions $(P1)-(P5)$ as $\overline{K}(\mathbb{R})$. For convenience in later discussion, we will refer to that subset of admissible kernel perturbations in $\overline{K}(\mathbb{R})$, which render the full kernel everywhere continuous, as \emph{admissible smoothing kernel perturbations}. The object of the paper is to study the evolution problem specified as (IBVP) in (NB), but now with a nontrivial admissible kernel perturbation $\overline{\phi}\in \overline{K}(\mathbb{R})$. Particular attention is paid to the robustness, or otherwise, of the qualitative and quantitative features identified for (IBVP) with top hat kernel (as reported in (NB)), for perurbed kernels when the perturbation has a suitably small but nonzero magnitude. In this respect, we will find it convenient to measure the magnitude of an admissible kernel perturbation $\overline{\phi}\in \overline{K}(\mathbb{R})$ in the norm $||\cdot||_1^m: \overline{K}(\mathbb{R}) \to [0,\infty)$, where,
\begin{equation}
 ||\overline{\phi}||_1^m = \text{max} \{||\overline{\phi}||_1,~||(\cdot) \overline{\phi}(\cdot)||_1 \}.   
\end{equation}
We remark that there are admissible smoothing kernel perturbations which have $||\overline{\phi}||_1^m$ arbitrarilly small. We will now refer to this evolution problem throughout as (IBVP)$_p$, which is precisely as given in (NB), except now with the perturbed kernel replacing the top hat kernel.

The paper should be read in conjunction with (NB), where the context for studying the problem at hand is set, much of the theory relied upon here is developed and its context in relation to other relevant works in discussed and expounded (for brevity we do not repeat this here, but note that the following works,......, in particular, are contextualised in (NB)). The approach adopted in the paper is first to fix a given admissible kernel perturbation $\overline{\phi}\in \overline{K}(\mathbb{R})$, with $||\overline{\phi}||_1^m$ small, and then consider each of the key features identified in (NB) in relation to (IBVP), investigating their robustness, both qualitatively and quantitatively, in relation to the perturbed problem (IBVP)$_p$, and how this robustness depends upon the key parameter $D$. Firstly, without repeating the details, it is readily established that all of the well-posed (\emph{with respect to initial data}) and basic qualitative results reported in section 2 of (NB) for (IBVP) with the top hat kernel, carry over without modification, to (IBVP)$_p$ with the admissible perturbed kernel, when $\overline{\phi}\in \overline{K}(\mathbb{R})$. The key features which we address, and the associated results which we establish, are summarised as follows: \\\\
$\mathbf{In~Section~2}$ we confirm that, for any admissible kernel perturbation $\overline{\phi}\in \overline{K}(\mathbb{R})$, the equilibrium states associated with (IBVP)$_p$ remain as the \emph{unreacted} state $u=0$ and the \emph{fully reacted} state $u=1$. We then address the linearised stability and dispersion relation associated with each of these equilibrium states. Concerning the equilibrium state $u=0$, the linearised operator associated with the perturbed problem (IBVP)$_p$ is identical to that associated with the unperturbed problem (IBVP), and so the linearised stability characteristics and the dispersion relation remain unchanged. This is not the case when considering the linearised stability and the corresponding dispersion relation associated with the equilibrium state $u=1$. However, we readily establish that in the non-negative quadrant of the $(k,D)$-plane, with $k$ representing wavenumber, the neutral curve corresponding to (IBVP)$_p$ is uniformly close to the neutral curve for (IBVP) \emph{in both location and slope}, with the displacments bounded by an error of $O(||\overline{\phi}||_1^m)$ as $||\overline{\phi}||_1^m \to 0$, \emph{uniformly} throughout the non-negative quadrant of the $(k,D)$-plane. It can therefore be concluded that the linearised stability characteristics of the equilibrium state $u=1$ associated with (IBVP)$_p$ are simply a regular perturbation from the corresponding linearised stability characteristics for the equilibrium state $u=1$ associated with (IBVP), when $||\overline{\phi}||_1^m$ is small, with this regular perturbation remaining uniform throughout $D>0$. The qualitative and quantitative estimates involved are developed throughout section 2, and summarised in detail at the end of the section, where it is concluded that the fundamental conjectures [P1] and [P2] formulated in (NB) (at the end of section 3) for (IBVP) will continue to hold for (IBVP)$_p$ (uniformly in $D>0$) when $||\overline{\phi}||_1^m$ is sufficiently small.\\\\
$\mathbf{In~Section~3}$ we begin by seeking to identify steady state bifurcations from the equilibrium state $u=1$ which give rise to positive periodic states associated with (IBVP)$_p$. For clarity (and avoiding repetition) we adopt all of the notation introduced in sections 3 and 4 of (NM) for (IBVP), when referring to (IBVP)$_p$, except where confusion may arise, when we introduce a superscript $p$ to indicate that (IBVP)$_p$ is being referred to. In the positive quadrant of the $(\lambda,D)$-plane (with $\lambda$ again representing fundamental wavelength) we demonstrate that, as in the unperturbed case (IBVP) (and detailed in section 4 of (NM)), there remains a smooth bifurcation locus, with a countable multiple of disjoint smooth components, across which there is a steady state pitchfork bifurcation generating a unique (up to spatial translation) nontrivial positive periodic steady state. This bifurcation locus is contained, in the positive quadrant of the $(\lambda,D)$-plane, within the rectangle $(0,\lambda_p] \times (0,\Delta_1^p]$, where,
\begin{equation}
    \lambda_p = 1 + O(||\overline{\phi}||_1^m)
\end{equation}
and
\begin{equation}
    \Delta_1^p = \Delta_1 + (||\overline{\phi}||_1^m) 
\end{equation}
as $||\overline{\phi}||_1^m \to 0$. Moreover, the bifurcation locus is uniformly close to that of (IBVP), in both location and slope, with the displacement being uniformly of $O(||\overline{\phi}||_1^m)$ as $||\overline{\phi}||_1^m \to 0$ throughout $(0,\lambda_p] \times (0,\Delta_1^p]$, whilst the criticallity of the bifurcation remains unchanged. Thus, as in (NB), the bifurcation locus is again in the form of a sequence of 'tongue like' curves based on the $\lambda$-axis, with decreasing height, and with each 'tongue' located within an $O(||\overline{\phi}||_1^m)$ neighbourhood of the corresponding tongue, $\partial\Omega_i (i=1,2,..)$, for the case of the unperturbed top hat kernel. For simplicity, we adopt directly, and without adaptation, the notion for, and relating to, these 'tongues` introduced in (NB), always bearing in mind, from above, that associated locations in the $(\lambda,D)$-plane are perturbed by $O(||\overline{\phi}||_1^m)$ when the kernel is perturbed with $||\overline{\phi}||_1^m$ small. A straightforward, but lengthy, calculation also establishes  that, at each bifurcation point, the nature and criticality of the bifurcation is unchanged when $||\overline{\phi}||_1^m$ is small, with each bifurcation again being a pitchfork bifurcation, generating a unique (up to translation in $x$) periodic steady state with small amplitude relative to the equilibrium state $u=1$, on traversing through the bifurcation point in the $(\lambda,D)$-plane, from the \emph{exterior to the interior} of the tongue. As in (NB), translation invariance can be most conveniently fixed by selecting that representative which is an even function of $x$.
With the unperturbed top hat kernel, we demonstrated in (NB) that at each point, within each tongue, there is a unique, even and positive periodic steady state, which we labelled as $u = F_p(x,\lambda,D)$ for $x\in \mathbb{R}$ and $(\lambda,D)\in \Omega$. At each such point $(\lambda,D)\in \Omega$, we next consider how the structure of this periodic steady state is affected by inclusion of an admissible kernel perturbation $\overline{\phi} \in \overline{K}(\mathbb{R})$, when $||\overline{\phi}||_1^m$ is small. We focus attention on the first tongue $\Omega_1$, with similar results, which are not repeated, for the  subsequent tongues $\Omega_i, i=2,3,..$ . Thus we fix $(\lambda,D)\in \Omega_1$ and explore the existence and structure of an even, positive, periodic steady state when $||\overline{\phi}||_1^m$ is small and positive, and, formally $D=O(1)$ as $||\overline{\phi}||_1^m\to 0$. Via developing a natural perturbation theory, formally based upon small $||\overline{\phi}||_1^m$, with $D=O(1)$ as $||\overline{\phi}||_1^m\to 0$, we demonstrate that, for a given  admissible kernel perturbation with $||\overline{\phi}||_1^m$ small, then for each $(\lambda,D) \in \Omega_1$ satisfying $D=O(1)$ as $||\overline{\phi}||_1^m\to 0$, there exists a unique, even and positive periodic steady state associated with (IBVP)$_p$.  Moreover, this periodic steady state is a regular perturbation of the corresponding periodic steady state associated with (IBVP). The magnitude of the perturbation is determined to be of $O(D^{-1}||\overline{\phi}||_1^m))$ with $D=O(1)$ as $||\overline{\phi}||_1^m\to 0$, uniformly for $(x,\lambda) \in \mathbb{R} \times (1/2,1)$. This estimate of the perturbation magnitude determines that the regular perturbation structure fails when $D$ becomes small in $\Omega_1$, and in particular when $D=O(||\overline{\phi}||_1^m)$ as $||\overline{\phi}||_1^m \to 0$.

Naturally, we next examine points in $\Omega_1$ which have $D=O(||\overline{\phi}||_1^m)$ as $||\overline{\phi}||_1^m \to 0$. We develop an approximate approach ultimately based on the classical Schauder Fixed Point Theorem on a suitably defined compact subset of the real sequence space $l_1$. With a fixed admissible kernel perturbation which has $||\overline{\phi}||_1^m$ small, this first reveals that at each such point in $\Omega_1$, there exits at least one singularly (has O(1) perturbation as $||\overline{\phi}||_1^m\to0$) perturbed positive periodic state. However now the magnitude of the perturbation is no longer small, but of $O(1)$ as $||\overline{\phi}||_1^m \to 0$. We anticipate that there will be a principal branch, which is now a singularly perturbed continuation of the unique regularly perturbed branch, as $D$ decreases in $\Omega_1$ from being $O(1)$ to being $O(||\overline{\phi}||_1^m$. In addition there is the possibility of further secondary bifurcations as $D$ decreases further, creating multiple branches of singularly perturbed positive periodic states. To investigate this in more detail, we next focus attention on two specific, but representative, kernel perturbations, both of which retain their support on $[-1/2,1/2]$. The first (which we refer to as the positive kernel perturbation) has focus near the origin, whilst, conversely, the second (which we refer to as the negative kernel perturbation) has focus towards the ends of the support interval. This detailed investigation reveals that according to the approximate theory, in each case, the principal branch continues throughout this boundary layer region of $\Omega_1$, without identifying any secondary bifurcations. However, it does identify a 'spike-formation' process in both cases, and a 'peak-splitting' process which is solely restricted to the positive kernel perturbation case. The 'peak-splitting' process turns out to be a significant indicator - it leads to a local break-down in the approximate theory developed this far, and the nature of this break-down is investigated in detail. The break-down is shown to be caused by an incipient structure change in the positive periodic state on the principle branch, and this drives a secondary bifurcation from the principle branch. To study the secondary bifurcations and beyond, we must abandon the approximate approach, and return to the full problem determining positive periodic steady states.

$\mathbf{In~Section~4}$ the secondary bifurcations, and beyond, are investigated in the remaining sublayers of the boundary layer in $\Omega_1$ via careful numerical investigation (which requires a significant degree of sophistication and care to perform with confidence and accuracy). These two complementary investigations enable us to determine the complete outcomes in both cases of the positive and negative kernel perturbations which, in turn, enables us to draw some general conclusions and conjectures for the cases of more general admissible kernel perturbations which also have the 'positive' or 'negative' type of structure.

$\mathbf{In~Section~5}$ we return to the evolution problem (IBVP)$_p$ with the aim of using the theory developed in sections 2-4 to investigate its large-t structure. We begin performing some direct numerical solutions of (IBVP)$_p$. In fact, we observe that a number of key mechanisms in the spatio-temporal evolution are preserved from the unperturbed evolution problem (IBVP), even in the limit $D=O(||\overline{\phi}||_1^m)$ as $||\overline{\phi}||_1^m \to 0$. 

\section{Stability of Equilibrium States} \label{sec2}
The nonlocal PDE (\ref{eqn1.1}) again has exactly two equilibrium states, namely the \emph{unreacted} equilibrium state $u=0$ and the \emph{fully reacted} equilibrium state $u=1$. As seen in (NB) for (IBVP), these equilibrium states, and in particular, their temporal stability characteristics, will play a key role in the large-$t$ development of the solution to (IBVP)$_p$. The local temporal stability of these equilibrium states is investigated via linearised theory, and follows that in (NB) (see section 3). In particular, as the linearisation of (\ref{eqn1.1}) about the equilibrium state $u=0$ does not involve the nonlocal terms, then the linearised problem (LIVP)$_0$, its analysis, and the conclusions thereof, are identical to those obtained in (NB) for this equilibrium state; we conclude that, for all admissible kernel perturbations $\overline{\phi}\in \overline{K}(\mathbb{R})$, the equilibrium state $u=0$ is always temporally unstable, and small disturbances evolve precisely as detailed in (NB). However, for the equilibrium state $u=1$, the associated linearised problem (LIVP)$_1$ now has the modified linear evolution equation,
\begin{equation}
    \overline{u}_t = D\overline{u}_{xx} - \int_{x- \frac{1}{2}}^{x+\frac{1}{2}}{\overline{u}(y,t)}dy - \int_{-\infty}^{\infty}{\overline{\phi}(x-y)\overline{u}(y,t)}dy~~\forall~~(x,t)\in D_{\infty}. \label{eqn2.1}
\end{equation}
with the formulation and notation as in (NB). We now consider the dispersion relation  of the perturbed linearised equation (\ref{eqn2.1}). Without repeating details, we obtain the dispersion relation as
\begin{equation}
    \omega(k,D) = \omega_T(k,D) + \hat{\phi}(k)~~\forall~~k\in \overline{\mathbb{R}}^+, \label{eqn2.2}
\end{equation}
where
\begin{equation}
    \omega_T(k,D) = Dk^2 + \frac{2}{k}\text{sin}\frac{1}{2}k~~\forall~~k\in \overline{\mathbb{R}}^+ \label{eqn2.3}
\end{equation}
is the corresponding dispersion relation for the top hat kernel, as in (NB), and $\hat{\phi}:\overline{\mathbb{R}}^+\to \mathbb{R}$ is given by
\begin{equation}
    \hat{\phi}(k) = 2\int_0^{\infty}{\overline{\phi}(s)\text{cos}(ks)}ds~~\forall~~k\in \overline{\mathbb{R}}^+. \label{eqn2.4} 
\end{equation}
With $\overline{\phi}\in \overline{K}(\mathbb{R})$, it is readily established that
\begin{equation}
\hat{\phi}\in C^R(\overline{\mathbb{R}}^+), \label{eqn2.5}
\end{equation}
\begin{equation}
    \hat{\phi}(k) \to 0 ~~\text{as}~~k\to \infty, \label{eqn2.6}
\end{equation}
\begin{equation}
    ||\hat{\phi}^{(r)}||_{\infty}\le ||(\cdot)^r \overline{\phi}(\cdot)||_1~~\text{for}~~r = 0,1,2,..,R \label{eqn2.7}
\end{equation}
\begin{equation}
    \hat{\phi}(0)=0.  \label{eqn2.8}
\end{equation}
It follows from (\ref{eqn2.2})-(\ref{eqn2.8}) that,
\begin{equation}
||(\omega,\omega') - (\omega_T,\omega_{T}')||_{\infty} \le \sqrt{2} \, \text{max} \{||\overline{\phi}||_1, ||(\cdot)\overline{\phi}(\cdot)||_1\},\label{eqn2.9}
\end{equation}
with $'$ representing differentiation with respect to $k$ and the supremum in $||\cdot||_{\infty}$ taken over $(k,D)\in \overline{\mathbb{R}}^+ \times \mathbb{R}^+$. Thus the perturbed dispersion relation $\omega(k,D)$, and its $k$-derivative $\omega'(k,D)$, are respectively \emph{uniformly close}
to the top hat dispersion relation $\omega_T(k,D)$, and its $k$-derivative $\omega_T'(k,D)$, for all $(k,D)\in \overline{\mathbb{R}}^+\times 
\mathbb{R}^+$, as $\text{max}\{||\overline{\phi}||_1, ||(\cdot)\overline{\phi}(\cdot)||_1\}\to 0$. As an immediate consequence we may infer that the linearised temporal stability characteristics of the equilibrium state $u=1$ are qualitatively unaffected, and only uniformly regularly perturbed quantitatively, in the first quadrant of the $(k,D)$-plane, when $\text{max}\{||\overline{\phi}||_1, ||(\cdot)\overline{\phi}(\cdot)||_1\}$ is sufficiently small. We may also conclude that, specifically, in the positive quadrant of the $(k,D)$-plane the level curves of $\omega(k,D)$ and $\omega_T(k,D)$ (and hence the associated neutral curves) are uniformly close in both location and slope as $||\overline{\phi}||_1^m\equiv \text{max}\{||\overline{\phi}||_1, ||(\cdot)\overline{\phi}(\cdot)||_1\}\to 0$, with a uniform displacement of $O(||\overline{\phi}||_1^m)$.\\\\
 The key conclusions from this analysis can usefully and formally be summarised in:
\begin{theorem}
Let $\overline{\phi}\in \overline{K}(\mathbb{R})$ be an admissible kernel perturbation with $||\overline{\phi}||_1^m$ small. Then the linearised stability properties for the unreacted  equilibrium state $u=0$ associated with (IBVP)$_p$ and (IBVP) are identical. For the linearised stability properties of the fully reacted equilibrium state $u=1$, the dispersion relation (as a smooth surface above the positive quadrant of the $(k,D)$-plane) and the neutral curve (as a smooth curve on the positive quadrant of the $(k,D)$-plane) associated with (IBVP)$_p$ are both regular perturbations (for both location and gradient) in $||\overline{\phi}||_1^m$ from the corresponding unperturbed forms associated with (IBVP), uniform over the positive quadrant of the $(k,D)$-plane. Thus the principal linearised stability characteristics of both equilibrium states are unaltered by inclusion of the admissible kernel perurbation when $||\overline{\phi}||_1^m$ is sufficiently small.
\end{theorem}

In relation to (IBVP)$_p$, the above analysis determines that the conclusions laid out in the final part of section 3 in (NB) for (IBVP), and in particular the conjectures [P1] and [P2], continue to hold for (IBVP)$_p$ when $||\overline{\phi}||_1^m$ is sufficiently small. Specifically, the critical value of $D$ given by $\Delta_1$ is now replaced by $\Delta_1^p = \Delta_1 + O(||\overline{\phi}||_1^m)$ as $||\overline{\phi}||_1^m\to 0$. The next natural step in investigating [P1] and [P2] is to examine how the existence of, bifurcation to and structural features of positive periodic states of the nonlocal Fisher-KPP equation with top hat kernel (as studied in detail in section 4 of (NB)) are first qualitatively, and then quantitatively, affected when considering the nonlocal Fisher-KPP equation with admissible perturbed kernel, as given in (\ref{eqn1.1}), having $\overline{\phi}\in \overline{K}(\mathbb{R})$, with $||\overline{\phi}||_1^m$ small. 

\section{Positive Periodic Steady States} \label{sec3}
As in section 4 of (NB), and in relation to [P2], we seek to again identify steady state bifurcations from the equilibrium state $u=1$ which give rise to periodic steady states. Without repeating the detailed calculations, we find that, in the positive quadrant of the $(\lambda,D)$-plane, with $\lambda$ again representing fundamental wavelength, there remains the bifurcation locus across which there is a pitchfork bifurcation to periodic steady states from the equilibrium state $u=1$. The bifurcation locus is now given as the intersection of the curve
\begin{equation}
    D = D_p(\lambda) \equiv  D_T(\lambda) - \frac{\lambda^2}{4\pi^2}\hat{\phi}\left(\frac{2\pi}{\lambda}\right) \label{eqn3.1}
\end{equation}
with the positive quadrant of the $(\lambda,D)$-plane. Here
\begin{equation}
    D_T(\lambda) = -\frac{\lambda^3}{4\pi^3}\text{sin}\left(\frac{\pi}{\lambda}\right) \label{eqn3.2} 
\end{equation}
is the corresponding curve with the unperturbed top hat kernel, as obtained in (NB). It follows from (\ref{eqn2.5})-(\ref{eqn2.8}), (\ref{eqn3.1}) and (\ref{eqn3.2}) that,
\begin{equation}
    D_p(\lambda) = D_T(\lambda) + O(\lambda^2||\overline{\phi}||_1^m) ~~\text{as}~~||\overline{\phi}||_1^m\to 0 \label{eqn3.3}
\end{equation}
uniformly for $\lambda >0$, whilst
\begin{equation}
    D'_p(\lambda) = D'_T(\lambda)  + O((\lambda + \pi)||\overline{\phi}||_1^m) ~~\text{as}~~||\overline{\phi}||_1^m\to 0 \label{eqn3.4} 
\end{equation}
uniformly for $\lambda>0$. The bifurcation locus is contained in the rectangle $[0,\lambda^p] \times [0,\Delta_1^p]$, with $\Delta_1^p = \Delta_1 + O(||\overline{\phi}||_1^m)$ as $||\overline{\phi}||_1^m\to 0$, and $\lambda^p$ being the largest positive zero of $D_p(\lambda)$, so that $\lambda^p = 1 + O(||\overline{\phi}||_1^m)$. It follows from (\ref{eqn3.3}) and (\ref{eqn3.4}) that the bifurcation locus, in the first quadrant of the $(\lambda,D)$-plane, when the admissible kernel perturbation $\overline{\phi}\in \overline{K}(\mathbb{R})$ has $||\overline{\phi}||_1^m$ small, will be uniformly close, of $O(||\overline{\phi}||_1^m)$, in both position and slope, to that when the kernel perturbation is absent, \emph{except in a neighbourhood} $\mathcal{N}(||\overline{\phi}||_1^m)$ of $(\lambda,D)=(0,0)$ in which $(\lambda,D) = (O(||\overline{\phi}||_1^m) , O((||\overline{\phi}||_1^m)^3)$. Thus, as in (NB), and except possibly in the small neighbourhood $\mathcal{N}((||\overline{\phi}||_1^m)$ of the origin, the bifurcation locus is again in the form of a sequence of 'tongue like' curves based on the $\lambda$-axis, with decreasing height, and with each 'tongue' located within an $O(||\overline{\phi}||_1^m)$ neighbourhood of the corresponding tongue, $\partial\Omega_i (i=1,2,..)$, for the case of the unperturbed kernel top hat kernel. For simplicity, we will adopt directly, and without adaptation, the notion for, and relating to, these 'tongues` introduced in (NB), always bearing in mind, from above, that associated locations in the $(\lambda,D)$-plane are perturbed by $O(||\overline{\phi}||_1^m)$ when the kernel is perturbed with $||\overline{\phi}||_1^m$ small. A straightforward, but lengthy, calculation also establishes  that, at each bifurcation point outside of $\mathcal{N}(||\overline{\phi}||_1^m)$, the nature and criticality of each bifurcation is unchanged when $||\overline{\phi}||_1^m$ is small, with each bifurcation being a pitchfork bifurcation, generating a unique (up to translation in $x$) periodic steady state with small amplitude relative to the equilibrium state $u=1$, on traversing through the bifurcation point in the $(\lambda,D)$-plane, from the \emph{exterior to the interior} of the tongue. As in (NB), translation invariance can be most conveniently fixed by selecting that representative which is an even function of $x$.
With the unperturbed top hat kernel, we have demonstrated in (NB) that at each point, within each tongue, there is a unique, symmetric and positive periodic steady state, which we labelled as $u = F_p(x,\lambda,D)$ for $x\in \mathbb{R}$ and $(\lambda,D)\in \Omega$. At each such point $(\lambda,D)\in \Omega \setminus{\mathcal{N}(||\overline{\phi}||_1^m)}$, we next consider how the structure of this periodic steady state is affected by inclusion of an admissible kernel perturbation $\overline{\phi} \in \overline{K}(\mathbb{R})$ when $||\overline{\phi}||_1^m$ is small. We focus attention on the first tongue $\Omega_1$, with similar results, which are not repeated, for the  subsequent tongues $\Omega_i, i=2,3,..$ . First we fix $(\lambda,D)\in \Omega_1$ and explore the existence and structure of an even, positive, periodic steady state when, formally,  $D=O(1)$ as $||\overline{\phi}||_1^m\to 0$.

\subsection{$(\lambda,D)\in \Omega_1$ with $D=O(1)$ as $||\overline{\phi}||_1^m\to 0$}
We fix $(\lambda,D)\in \Omega_1$ and look for an even, positive periodic steady state, with fundamental period $\lambda$, of equation (\ref{eqn1.1}) with nontrivial admissible kernel perturbation $\overline{\phi}\in \overline{K}(\mathbb{R})$, in the limit $||\overline{\phi}||_1^m\to 0$ with $D=O(1)$. An examination of the terms in equation (\ref{eqn1.1}) indicates that this limit is a regular perturbation of the unperturbed top hat kernel structure. To investigate this, we look for an even, positive periodic steady state, with fundamental period $\lambda$, in the form,
\begin{equation}
    \mathcal{F}_p(x,\lambda,D;\overline{\phi}) = F_p(x,\lambda,D) + ||\overline{\phi}||_1^m\overline{F}_p(x,\lambda,D;\overline{\phi}_N) + o(||\overline{\phi}||_1^m)~~\text{as}~~||\overline{\phi}||_1^m\to 0  \label{eqn3.5}
\end{equation}
with $x\in \mathbb{R}$. Here $\overline{\phi}_N = (||\overline{\phi}||_1^m)^{-1}\overline{\phi}$ is the normalised kernel perturbation in $||\cdot||_1^m$. On substitution from (\ref{eqn3.5}) into equation (\ref{eqn1.1}), we obtain the following problem for $\overline{F}_p$, namely,
\begin{equation}
   L[\overline{F}_p]\equiv \overline{F}_p'' + D^{-1}(1 - J(F_p))\overline{F}_p - D^{-1}F_pJ(\overline{F}_p) = D^{-1}F_p\overline{J}(F_p),~~x\in \mathbb{R}, \label{eqn3.6}
\end{equation}
subject to the conditions
\begin{equation}
\overline{F}_p \in P_{\lambda}(\mathbb{R})\cap C^2(\mathbb{R}),~~\overline{F}_p~~\text{is even}, \label{eqn3.7}   
\end{equation}
 and where $P_{\lambda}(\mathbb{R})$ represents the set of $\lambda$-periodic real-valued functions on $\mathbb{R}$. Here, for any function $h\in C(\mathbb{R})\cap L^{\infty}(\mathbb{R})$,
\begin{equation}
    J(h)(x) = \int_{x-\frac{1}{2}}^{x+\frac{1}{2}}{h(s)}ds, \label{eqn3.8}
\end{equation}
and
\begin{equation}
    \overline{J}(h)(x) = \int_{-\infty}^{\infty}{h(s)\overline{\phi}_N(x-s)}ds, \label{eqn3.9}
\end{equation}
for all $x\in \mathbb{R}$. We observe that, since $F_p\in P_{\lambda}(\mathbb{R})\cap C^2(\mathbb{R})$ and even, then,
\begin{equation}
    J(F_p), \overline{J}(F_p)\in P_{\lambda}(\mathbb{R})\cap C^2(\mathbb{R})~~\text{with both even}.  \label{eqn3.10}
\end{equation}
In order to analyse the inhomogeneous, $\lambda$-periodic, linear, boundary value problem for $L$, given by (\ref{eqn3.6}) and (\ref{eqn3.7}), we must digress into the spectral properties of $L$. With the linear operator $L:P_{\lambda}(\mathbb{R})\cap C^2(\mathbb{R})\to P_{\lambda}(\mathbb{R})\cap C(\mathbb{R})$ given by
\begin{equation}
    L[\psi] = \psi'' + D^{-1}(1 - J(F_p))\psi - D^{-1}F_pJ(\psi) ~~\forall~~ \psi\in P_{\lambda}(\mathbb{R})\cap C^2(\mathbb{R}) , \label{eqn3.11} 
\end{equation}
it is straightforward to determine the associated adjoint operator $L^*:P_{\lambda}(\mathbb{R})\cap C^2(\mathbb{R})\to P_{\lambda}(\mathbb{R})\cap C(\mathbb{R})$, as
\begin{equation}
    L^* [\psi] = \psi'' + D^{-1}(1 - J(F_p))\psi - D^{-1}J(F_p\psi) ~~\forall~~ \psi\in P_{\lambda}(\mathbb{R})\cap C^2(\mathbb{R}). \label{eqn3.12}
\end{equation}
We conclude that $L$ is not self-adjoint. However, we may decompose $L$ by writing $L\equiv L_S + L_N$, where $L_S, L_N:P_{\lambda}(\mathbb{R})\cap C^2(\mathbb{R})\to P_{\lambda}(\mathbb{R})\cap C(\mathbb{R})$ are the self-adjoint differential operator given by
\begin{equation}
    L_S[\psi] = \psi'' + D^{-1}(1 - J(F_p))\psi ~~\forall~~ \psi\in P_{\lambda}(\mathbb{R})\cap C^2(\mathbb{R}) \label{eqn3.13}
\end{equation}
and the non-self-adjoint integral operator given by
\begin{equation}
    L_N[\psi] = - D^{-1}F_pJ(\psi) ~~\forall~~ \psi\in P_{\lambda}(\mathbb{R})\cap C^2(\mathbb{R}). \label{eqn3.14}
\end{equation}
We now observe that the integral operator $L_N$ is compact and, for any regular point $\sigma \in \mathbb{C}$ of the differential operator $L_S$ (of which there is at least one), that the integral operator $L_N(L_S + \sigma I)^{-1}$ is also compact. Thus we may conclude that $L_N$ is $L_S$-completely continuous. In addition, the classical periodic Sturm-Liouville theory immediately determines that $L_S$ has a discrete and countably infinite spectrum, with each eigenvalue being real and isolated in the complex-plane. Under these conditions we may conclude that $L$ falls into the class of \emph{weakly perturbed self-adjoint operators}. We can then follow the general spectral theory developed by Gohberg and Krein in \cite{selfadjoint} to conclude that the eigenvalues of $L$ form a countably infinite sequence of isolated points in the complex-plane, say $\mu_r,~r=0,1,2,..$, with $|\mu_r|\to \infty$ and arg$(\mu_r)$ approaching $0$ or $\pi$ as $r\to \infty$. The eigenvalue $\mu_0$ can be chosen so that its eigenspace is span$\{e_0\}$ whilst, for $r=1,2,3...$, the eigenspace of $\mu_r$ is span$\{e_r^a,e_r^s\}$. Here the functions $e_0, e_r^s, e_r^a\in P_{\lambda}(\mathbb{R})\cap C^2(\mathbb{R})$ are possibly complex-valued, and normalised in the $\lambda$-periodic $L^2$ inner product. In the specific case at hand, we can verify directly that 
\begin{equation}
    \mu_0=0~~\text{and}~~e_0(x)=F_p'(x)~~\forall~~x\in \mathbb{R}.\label{eqn3.15}
\end{equation}
which we observe is an odd function of $x\in \mathbb{R}$. In addition, the invariance of $L$ to the transformation $x\to (-x)$ and to conjugation, determines that we may take, whenever Im$(\mu_j) \neq 0$, one or the other of
\begin{equation}
  \mu_{j\pm1}= \overline{\mu}_{j},\label{eqn3.16}
\end{equation}
and, for each $r=1,2,3,..$, choose $e_r^a$ and $e_r^s$ to be odd and even functions of $x\in \mathbb{R}$, respectively. The adjoint operator $L^*$ has corresponding spectrum with eigenvalues $\overline{\mu}_r,~r=0,1,2,..$. The corresponding eigenfunctions $f_0$ and $f_r^a, f_r^s, r=1,2,3,...$ share the same properties as the corresponding eigenfunctions for $L$. Moreover, in the $\lambda$-periodic $L^2$ inner product, 
\begin{equation}
    <e_0,f_0> = 1,~~<e_0,f_r^a>=<e_0,f_r^s>=<e_r^a,f_0>=<e_r^s,f_0>=0, \label{eqn3.17}
\end{equation}
for $r=1,2,3...$, whilst,
\begin{equation}
    <e_i^a,f_j^s>=<e_i^s,f_j^a>=0~~\forall~~ i,j=1,2,3.., \label{eqn3.18}
\end{equation}
and,
\begin{equation}
    <e_i^a,f_j^a>=<e_i^s,f_j^s>=\delta_{ij}~~\forall~~ i,j=1,2,3.., \label{eqn3.19}
\end{equation}
with $\{e_0,e_1^a,e_1^s,e_2^a,e_2^s,....\}$ forming a basis for the $\lambda$-periodic $L^2$ inner product space. In addition, for the specific operator $L$ under consideration here, we can employ a straightforward WKB asymptotic analysis to establish that the large eigenvalues are all real and have
\begin{equation}
\mu_n = \frac{4\pi^2n^2}{\lambda^2} + O(1)~~\text{as}~~n\to \infty, \label{eqn3.20}
\end{equation}
with,
\begin{equation}
 e_n^a(x),f_n^a(x) \sim \frac{\sqrt{2}}{\sqrt{\lambda}} \text{sin}\left(\frac{2n\pi x}{\lambda}\right)
 \end{equation}

\begin{equation}
 e_n^s(x),f_n^s(x) \sim \frac{\sqrt{2}}{\sqrt{\lambda}} \text{cos}\left(\frac{2n\pi x}{\lambda}\right)   
\end{equation}
as $n\to \infty$ uniformly for $x\in \mathbb{R}$.\\
We can now return to the inhomogeneous, $\lambda$-periodic, linear, boundary value problem for $L$ given by (\ref{eqn3.6}) and (\ref{eqn3.7}), and which we henceforth refer to as [NSP].  First we observe, via (\ref{eqn3.10}), that $D^{-1}F_p\overline{J}(F_p) \in P_{\lambda}(\mathbb{R})\cap C^2(\mathbb{R})$, which therefore may be written (uniquely) as the eigenfunction expansion,
\begin{equation}
 D^{-1}F_p(x,\lambda,D)\overline{J}(F_p)(x,\lambda,D) = a_0e_0(x) + \sum_{r=1}^{\infty}(a_re_r^a(x) + b_re_r^s(x))~~\forall~~x\in \mathbb{R}, \label{eqn3.24}   
\end{equation}
with the coefficients given by,
\begin{equation}
    a_0 = <D^{-1}F_p\overline{J}(F_p),f_0> = 0,~~a_r = <D^{-1}F_p\overline{J}(F_p),f_r^a> = 0  \label{eqn3.25}
\end{equation}
for $r=1,2,3...$, since $D^{-1}F_p\overline{J}(F_p)$ is an even function of $x\in \mathbb{R}$, whilst,
\begin{equation}
 b_r =<D^{-1}F_p\overline{J}(F_p),f_r^s>.  \label{eqn3.26}   
\end{equation}
We should observe here that all  eigenvalues, eigenfunctions and coefficients introduced above have a continuous dependence on the parameter $D$, \emph{when $D$ is strictly positive}, although this is suppressed in some of the notation for convenience. Now, via conditions (\ref{eqn3.7}), \emph{any solution} to [NSP] may be written as the eigenfunction expansion,
\begin{equation}
    \overline{F}_p(x,\lambda,D;\overline{\phi}_N) = \sum_{r=1}^{\infty}{c_re_r^s(x)}~~\forall~~x\in \mathbb{R}, \label{eqn3.27}
\end{equation}
with the coefficients $c_r, r=1,2,...$ to be determined. On substitution from (\ref{eqn3.24}) to (\ref{eqn3.26}) into [NSP] it is readily established that (\ref{eqn3.27}) provides a solution if and only if the coefficients are chosen as (recalling that $\mu_r \neq 0~~\forall~~r=1,2,..$)
\begin{equation}
c_r = - \frac{b_r}{\mu_r}\label{eqn3.28}
\end{equation}
for all $r=1$, $2$, $3\ldots$. We may conclude that [NSP] has a unique solution, and specifically, that solution is given by (\ref{eqn3.27}) with (\ref{eqn3.28}).

This formally confirms that when $(\lambda,D)\in \Omega_1$ then there is a unique even, positive, $\lambda$-periodic steady state when $D=O(1)$ as $||\overline{\phi}||_1^m \to 0$, and this is a regular perturbation of the associated even, positive, $\lambda$-periodic steady state for the unperturbed top hat kernel. However, an examination of equation (\ref{eqn1.1}) immediately suggests that this regular perturbation may fail when $D$ becomes small, with formally $D=o(1)$ as $||\overline{\phi}||_1^m\to 0$. In particular, a balancing of terms in equation (\ref{eqn1.1}) suggests that the magnitude of the perturbation to the $\lambda$-periodic steady state for the top hat kernel increases rapidly to $O(1)$ when $D$ becomes small enough to satisfy $D=O(||\overline{\phi}||_1^m)$ as $D\to 0$. A more precise evaluation of this observation can be made by examining the detailed structure of $\overline{F}_p(x,\lambda,D;\overline{\phi}_N)$ when $D$ is small. This could be achieved by examining and approximating all of the components of the above spectral theory when $D$ is small, and then putting that together to obtain the approximate form of the solution to [NSP] for $D$ small. However, we find that constructing the (unique) solution to [NSP] directly from (\ref{eqn3.6}) and (\ref{eqn3.7}), in the limit $D\to 0$, is most convenient, and we turn to this now.

To start, we recall the approximation to $F_p(x,\lambda,D)$ as $D\to 0$, which was constructed in detail in subsection 4.1 of (NB). This determines that we have, for each wavelength $\lambda\in \left(\frac{1}{2}+O\left(\sqrt{D}\right),1-O\left(D\right)\right)$,
\begin{equation}
F_p(x,\lambda,D)\sim 
\begin{cases}
\dfrac{\pi}{4a(\lambda)}\cos \left(\dfrac{\pi x}{2a(\lambda)}\right),~~0\le|x|<a(\lambda)\\
E(x,\lambda,D),~~a(\lambda)< |x| \le \frac{1}{2}\lambda, 
\end{cases} \label{eqn3.29}
\end{equation}
as $D\to 0$ uniformly over one spatial wavelength $\lambda$, with $E(x,\lambda,D)$, as given in (NB), being exponentially small in $D$ as $D\to 0$, and
\begin{equation}
a(\lambda) = \frac{1}{2}\left(\lambda - \frac{1}{2}\right) \label{eqn3.30}
\end{equation}
whilst
\begin{equation}
    \int_{-a(\lambda)}^{a(\lambda)}{F_p(y,\lambda,D)}dy = 1 - \frac{\pi^2}{4a(\lambda)^2}D + o(D^{\frac{5}{4}})   \label{eqn3.31}
\end{equation}
as $D\to 0$. It should also be noted, from (NB), that the gradient and curvature discontinuities at $|x|=a(\lambda)$ are smoothed out across passive thin edge regions of spatial thickness $O(D^{\frac{1}{4}})$. Without loss of generality, we restrict attention to equation (\ref{eqn3.6}) on the period $x\in (-\frac{1}{2}\lambda, \frac{1}{2}\lambda)$ and a balancing of terms in equation (\ref{eqn3.6}), with use of (\ref{eqn3.29})-(\ref{eqn3.31}), indicates that the solution to [NSP] should be developed in the form,
\begin{equation}
    \overline{F}_p(x,\lambda,D;\overline{\phi}_N) \sim 
    \begin{cases}
    D^{-1}\overline{F}(x,\lambda;\overline{\phi}_N),~~0\le|x|<a(\lambda)\\
D^{-1}\overline{E}(x,\lambda,D;\overline{\phi}_N),~~a(\lambda)< |x| \le \frac{1}{2}\lambda, 
\end{cases} \label{eqn3.32}
\end{equation}
as $D\to 0$ uniformly for $x\in (-\frac{1}{2}\lambda, \frac{1}{2}\lambda)$, with $\overline{E}(x,\lambda,D;\overline{\phi}_N)$ being exponentially small in $D$ as $D\to 0$. The objective now is to determine $\overline{F}(\cdot,\lambda;\overline{\phi}_N):[-a(\lambda),a(\lambda)]\to \mathbb{R}$. We first use (\ref{eqn3.30}) and (\ref{eqn3.31}) to approximate the periodic coefficients and inhomogeneous term in (\ref{eqn3.6}) as,
\begin{equation}
    D^{-1}(1 - J(F_p))(x,\lambda,D) \sim \frac{\pi^2}{4a(\lambda)^2}, \label{eqn3.33}
\end{equation}
\begin{equation}
    D^{-1}(F_p\overline{J}(F_p))(x,\lambda,D) \sim \frac{\pi^2}{16a(\lambda)^2D}\overline{I}(x,\lambda;\overline{\phi}_N)\cos \left(\dfrac{\pi x}{2a(\lambda)}\right), \label{eqn3.34}
\end{equation}
as $D\to 0$ uniformly for $x\in [-a(\lambda),a(\lambda)]$. Here $\overline{I}(\cdot,\lambda;\overline{\phi}_N)\in C^2([-a(\lambda),a(\lambda)]$ is given by,
\begin{equation}
\overline{I}(x,\lambda;\overline{\phi}_N) = \sum_{n=-\infty}^{\infty}\left(\int_{-a(\lambda)}^{a(\lambda)}{\overline{\phi}_N(x-s-n\lambda)\cos \left(\dfrac{\pi s}{2a(\lambda)}\right)}ds\right) \label{eqn3.35}
\end{equation}
for $x\in [-a(\lambda),a(\lambda)]$, and is an even function of $x$. On substitution from equations (\ref{eqn3.29})-(\ref{eqn3.34}) into equation (\ref{eqn1.1}), we obtain at leading order with $D$ small,
\begin{equation}
    \overline{F}'' + \frac{\pi^2}{4a(\lambda)^2}\overline{F} = H(x,\lambda;\overline{\phi}_N),~~x\in (-a(\lambda),a(\lambda)),\label{eqn3.36}  
\end{equation}
with,
\begin{equation}
  H(x,\lambda;\overline{\phi}_N) = \frac{\pi}{4a(\lambda)}\overline{\alpha}(\lambda;\overline{\phi}_N) \cos \left(\dfrac{\pi x}{2a(\lambda)}\right) + \frac{\pi^2}{16a(\lambda)^2}\overline{I}(x,\lambda;\overline{\phi}_N)\cos \left(\dfrac{\pi x}{2a(\lambda)}\right), \label{eqn3.36'}  
\end{equation}
for $x\in [-a(\lambda),a(\lambda)]$, which is also an even function of $x$. Here we have introduced the constant $\overline{\alpha}(\lambda;\overline{\phi}_N)$, which is to be determined, so that
\begin{equation}
     \int_{-a(\lambda)}^{a(\lambda)}{\overline{F}_p(y,\lambda,D;\overline{\phi}_N)}dy \sim \overline{\alpha}(\lambda;\overline{\phi}_N) \label{eqn3.37}
\end{equation}
as $D\to 0$. The boundary conditions are
\begin{equation}
    \overline{F}(-a(\lambda),\lambda;\overline{\phi}_N) = \overline{F}(a(\lambda),\lambda;\overline{\phi}_N) = 0,  \label{eqn3.38} 
\end{equation}
\begin{equation}
    \overline{F}(x,\lambda;\overline{\phi}_N)~~\text{is even in }~x,  \label{eqn3.39}
\end{equation}
whilst (\ref{eqn3.32}) and (\ref{eqn3.37}) require that
\begin{equation}
 \int_{-a(\lambda)}^{a(\lambda)}{\overline{F}(y,\lambda;\overline{\phi}_N)}dy = 0. \label{eqn3.40}   
\end{equation}
The homogeneous Dirichlet conditions (\ref{eqn3.38}) arise to enable asymptotic matching to the regions where $\overline{F}_p$ is exponentially small, via symmetric edge layers located at $x= \pm a(\lambda) + O(D^{\frac{1}{4}})$ (see (NB), section 4, for similar details). We will refer to this linear, inhomogeneous boundary value problem as [BVP]. The general solution to the linear, inhomogeneous ODE (\ref{eqn3.37}) is readily obtained, via variation of parameters, as,
\begin{multline}
    \overline{F}(x,\lambda;\overline{\phi}_N) = A\cos \left(\dfrac{\pi x}{2a(\lambda)}\right) + B\sin \left(\dfrac{\pi x}{2a(\lambda)}\right)\\
    + \frac{2a(\lambda)}{\pi}\sin \left(\dfrac{\pi x}{2a(\lambda)}\right) \int_{0}^{x}{H(y,\lambda;\overline{\phi}_N)\cos \left(\dfrac{\pi y}{2a(\lambda)}\right)}dy\\
    - \frac{2a(\lambda)}{\pi}\cos \left(\dfrac{\pi x}{2a(\lambda)}\right) \int_{0}^{x}{H(y,\lambda;\overline{\phi}_N)\sin \left(\dfrac{\pi y}{2a(\lambda)}\right)}dy  \label{eqn3.41} 
\end{multline}
for $x\in [-a(\lambda),a(\lambda)]$, with $A$ and $B$ arbitrary real constants. We first apply condition (\ref{eqn3.39}) which requires setting $B=0$. After this, conditions (\ref{eqn3.38}) require that
\begin{equation}
  \int_{0}^{a(\lambda)}{H(y,\lambda;\overline{\phi}_N)\cos \left(\dfrac{\pi y}{2a(\lambda)}\right)}dy = 0,  
\end{equation}
which, on using (\ref{eqn3.36'}), determines
\begin{equation}
  \overline{\alpha}(\lambda;\overline{\phi}_N) = -\frac{\pi}{2a(\lambda)^2}\int_{0}^{a(\lambda)}{\overline{I}(y,\lambda;\overline{\phi}_N)\cos^{2} \left(\dfrac{\pi y}{2a(\lambda)}\right)}dy. \label{eqn3.43}
\end{equation}
Finally, condition (\ref{eqn3.40}) determines
\begin{equation}
    A = A(\lambda;\overline{\phi}_N) \equiv  \frac{\pi}{2a(\lambda)}\int_{0}^{a(\lambda)}{\Psi(y,\lambda;\overline{\phi}_N)}dy, \label{eqn3.44}
\end{equation}
where
\begin{multline}
   \Psi(x,\lambda;\overline{\phi}_N) = \frac{2a(\lambda)}{\pi}\sin \left(\dfrac{\pi x}{2a(\lambda)}\right) \int_{0}^{x}{H(y,\lambda;\overline{\phi}_N)\cos \left(\dfrac{\pi y}{2a(\lambda)}\right)}dy\\
    - \frac{2a(\lambda)}{\pi}\cos \left(\dfrac{\pi x}{2a(\lambda)}\right) \int_{0}^{x}{H(y,\lambda;\overline{\phi}_N)\sin \left(\dfrac{\pi y}{2a(\lambda)}\right)}dy  \label{eqn3.45}  
\end{multline}
for $x\in [-a(\lambda),a(\lambda)]$. The solution to [BVP] is now given, uniquely, by (\ref{eqn3.41}) with $B=0$, (\ref{eqn3.43}) and (\ref{eqn3.44}). We observe from (\ref{eqn3.35}) that,
\begin{equation}
    ||\overline{I}||_{\infty} \le ||\overline{\phi}_N||_1 ^m= 1
\end{equation}
and then from (\ref{eqn3.36'}) that,
\begin{equation}
  ||H||_{\infty} \le \frac{\pi}{4a(\lambda)} |\overline{\alpha}| + \frac{\pi^2}{16a(\lambda)^2}.  
\end{equation}
As a consequence, we have from (\ref{eqn3.41}) that,
\begin{equation}
    ||\overline{F}||_{\infty} \le |A| + \frac{4a(\lambda)^2}{\pi}||H||_{\infty}.
\end{equation}
It also follows from (\ref{eqn3.43}) and (\ref{eqn3.44}) that,
\begin{equation}
    |A| \le 2a(\lambda)^2||H||_{\infty},
\end{equation}
\begin{equation}
    |\overline{\alpha}| \le \frac{\pi}{2a(\lambda)}.
\end{equation}
Thus, combining the above inequalities, we have,
\begin{equation}
    ||\overline{F}(\cdot,\lambda;\overline{\phi}_N)||_{\infty} \le \frac{5}{8}\pi^2. \label{eqn3.50}
\end{equation}
with this upper bound independent of both $\lambda$ and $\overline{\phi}_N$. It follows from (\ref{eqn3.32}), with (\ref{eqn3.50}), that,
\begin{equation}
    \overline{F}_p(x,\lambda,D;\overline{\phi}_N) = O(D^{-1})~~\text{as}~~D\to 0, \label{eqn3.51}
\end{equation}
uniformly for $(x,\lambda,\overline{\phi}_N)\in \{(x,\lambda):x\in[-a(\lambda),a(\lambda)], \lambda\in \left(\frac{1}{2},1\right)\} \times \partial{B}_{\overline{K}}(1,0)$, with,
\begin{equation}
   \partial{B}_{\overline{K}}(1,0) = \{\overline{\phi}\in \overline{K}(\mathbb{R}): ||\overline{\phi}||_1^m=1\}.  
\end{equation}
We conclude that the regular perturbation expansion (\ref{eqn3.5}) for $\mathcal{F}_p(x,\lambda,D;\overline{\phi})$, with \\
$(x,\lambda,D,\overline{\phi})\in \mathbb{R}\times \Omega_1 \times \overline{K}(\mathbb{R})$, and based on $||\overline{\phi}||_1^m$ being small, becomes nonuniform when $D = O(||\overline{\phi}||_1^m)$ as $||\overline{\phi}||_1^m\to 0$, and when, in particular,
\begin{equation}
    \mathcal{F}_p(x,\lambda,D;\overline{\phi}) =
    \begin{cases}
O(1),~~0\le|x|<a(\lambda),\\
\mathcal{E}(x,\lambda,D;\overline{\phi}),~~a(\lambda)< |x| \le \frac{1}{2}\lambda, 
\end{cases} \label{eqn3.53}
\end{equation}
with $\mathcal{E}(x,\lambda, \overline{D}, \overline{\phi})$ being uniformly exponentially small as $||\overline{\phi}||_1^m\to 0$ with $D=O(||\overline{\phi}||_1^m)$. 

To bring together the above developments it is convenient to adopt (along with the notation introduced in (NB) and adapted above) the following notation relating to subsets of the positive quadrant of the $(\lambda,D)$-plane, namely,
with $||\overline{\phi}||_1^m$ small,
\begin{equation}
    \mathcal{B}(\overline{\phi}) \equiv \mathbb{R}^+ \times (0,\delta(||\overline{\phi}||_1^m))
\end{equation}
where $\delta(||\overline{\phi}||_1^m )= O(||\overline{\phi}||_1^m)$ as $||\overline{\phi}||_1^m\to 0$, and,
\begin{equation}
\Omega^-(\overline{\phi}) = \Omega_1 \cap \mathcal{B}(\overline{\phi}),~~
\Omega^+(\overline{\phi}) = \Omega_1 \setminus \Omega^-(\overline{\phi}),
\end{equation}
noting that $\Omega_1=\Omega^+(\overline{\phi}) \cup  \Omega^-(\overline{\phi})$. The results of this subsection can now be formally summarised in:
\begin{theorem}
   Let $\overline{\phi}\in \overline{K}(\mathbb{R})$ be an admissible kernel perturbation with $||\overline{\phi}||_1^m$ small. Then for each $(\lambda,D) \in \Omega^+(\overline{\phi})$ there exists a unique positive periodic steady state associated with (IBVP)$_p$, which is a regular perturbation of the corresponding unperturbed positive periodic steady state associated with (IBVP), with the magnitude of the perturbation spatially bounded in the supnorm, with bound of $O(D^{-1}||\overline{\phi}||_1^m)$ as $||\overline{\phi}||_1^m\to 0$, uniformly for $x\in [0,\lambda)$ and $(\lambda,D) \in \Omega^+(\overline{\phi})$. This regular perturbation structure fails when $D$ becomes sufficiently small so that it satisfies $D=O(||\overline{\phi}||_1^m)$, and $(\lambda,D)$ enters the boundary region $\Omega^-(\overline{\phi})$.
\end{theorem}

We now consider in detail the continuation of this family of periodic steady states into this singular boundary limit, when $D = O(||\overline{\phi}||_1^m)$ as $||\overline{\phi}||_1^m\to 0$.

\subsection{$(\lambda,D)\in \Omega_1$ with $D=O(||\overline{\phi}||_1^m)$ as $||\overline{\phi}||_1^m\to 0$}\label{subsec_smallD}
We first give a useful general formulation for considering this limit, and then address in detail two specific representative cases.
\subsubsection{General Formulation and Results}
To formalise this limit we write, for a chosen $\overline{\phi} \in \overline{K}(\mathbb{R})$,
\begin{equation}
    D = ||\overline{\phi}||_1^m \overline{D}
\end{equation}
with $\overline{D}=O(1)$ as $||\overline{\phi}||_1^m \to 0$. It then follows that, \emph{to continue this family of periodic steady states which have the structural form (\ref{eqn3.53}), from $D=O(1)$ (on $\Omega^+(\overline{\phi}) = \Omega_1 \setminus \Omega^-(\overline{\phi})$) into $D = O(||\overline{\phi}||_1^m)$ (on $\Omega^-(\overline{\phi}) = \Omega_1 \cap \mathcal{B}(\overline{\phi})$), within the tongue $\Omega_1$}, we should expand in the form,
\begin{equation}
    \mathcal{F}_p(x,\lambda,D;\overline{\phi}) = \mathcal{F}(x,\lambda,\overline{D}; \overline{\phi}_N) + o(1)~~\text{as}~~||\overline{\phi}||_1^m\to 0, \label{eqn3.55}
\end{equation}
with $\overline{D} = O(1)$ and $(x,\lambda,\overline{\phi}_N)\in \{(x,\lambda): \lambda\in \left(\frac{1}{2},1\right), x\in [-a(\lambda),a(\lambda)]\} \times \partial{B}_{\overline{K}}(1,0)$. Our objective is now to determine $\mathcal{F}(x,\lambda,\overline{D}; \overline{\phi}_N)$ for $x\in [-a(\lambda),a(\lambda)]$.  On substitution from (\ref{eqn3.55}) into (\ref{eqn1.1}) we obtain at leading order in $||\overline{\phi}||_1^m$ the nonlocal, nonlinear ODE,
\begin{equation}
    \mathcal{F}'' + \mathcal{F}\left(\mathcal{I} - \overline{D}^{-1}\overline{J}(\mathcal{F})\right) = 0,~~x\in (-a(\lambda),a(\lambda)),  \label{eqn3.56}
\end{equation}
with $\mathcal{I}$ being a constant (depending upon $\lambda$, $\overline{D}$ and $\overline{\phi}_N$) to be determined, so that,
\begin{equation}
   \int_{-a(\lambda)}^{a(\lambda)}{\mathcal{F}_p(y,\lambda,D;\overline{\phi})}dy = 1 - \mathcal{I}\overline{D}||\overline{\phi}||_1^m + o(||\overline{\phi}||_1^m)~~\text{as}~~||\overline{\phi}||_1^m\to 0. \label{eqn3.57} 
\end{equation}
The associated boundary conditions are,
\begin{equation}
    \mathcal{F}(-a(\lambda),\lambda,\overline{D};\overline{\phi}_N) = \mathcal{F}(a(\lambda),\lambda,\overline{D};\overline{\phi}_N) = 0,  \label{eqn3.58} 
\end{equation}
\begin{equation}
    \mathcal{F}(x,\lambda,\overline{D};\overline{\phi}_N)~~\text{is positive and even for }~x \in (-a(\lambda),a(\lambda)) ,  \label{eqn3.59}
\end{equation}
whilst (\ref{eqn3.55}) and (\ref{eqn3.57}) require that
\begin{equation}
 \int_{-a(\lambda)}^{a(\lambda)}{\mathcal{F}(y,\lambda,\overline{D};\overline{\phi}_N)}dy = 1. \label{eqn3.60}   
\end{equation}
The homogeneous Dirichlet conditions (\ref{eqn3.58}) again arise to enable asymptotic matching to the regions where $\mathcal{F}_p$ is exponentially small, via symmetric edge layers located at $x= \pm a(\lambda) + O(||\overline{\phi}||_1^m)^{\frac{1}{4}})$ (see (NB), section 4, for similar details). We will henceforth refer to the nonlinear, nonlocal boundary value problem (\ref{eqn3.56}) with (\ref{eqn3.58})-(\ref{eqn3.60}) as [NBVP]. 
\begin{remk}
At this stage, before we begin the detailed analysis of [NBVP], it is important to note that study of the asymptotic region where $\mathcal{F}_p(x,\lambda,D;\overline{\phi})$ is exponentially small as $||\overline{\phi}||_1^m \to 0$ (that is $a(\lambda)< |x| \le \frac{1}{2}\lambda$) is not required beforehand - however, it must be verified at the conclusion of the analysis of [NBVP] that this region can indeed be constructed, and asymptotically matched, via the aforementioned edge region, to the bulk region under consideration at present via [NBVP]. We return, and give consideration to, this observation at the end of the section.
\end{remk}
We begin our analysis of [NBVP] by considering the representation of the nonlocal term. Using the $\lambda$-periodicity of $\mathcal{F}_p$ and the structure given in (\ref{eqn3.53}) and (\ref{eqn3.55}) over one period, we may rewrite the nonlocal term in (\ref{eqn3.56}) as,
\begin{equation}
 \overline{J}(\mathcal{F})(x,\lambda,\overline{D};\overline{\phi}_N) = \sum_{r=-\infty}^{\infty}\left(\int_{-a(\lambda)}^{a(\lambda)}{\overline{\phi}_N(x-s-r\lambda)\mathcal{F}(s,\lambda,\overline{D};\overline{\phi}_N)}ds\right) \label{eqn3.61}
\end{equation}
for $x\in [-a(\lambda),a(\lambda)]$, and which we observe is an even function of $x$. When the normalised kernel perturbation decays sufficiently rapidly as $x\to \infty$ so that there are positive constants $M$ and $\gamma$ for which
\begin{equation}
    |\overline{\phi}_N(x)|,~|\overline{\phi}'_N(x)|,~|\overline{\phi}''_N(x)| \le \frac{M}{x^{1+\gamma}}~~\text{as}~~x\to \infty   \label{eqn3.62}
\end{equation}
then the series
\begin{equation}
    \Phi(X,\lambda;\overline{\phi}_N)\equiv  \sum_{r=-\infty}^{\infty}\overline{\phi}_N(X-r\lambda)  \label{eqn3.63}
\end{equation}
is absolutely and uniformly convergent for all $(X,\lambda)\in [-\frac{1}{2},\frac{1}{2}] \times (\frac{1}{2},1)$ (we also note that $\Phi(X,\lambda;\overline{\phi}_N)$ is even and, when extended onto the real line, periodic in $X$, with fundamental period $\lambda$). In addition, since $\overline{\phi}\in \overline{K}(\mathbb{R})$, it is readily observed that for $r\in \mathbb{Z} \setminus \{\pm1\}$ we have $\overline{\phi}_N(\cdot-r\lambda)\in PC^2([-\frac{1}{2},\frac{1}{2}]) \cap C([-\frac{1}{2},\frac{1}{2}])$ whilst $\overline{\phi}_N(\cdot\pm\lambda)\in PC^2([-\frac{1}{2},\frac{1}{2}] \cap C([-\frac{1}{2},\frac{1}{2}] \setminus \{\pm (\lambda - \frac{1}{2})\})$. It follows that 
\begin{equation}
    \Phi(\cdot,\lambda;\overline{\phi}_N)\in PC^2\left(\left[-\frac{1}{2},\frac{1}{2}\right]\right) \cap C\left(\left[-\frac{1}{2},\frac{1}{2}\right] \setminus \left\{\pm \left(\lambda - \frac{1}{2}\right)\right\}\right)   \label{eqn3.64}
\end{equation}
and we can write (\ref{eqn3.61}) in the form
\begin{equation}
  \overline{J}(\mathcal{F})(x,\lambda,\overline{D};\overline{\phi}_N) = \int_{-a(\lambda)}^{a(\lambda)}{\chi(x-s,\lambda;\overline{\phi}_N)\mathcal{F}(s,\lambda,\overline{D};\overline{\phi}_N)}ds~~\forall ~~x\in [-a(\lambda),a(\lambda)], \label{eqn3.65}
\end{equation}
where, for later convenience, we have introduced 
\[\chi(\cdot,\lambda;\overline{\phi}_N) \in PC^2\left(\left[-\frac{1}{2},\frac{1}{2}\right]\right) \cap C\left(\left[-\frac{1}{2},\frac{1}{2}\right] \setminus \left\{\pm \left(\lambda - \frac{1}{2}\right)\right\}\right),\] 
with,
\begin{equation}
    \chi(X,\lambda;\overline{\phi}_N) =
     \begin{cases}
\Phi(X,\lambda;\overline{\phi}_N),~~\text{supp}(\overline{\phi}_N) \nsubseteq [-1/2,1/2],\\
\overline{\phi}_N(X),~~\text{supp}(\overline{\phi}_N) \subseteq [-1/2,1/2], 
\end{cases} \label{eqn3.65'} 
\end{equation}
 We also observe from (\ref{eqn3.63}) that,
\begin{equation}
 \int_{-\frac{1}{2}\lambda}^{\frac{1}{2}\lambda}{\Phi(s,\lambda;\overline{\phi}_N)}ds = 0,  \label{eqn3.66}
\end{equation}
and
\begin{equation}
 \int_{-\frac{1}{2}\lambda}^{\frac{1}{2}\lambda}{|\Phi(s,\lambda;\overline{\phi}_N)|}ds \le 1.  \label{eqn3.67}
\end{equation}
 We can now use Fourier's Theorem, regularity and evenness to write, 
\begin{equation}
 \chi(X,\lambda;\overline{\phi}_N) = \frac{1}{2}\chi_0(\lambda;\overline{\phi}_N) + \sum_{n=1}^{\infty}{\chi_n(\lambda;\overline{\phi}_N)\cos(2n\pi X)}~~\forall~~X\in \left[-\frac{1}{2},\frac{1}{2}\right],  \label{eqn3.69}
\end{equation}
with coefficients,
\begin{equation}
   \chi_n(\lambda;\overline{\phi}_N) = 4\int_{0}^{\frac{1}{2}}{\chi(s,\lambda;\overline{\phi}_N)\cos(2n\pi s)}ds,~~n=0,1,2..~.\label{eqn3.70}  
\end{equation}
For admissible kernel perturbations in $\overline{K}(\mathbb{R})$ that are continuous it follows from the bounds in (\ref{eqn3.62}) and the above regularity on $\chi$ that the real sequence $\{\chi_n(\lambda;\overline{\phi}_N)\}\in l^1$ with its norm bounded in terms of $M$ and $\gamma$, independent of $\lambda\in (1/2,1)$. In particular there is a positive sequence $\{c_r\} \in l_1$, depending only upon $M$ and $\gamma$, such that $|\chi_n(\lambda;\overline{\phi}_N)| \le c_r$ for all $r=1,2,3..$ and $\lambda \in (1/2,1)$.
 Therefore the convergence in (\ref{eqn3.69}) will be absolute and uniform for $X\in \left[-\frac{1}{2},\frac{1}{2}\right]$. Under these circumstances we may substitute from (\ref{eqn3.69}) into (\ref{eqn3.65}) to obtain, using (\ref{eqn3.59}) and (\ref{eqn3.60}), the following representation of the nonlocal term in equation (\ref{eqn3.56}),
\begin{multline}
  \overline{J}(\mathcal{F})(x,\lambda,\overline{D};\overline{\phi}_N) = \frac{1}{2}\chi_0(\lambda;\overline{\phi}_N)~ + \\  \sum_{n=1}^{\infty}{\left(\int_{-a(\lambda)}^{a(\lambda)}{\mathcal{F}(s,\lambda,\overline{D};\overline{\phi}_N)\cos(2n\pi s)}ds\right)\chi_n(\lambda;\overline{\phi}_N)\cos(2n\pi x)}  \label{eqn3.71}
\end{multline}
for all $x\in [-a(\lambda),a(\lambda)] $. Having obtained a suitable representation of the nonlocal term, we now return to considering [NBVP]. We begin by introducing  the following linear, regular Sturm-Liouville eigenvalue problem, which we refer to as [SL$(\mathbf{I},a)$],
\begin{equation}
  \mathcal{F}_L'' + \left(\alpha_L - \sum_{n=1}^{\infty}{I_n\cos(2n\pi x)}\right)\mathcal{F}_L = 0,~~x\in (-a,a),  \label{eqn3.72}
\end{equation}
\begin{equation}
    \mathcal{F}_L(-a) = \mathcal{F}_L(a)=0. \label{eqn3.73}
\end{equation}
Here, $a=a(\lambda)$ throughout, and the eigenvalue is $\alpha_L \in \mathbb{R}$, whilst $\mathbf{I}=(I_1,I_2,I_3,...,..)\in l^1$ is a given sequence of real parameters. In addition we impose the conditions,
\begin{equation}
    \mathcal{F}_L(x)>0~~\text{and  even}~~\forall~~x\in (-a,a) \label{eqn3.74}
\end{equation}
and
\begin{equation}
    \int_{-a}^a{\mathcal{F}_L(s)}ds = 1.   \label{eqn3.75}
\end{equation}
We first observe that $Q:[-a,a]\times l^1\to \mathbb{R}$, given by
\begin{equation}
   Q(x,\mathbf{I}) = \sum_{n=1}^{\infty}{I_n\cos(2n\pi x)}  \label{eqn3.75'} 
\end{equation}
has $Q\in C([-a,a]\times l^1)$ (and is  uniformly Lipschitz continuous in $l^1$, with a Lipschitz constant of unity) and so the classical Sturm-Liouville Theory (see, for example, Coddington and Levinson \cite{CoddingtonLevinson}, chapter 5) establishes that, for any given $\mathbf{I}\in l^1$, [SL$(\mathbf{I},a)$] has a unique solution, and this corresponds to taking the principal eigenvalue and its associated eigenfunction (normalised in $L^1([-a,a])$), which we write as,
\begin{equation}
    \mathcal{F}_L=\mathcal{F}_L(x,\mathbf{I},a)\in C^{2,0,0}\left(\{(x,\mathbf{I},a):x\in [-a,a],\mathbf{I}\in l^1, a\in \left(0,\frac{1}{4}\right)\}\right),  \label{eqn3.76}
\end{equation}
\begin{equation}
  \alpha_L=\alpha_L(\mathbf{I},a)\in C\left(l^1 \times\left(0,\frac{1}{4}\right)\right).  \label{eqn3.76'}  
\end{equation}
We note that it is readily established directly, using the turning point properties of principal eigenfunctions, that,
\begin{equation}
    \alpha_L(\mathbf{I},a) \ge \sum_{n=1}^{\infty}{I_n}. \label{eqn3.76''}
\end{equation}
We next introduce the mapping $\mathbf{G}:l^1\times (0,\frac{1}{4}) \to l^1$, with $\mathbf{G} = (G_1,G_2,G_3,....)$, and, for each $r=1,2,3,..$, $G_r:l^1\times (0,\frac{1}{4}) \to \mathbb{R}$ is given by,
\begin{equation}
    G_r(\mathbf{I},a) = \int_{-a}^a{\mathcal{F}_L(y,\mathbf{I},a)\cos(2r\pi y)}dy.   \label{eqn3.77}
\end{equation}
It follows from (\ref{eqn3.73})-(\ref{eqn3.75}), after two integrations by parts, that, for each $r=1,2,3..$,
\begin{equation}
   |G_r(\mathbf{I},a)| \le \min\{1, \frac{1}{2r^2\pi^2} (||\mathbf{I}||_1 + |\alpha_L(\mathbf{I},a)|)\}~~\forall~~(\mathbf{I},a)\in l^1\times (0,\frac{1}{4}), \label{eqn3.78}  
\end{equation}
which ensures that $\text{Im}(\mathbf{G}) \subseteq l^1$, with,
\begin{equation}
    ||\mathbf{G}(\mathbf{I},a)||_1 \le \frac{1}{12} (||\mathbf{I}||_1 + |\alpha_L(\mathbf{I},a)|)~~\forall~~(\mathbf{I},a)\in l^1\times (0,\frac{1}{4}).  \label{eqn3.79}
\end{equation}
In addition, it follows from (\ref{eqn3.76}) and (\ref{eqn3.78}) (specifically with the quadratic decay rate in r) that $\mathbf{G}$ is uniformly continuous and bounded on any compact subset of $l^1\times \left(0,\frac{1}{4}\right)$. It will be of use later to note that simple direct calculations determine that,
\begin{equation}
    \mathcal{F}_L(x,\mathbf{0},a) = \frac{\pi}{4a}\cos\left(\frac{\pi x}{2a}\right), \label{eqn3.79'}
\end{equation}
\begin{equation}
\alpha_L(\mathbf{0},a) = \frac{\pi^2}{4a^2}, \label{eqn3.79''}
\end{equation}
\begin{equation}
    G_r(\mathbf{0},a) = -(16a^2r^2-1)^{-1}\cos(2ar\pi), \label{eqn3.79'''}
\end{equation}
with $x\in [-a,a]$, $a\in (0,\frac{1}{4})$ and $r=1,2,3..$.

It is now straightforward to establish that there is a one to one correspondence between solutions to [NBVP], and solutions $\mathbf{I}=\mathbf{I}(a,\overline{D})\in l^1$ to the nonlinear equation in $l^1$ given by,
\begin{equation}
  \overline{D}^{-1}  \text{diag}[\chi_1,\chi_2,...](a;\overline{\phi}_N)\mathbf{G}(\mathbf{I},a) = \mathbf{I},\label{eqn3.80}
\end{equation}
with $(\mathbf{I},a,\overline{D})\in l^1\times\left(0,\frac{1}{4}\right)\times\mathbb{R}^+$.  Recalling that $\{\chi_r(a,\overline{\phi}_N)\}\in l_1$ is bounded independently of $a\in(0,1/4)$, together with the inequalities (\ref{eqn3.78}), we can see directly, that the image, in $l_1$, of the left side of equation (\ref{eqn3.80}) is contained in the compact (closed, bounded and equismall) subset 
\begin{equation}
   K(\overline{D}) = \{ \{\alpha_r\} \in l_1: |\alpha_r| \le \overline{D}^{-1}c_r,~r=1,2,..\},
\end{equation}
where the sequence $\{c_r\}\in l_1$ is as introduced earlier, depending only upon the constants $M$ and $\gamma$. Thus, at any $(a,\overline{D}) \in (0,1/4)\times \mathbb{R}^+$, we need only seek solutions of (\ref{eqn3.80}) which reside in the compact subset $K(\overline{D})$.The corresponding solution to [NBVP] is then given by,
\begin{equation}
  \mathcal{F}(x,\lambda,\overline{D};\overline{\phi}_N) = \mathcal{F}_L(x,\mathbf{I}(a,\overline{D}),a), \label{eqn3.81}
\end{equation}
\begin{equation}
    \mathcal{I} = \alpha_L(\mathbf{I}(a,\overline{D}),a) + \frac{1}{2\overline{D}}\chi_0(a;\overline{\phi}_N), \label{eqn3.82}
\end{equation}
with $x\in [-a(\lambda),a(\lambda)]$, $\lambda\in \left(\frac{1}{2},1\right)$ and $\overline{D}\in \mathbb{R}^+$. We note that since any solution $\mathbf{I}(a,\overline{D}) \in l_1$ to (\ref{eqn3.80}) must be in   $K(\overline{D})$, then,
\begin{equation}
    ||\mathbf{I}(a,\overline{D})||_1 \le \overline{D}^{-1}\sum_{r=1}^{\infty}{|\chi_r(a,\overline{\phi}_N)|} \le \overline{D}^{-1} \sum_{r=1}^{\infty}c_r, \label{eqn3.82'}
\end{equation}
for $(a,\overline{D})\in (0,\frac{1}{4})\times \mathbb{R}^+$. We now fix $(a,\overline{D}) \in (0,1/4)\times \mathbb{R}^+$, and regard the left side of equation (\ref{eqn3.80}) as a mapping from the subset $K(\overline{D})$ to $l_1$. It has been demonstrated above that the image of this map is, in fact, contained in the subset $K(\overline{D})$. Also, as $K(\overline{D})$ is compact, we have also demonstrated above that the mapping is continuous on $K(\overline{D})$. Finally, it is an immediate consequence of its definition, that the subset $K(\overline{D})$ is convex. We can then conclude from the Schauder Fixed Point Theorem (see, for example, Kantorovitch and Akilov \cite{Kant}) that the mapping has at least one fixed point in the subset $K(\overline{D})$. Therefore we may conclude that at each fixed $(a,\overline{D}) \in (0,1/4)\times \mathbb{R}^+$, equation (\ref{eqn3.80}) has at least one solution $\mathbf{I}=\mathbf{I}(a,\overline{D})\in l^1$ in the subset $K(\overline{D})$ (and no solutions outside of this subset). The continuity of the mapping in $(a,\overline{D})$ on compact subsets of $(0,1/4)\times \mathbb{R}^+$ also enables us to deduce that there will be a principal  solution branch, say $\mathbf{I} = \mathbf{I}_p(a,\overline{D})$, which is defined and continuous on $(0,1/4)\times \mathbb{R}^+$, with the possibility of additional continuous solution branches bifurcating upon crossing possible continuous bifurcation curves on $(0,1/4)\times \mathbb{R}^+$. This completes the general formulation and results. We next delve further into equation (\ref{eqn3.80}) by examining its structure, when $(a,\overline{D}) \in (0,1/4) \times \mathbb{R}^+$ , in the cases of two specific, but representative admissible kernel perturbations.

First, it is interesting to examine solutions to equation (\ref{eqn3.80}) in the limit $\overline{D} \to \infty$. For a fixed admissible kernel perturbation with $||\overline{\phi}||_1^m$ small, this has the geometric interpretation, in the domain $\Omega_1$ of the positive quadrant of the $(\lambda, D)$-plane, of passing from the upper region of the subdomain $\Omega^-(\overline{\phi})$ to the lower region of the complementary subdomain $\Omega^+(\overline{\phi})$ with $\lambda \in (1/2,1)$ (or correspondingly $a\in(0,1/4)$) fixed. We observe directly from (\ref{eqn3.82'}) that any solution $\mathbf{I}(a,\overline{D}) \in K(\overline{D})$ must have $||\mathbf{I}(a,\overline{D})||_1 \to 0$ as $\overline{D} \to \infty$, uniformly for $a\in (0,1/4)$. In particular, equation (\ref{eqn3.80}) then determines that there is exactly one solution at each point $(a,\overline{D}) \in (0,1/4) \times \mathbb{R}^+$ when $\overline{D}$ is sufficiently large, with this solution forming the genesis of the principal branch, and it has the asymptotic form,
\begin{equation}
  \mathbf{I}(a,\overline{D}) \sim  \overline{D}^{-1} \text{diag}[\chi_1,\chi_2,...](a;\overline{\phi}_N)\mathbf{G}(\mathbf{0},a)~~\text{as}~~\overline{D}\to \infty,  \label{eqnabc}
\end{equation}
uniformly for $a$ in compact subsets of $(0,1/4)$, with $\mathbf{G}(\mathbf{0},a)$ as given in (\ref{eqn3.79'''}). It then follows from (\ref{eqn3.79'}) and (\ref{eqn3.81}) that,
\begin{equation}
  \mathcal{F}(x,\lambda,\overline{D};\overline{\phi}_N) = \mathcal{F}_L(x,\mathbf{I}(a,\overline{D}),a) = \frac{\pi}{4a}\cos\left(\frac{\pi x}{2a}\right) + O(\overline{D}^{-1})~~\text{as}~~\overline{D} \to \infty, \label{eqn3.81'}
\end{equation}
uniformly for $x\in [-a,a]$ and $a$ in compact subsets of $(0,1/4)$. Therefore we have a unique positive periodic state at each $(a,\overline{D}) \in (0,1/4) \times \mathbb{R}^+$, with $\overline{D}$ sufficiently large, and, via
(\ref{eqn3.55}) with (\ref{eqn3.81'}),this is given by,
\begin{equation}
  \mathcal{F}_p(x,\lambda,D;\overline{\phi}) = \frac{\pi}{4a}\cos\left(\frac{\pi x}{2a}\right) + O(\overline{D}^{-1},||\overline{\phi}||_1^m)~~\text{as}~~||\overline{\phi}||_1^m \to 0, \label{eqn3.81''}
\end{equation}
with $\overline{D}$ large, and the approximation being uniform for $x\in [-a,a]$ and $a$ in compact intervals of $(0,1/4)$. Conversely, returning to subsection 3.1, we saw that when $(\lambda,D)\in \Omega^+(\overline{\phi})$, we constructed a unique periodic steady state with  $||\overline{\phi}||_1^m$ small, which took on the following form (via (\ref{eqn3.5}), (\ref{eqn3.29}) and (\ref{eqn3.50})) as $D$ became small, namely,
\begin{equation}
  \mathcal{F}_p(x,\lambda,D;\overline{\phi}) = \frac{\pi}{4a}\cos\left(\frac{\pi x}{2a}\right) + O(D^{-1}||\overline{\phi}||_1^m)~~\text{as}~~||\overline{\phi}||_1^m \to 0, \label{eqn3.81'''}
\end{equation}
 with the approximation being uniform again for $x\in [-a,a]$ and $a$ in compact intervals of $(0,1/4)$. It follows directly from the form of $\mathcal{F}_p$ in $\Omega^+(\overline{\phi})$ when $D$ is small (in (\ref{eqn3.81''})) compared with its form in $\Omega^-(\overline{\phi})$ when $\overline{D} = D ||\overline{D}||_1^m$ is large (in(\ref{eqn3.81''}) that they agree according to the general asymptotic matching principle of Van Dyke[], which determines that the unique branch of periodic steady states appearing as a regular perturbation in $\Omega^+(\overline{\phi})$ is continued into the principal branch of periodic steady states identified above as a singular perturbation  in $\Omega^-(\overline{\phi})$, and thereafter, further reductions in $\overline{D}$ may lead to more exotic secondary global bifurcations, features which we investigate in the next section for a number of specific representative admissible kernel perturbations. Before doing this it is useful to summarise the key conclusions of the subsection in the following:

 \begin{theorem}
      Let $\overline{\phi}\in \overline{K}(\mathbb{R})$ be an admissible kernel perturbation (which is also either everywhere continuous and satisfies (\ref{eqn3.62}), or has support contained in $[-1/2,1/2]$) with $||\overline{\phi}||_1^m$ small. Then for each $(\lambda,D) \in \Omega^-(\overline{\phi})$ there is at least one positive periodic steady state associated with (IBVP)$_p$, with each being a singular perturbation of the unique unperturbed periodic steady state, at the same point $(\lambda,D)$, associated with (IBVP), in the sense that the magnitude of the perturbation is now of $O(1)$ in the supnorm as $||\overline{\phi}||_1^m \to 0$, uniformly on compact subsets of $\Omega^-(\overline{\phi})$. There is a principle branch which is the natural continuation of the unique branch identified in $\Omega^+(\overline{\phi})$. In addition there is the possibility of further branches emerging at secondary bifurcations as we decrease $\overline{D} = D ||\overline{\phi}||_1^m$ in $\Omega^-(\overline{\phi})$.      
 \end{theorem}

\subsubsection{Simple Representative Cases}
We next examine equation (\ref{eqn3.80}) in a number of simple, yet representative and tractable, cases. We restrict attention to those cases when the top hat kernel is perturbed only within its support, so that $\text{supp}(\overline{\phi}_N)\subseteq [-\frac{1}{2},\frac{1}{2}]$, and then, from (\ref{eqn3.65'}) and (\ref{eqn3.70}),
 \begin{equation}
   \chi_r(\lambda;\overline{\phi}_N) \equiv \overline{\phi}_{N,r} =  4\int_{0}^{\frac{1}{2}}{\overline{\phi}_N(s)\cos(2r\pi s)}ds,~~n=0,1,2..~.\label{eqn3.83}  
\end{equation}
which, in this case are now \emph{independent of $\lambda$, and hence $a$}. Also,
  \begin{equation}
  \chi_0(\lambda;\overline{\phi}_N) \equiv \overline{\phi}_{N,0} =  4\int_{0}^{\frac{1}{2}}{\overline{\phi}_N(s)}ds = 0, \label{eqn3.84}   
\end{equation}
via (P5) in section 1. The two cases of primary interest, are those contrasting cases when the kernel perturbation focuses nonlocalisation towards the centre of the support of the top hat kernel, or disperses nonlocalisation towards the edges of the support of the top hat kernel. The two most simple cases encompassing these two features are realised by taking,
\begin{equation}
    \overline{\phi}_N(x) \equiv \overline{\phi}_{\pm}(x) = \pm\frac{1}{2} \pi \cos(2\pi x)~~\forall~~x\in \left[-\frac{1}{2},\frac{1}{2}\right],  \label{eqn3.85}
\end{equation}
and so,
\begin{equation}
    \overline{\phi}_{N,1} = \pm \frac{1}{2} \pi   \label{eqn3.86}
\end{equation}
whilst
\begin{equation}
    \overline{\phi}_{N,r} = 0,~~r=2,3,4,... .\label{eqn3.87}
\end{equation}
We now consider these two cases in detail, and note that in each of these cases the infinite dimensional Banach space $l^1$ can be replaced at appropriate stages by simply $\mathbb{R}$.

In these two cases it follows immediately from (\ref{eqn3.86}) and(\ref{eqn3.87}) that equation (\ref{eqn3.80}) requires any solution to have,
\begin{equation}
    \mathbf{I}(a,\overline{D})=(I(a,\overline{D}),0,0,..,0,..), \label{eqn3.88}
\end{equation}
with $I(a,\overline{D})\in \mathbb{R}$ now satisfying the nonlinear scalar equation,
\begin{equation}
  G_1(I,a) = \pm\frac{2}{\pi}\overline{D}I. \label{eqn3.89}
\end{equation}
As a consequence, in $[SL(\mathbf{I},a)]$, we have the simplification,
\begin{equation}
    Q(x,\mathbf{I}) = I\cos(2\pi x),~~x\in [-a,a], \label{eqn3.90}
\end{equation}
so that now $Q$ is a \emph{holomorphic} function of $(x,I)\in [-a,a]\times \mathbb{R}$, from which it follows that we can guarantee, at least,
\begin{equation}
    \mathcal{F}_L= \mathcal{F}_L(x,I,a)\in C^{2,1,1}\left(\left\{(x,I,a): x\in [-a,a], I\in \mathbb{R}, a\in \left(0,\frac{1}{4} \right)\right\}\right). \label{eqn3.91}
\end{equation}
Thus, via (\ref{eqn3.77}), we now have,
\begin{equation}
 G_1 = G_1(I,a)\in C^{1,1}\left(\mathbb{R},\left(0,\frac{1}{4}\right)\right). \label{eqn3.92}
\end{equation}
Also, from (\ref{eqn3.78}) and (\ref{eqn3.79'''}), we observe that,
\begin{equation}
   0< G_1(I,a)<1~~\forall~~(I,a)\in \mathbb{R}\times \left(0,\frac{1}{4}\right)\label{eqn3.93}
\end{equation}
and
\begin{equation}
    G_1(0,a) = (1-16a^2)^{-1}\cos(2\pi a).\label{eqn3.94}
\end{equation}
Using the above increased regularity, a lengthy, but straightforward calculation (via differentiating through [SL$(I,a)$] with respect to $I$ or $a$) also allows us to establish that
\begin{equation}
    (G_1)_I(I,a),~ (G_1)_a(I,a) < 0~~\forall~~(I,a)\in \mathbb{R}\times \left(0,\frac{1}{4}\right). \label{eqn3.96}
\end{equation}
With these monotonicity conditions, the general theory of the last section, together now with the local Implicit Function Theorem, determines that, at each $(a,\overline{D}) \in (0,1/4) \times \mathbb{R}^+$:
\begin{itemize}
    \item When the kernel perturbation has the positive sign, equation (\ref{eqn3.89}) has exactly one solution, $I(a,\overline{D})$, and this lies in the interval $(0,\pi/2\overline{D})$, and depends continuously and continuously differentiably, on $(a,\overline{D})$. It can be viewed as forming a smooth single-valued surface, above the open rectangle $(0,1/4) \times \mathbb{R}^+$, in $(a,\overline{D},I)$ space. The partial derivatives are given by
    \begin{equation}
        I_a(a,\overline{D}) = - \mathcal{J}_+(a,I(a,\overline{D}),\overline{D})^{-1}(G_1)_a(I(a,\overline{D}),a)<0,
        \label{a}
    \end{equation}
    and
    \begin{equation}
        I_{\overline{D}}(a,\overline{D}) = \frac{2}{\pi} \mathcal{J}_+(a,I(a,\overline{D}),\overline{D})^{-1} I(a,\overline{D})<0, \label{aa}
    \end{equation}
    where the Jacobian function is given by,
    \begin{equation}
        \mathcal{J}_+(a,I,\overline{D}) = (G_1)_I(I,a) - \frac{2}{\pi}\overline{D} <0.
    \end{equation}
    \item When the kernel perturbation has the negative sign equation (\ref{eqn3.89}) has one or more solutions, each of which must lie in the interval $(-\pi/2\overline{D},0)$. In this case the associated Jacobian function is given by
    \begin{equation}
        \mathcal{J}_-(a,I,\overline{D}) = (G_1)_I(I,a) + \frac{2}{\pi}\overline{D} \label{eqn***}
    \end{equation}
    which does have the possibility of vanishing, and, if and when this occurs, there will be local secondary bifurcations which create the possibility of multiple solutions. Again, these can be viewed as forming families of smooth single-valued surfaces, joined together at folding bifurcation curves, above the open rectangle $(0,1/4) \times \mathbb{R}^+$, in $(a,\overline{D},I)$ space. Further details in this case will now have to wait until additional prerequisite information has been obtained. However, we do recall, from the general theory above, that when $\overline{D}$ is sufficiently large,equation (\ref{eqn3.98}) has a unique solution for each $a\in (0,1/4)$, with any secondary bifurcations occurring across curves at lower values of $\overline{D}$.
    \end{itemize}
    To examine these cases in more detail, we next consider the asymptotic forms of $\mathcal{F}_L$, $\alpha_L$ and hence $G_1$ in a number of useful limits in the $(I,a)$ parameter plane. In each of these limits [SL$(I,a)$] can be examined directly by a WKB-type analysis:\\
    
$(i)~ a=O(1)$ \emph{with} $|I|$ \emph{large}

In this case we first address the limit $I\to \infty$ with $a\in \left(o(1),\frac{1}{4}-o(1)\right)$. In this limit, for $x\in [-a,a]$, $Q(x,I)$ has its global minimum at $x=\pm a$ and so the structure of [SL$(I,a)$] determines that $\alpha_L(I,a)\sim I\cos(2\pi a) + I^{(1-q)}\overline{\alpha}(a)$ as $I\to \infty$, with
\begin{equation}
 \mathcal{F}_L(x,I,a) \sim
 \begin{cases}
 f_L(x,a)I^s,~~x\in [-a, -a + O(I^{-r})),\\
 E_L(x,I,a),~~x\in (-a + O(I^{-r}), a - O(I^{-r})),\\
 f_L(-x,a)I^s,~~x\in (a - O(I^{-r}), a], \label{eqn3.97}
 \end{cases}
\end{equation}
as $I\to \infty$, $q,r,s>0$ to be determined, and $E_L(x,I,a)$ being exponentially small in $I$ as $I\to \infty$ uniformly for $x\in (-a + O(I^{-r}), a - O(I^{-r}))$. Our objective is now to determine $f_L$ and $\overline{\alpha}$ together with $q,r$ and $s$. First we introduce the scaled coordinate $\overline{x}$ by,
\begin{equation}
    x = -a + {I^{-r}}\overline{x}, \label{eqn3.98}
\end{equation}
with $\overline{x}=O(1)^+$ as $I\to \infty$ in the edge region at $x=-a$. On substitution from (\ref{eqn3.97}) and (\ref{eqn3.98}) into equation (\ref{eqn3.72}) we require, for a nontrivial balance at leading order,
\begin{equation}
r = q = s = \frac{1}{3}. \label{eqn3.99}
\end{equation}
The leading order problem for $f_L$ in this edge region is then obtained from [SL$(I,a)$] as,
\begin{equation}
    f_L'' + (\overline{\alpha} - 2\pi \overline{x} \sin(2\pi a))f_L = 0,~~ \overline{x}>0, \label{eqn3.100}
\end{equation}
\begin{equation}
f_L(0,a)=0,~~f_L(\overline{x},a)\to 0~\text{as}~\overline{x}\to \infty,
\label{eqn101}
\end{equation}
\begin{equation}
    f_L(\overline{x},a)>0~~\forall~~\overline{x}>0, \label{eqn102}
\end{equation}
\begin{equation}
    \int_0^{\infty}{f_L(s,a)}ds = \frac{1}{2}. \label{eqn103}
\end{equation}
This eigenvalue problem can be solved exactly by applying an affine coordinate transformation to reduce the ODE to Airy's equation. The solution is,
\begin{equation}
    f_L(\overline{x},a) = \frac{(2\pi \sin(2\pi a))^{\frac{1}{3}}}{2J}A_i\left((2\pi\sin(2\pi a))^{\frac{1}{3}}\overline{x} + a_i\right),~\overline{x}\ge0, \label{eqn3.104}
\end{equation}
with
\begin{equation}
\overline{\alpha}(a) = - (2\pi \sin(2\pi a))^{\frac{2}{3}}a_i, \label{eqn3.105}
\end{equation}
where
\begin{equation}
J = \int_{a_i}^{\infty}{A_i(s)}ds \approx 1.274, \label{eqn3.106}
\end{equation}
$A_i$ is the usual Airy function and $a_i \approx -2.338$ is the largest zero of $A_i$. We note that we may write,
\begin{equation}
  \mathcal{F}_L(x,I,a)\sim \delta_A(x+a) + \delta_A(x-a)~~\text{as}~~I\to \infty, \label{eqn3.106'}  
\end{equation}
 with $x\in [-a,a]$, where we interpret $\delta_A(\cdot)$ as the Dirac delta function, limiting from the above Airy function sequence as $I\to \infty$. This representation will be convenient at later stages. It now follows from (\ref{eqn3.97}), (\ref{eqn3.77}) and (\ref{eqn3.104}), after some calculation, that,
\begin{equation}
    G_1(I,a)\sim \cos(2\pi a)~~\text{as}~~I\to \infty. \label{eqn3.107}
\end{equation}

We now consider the complementary limit $I\to -\infty$ with $a\in \left(o(1),\frac{1}{4}-o(1)\right)$.
In this limit, for $x\in [-a,a]$, $Q(x,I)$ now has its global minimum at $x=0$ and so the structure of [SL$(I,a)$] determines that $\alpha_L(I,a)\sim I+ |I|^{(1-q')}\tilde{\alpha}(a)$ as $I\to -\infty$, with,
\begin{equation}
 \mathcal{F}_L(x,I,a) \sim
 \begin{cases}
 \tilde{E}_L(x,a),~~x\in [-a, - O(|I|^{-r')})),\\
 \tilde{f}_L(x,I,a)|I|^{s'},~~x\in (- O(|I|^{-r'}),+ O(|I|^{-r'})),\\
 \tilde{E}_L(-x,a),~~x\in (+ O(|I|^{-r'}), a], \label{eqn3.108}
 \end{cases}
\end{equation}
as $I\to -\infty$, $q',r',s'>0$ to be determined, and $\tilde{E}_L(x,I,a)$ being exponentially small in $I$ as $I\to -\infty$ uniformly for $x\in [-a, -O(|I|^{-r'}))$. Our objective is now to determine $\tilde{f}_L$ and $\tilde{\alpha}$ together with $q',r'$ and $s'$. First we introduce the scaled coordinate $\tilde{x}$ by,
\begin{equation}
    x = {|I|^{-r'}}\tilde{x}, \label{eqn3.109}
\end{equation}
with $\tilde{x}=O(1)$ as $I\to -\infty$ in the centre region at $x=0$. On substitution from (\ref{eqn3.97}) and (\ref{eqn3.98}) into equation (\ref{eqn3.72}) we require, for a nontrivial balance at leading order,
\begin{equation}
r'= s'= \frac{1}{4},~q'=\frac{1}{2}. \label{eqn3.110}
\end{equation}
The leading order problem for $\tilde{f}_L$ in this centre region is then obtained from [SL$(I,a)$] as,
\begin{equation}
    \tilde{f}_L'' - (\tilde{\alpha} + 2\pi^2 \tilde{x}^2)\tilde{f}_L = 0,~~ \tilde{x}\in \mathbb{R}, \label{eqn3.111}
\end{equation}
\begin{equation}
\tilde{f}_L(\tilde{x},a)=\tilde{f}_L(-\tilde{x},a)~~\forall~~\tilde{x}\in \mathbb{R},~~\tilde{f}_L(\tilde{x},a)\to 0~\text{as}~|\tilde{x}|\to \infty,
\label{eqn112}
\end{equation}
\begin{equation}
    \tilde{f}_L(\tilde{x},a)>0~~\forall~~\tilde{x}\in \mathbb{R}, \label{eqn113}
\end{equation}
\begin{equation}
    \int_{-\infty}^{\infty}{\tilde{f}_L(s,a)}ds = 1. \label{eqn114}
\end{equation}
This eigenvalue problem can be solved exactly, and the solution is readily obtained as,
\begin{equation}
    \tilde{f}_L(\tilde{x},a) = \frac{1}{2^{\frac{1}{4}}}\exp\left(-\frac{\pi}{\sqrt{2}}\tilde{x}^2\right)~~\forall~~\tilde{x}\in \mathbb{R}, \label{eqn3.115}
\end{equation}
with,
\begin{equation}
    \tilde{\alpha}(a) = -\sqrt{2}\pi. \label{eqn3.116}
\end{equation}
In this case we note that we may write,
\begin{equation}
    \mathcal{F}_L(x,I,a)\sim \delta_G(x)~~\text{as}~~I\to -\infty, \label{eqn3.116'}
\end{equation}
 with $x\in [-a,a]$, where we now interpret $\delta_G(\cdot)$ as the Dirac delta function, limiting from the above Gaussian function sequence as $|I|\to \infty$. It now follows from (\ref{eqn3.97}), (\ref{eqn3.77}) and (\ref{eqn3.115}) that,
\begin{equation}
    G_1(I,a)\sim 1~~\text{as}~~I\to -\infty. \label{eqn3.117}
\end{equation}
We now move on to the next case.\\

$(ii)~a \to 0^+$

In this case the natural scalings in [SL$(I,a)$] lead us to write,
\begin{equation}
    I = a^{-4}\hat{I},   \label{eqn3.118}
\end{equation}
and introduce the scaled coordinate,
\begin{equation}
    X = a^{-1}x.
\end{equation}
We then determine that,
\begin{equation}
 \mathcal{F}_L(x,I,a) = a^{-1}\hat{f}(X,\hat{I}) + O(a^2\hat{I}),  \label{eqn3.119}   
\end{equation}
as $a\to 0^{+}$, with $X\in [-1,1]$, and uniformly for $\hat{I}\in \mathbb{R}$, whilst,
\begin{equation}
    \alpha_L(I,a) = a^{-4}(\hat{I} +a^2\hat{\alpha}(\hat{I})) + O(a^2\hat{I}),  \label{eqn3.120}
\end{equation}
as $a\to 0^+$ uniformly for $\hat{I}\in \mathbb{R}$. Here $\hat{f}$ and $\hat{\alpha}$ are respectively the principal, $L^1$ normalised, eigenfunction and its associated principal eigenvalue, of the reduced Sturm-Liouville eigenvalue problem,
\begin{equation}
    \hat{f}_{XX} + (\hat{\alpha} + 2\pi^2\hat{I}X^2)\hat{f} = 0,~~X\in (-1,1), \label{eqnA}
\end{equation}
\begin{equation}
   \hat{f}(\pm1,\hat{I})=0,\label{bcA}
\end{equation}
\begin{equation}
    \hat{f}(X,\hat{I})>0~~\forall~~X\in (-1,1),
\end{equation}
\begin{equation}
    \int_{-1}^{1}{\hat{f}(Y,\hat{I})}dY = 1. \label{eqn3.120'}
\end{equation}
When $|\hat{I}|\ll1$, we have,
\begin{equation}
    \hat{f}(X,\hat{I}) \sim \frac{1}{4}\pi \cos\left(\frac{1}{2}\pi X\right),~~X\in [-1,1], \label{eqn3.121}
\end{equation}
and
\begin{equation}
    \hat{\alpha}\sim \frac{1}{4}\pi^2.
\end{equation}
Also, when $\hat{I}\gg1$, we have,
\begin{equation}
  \hat{f}(X,\hat{I})\sim \delta_A(X-1) + \delta_A(X+1),~~ X\in [-1,1],   \label{eqn3.122}
\end{equation}
and
\begin{equation}
    \hat{\alpha}(\hat{I})\sim -2\pi^2\hat{I}
\end{equation}
whilst, when $(-\hat{I})\gg1$
\begin{equation}
\hat{f}(X,\hat{I}) \sim \delta_G(X),~~X\in [-1,1],  \label{eqn3.123}
\end{equation}
and
\begin{equation}
    \hat{\alpha}(\hat{I}) = o(\hat{I}).
\end{equation}
It now follows from (\ref{eqn3.119}), (\ref{eqn3.77}) and (\ref{eqn3.120'}) that,
\begin{equation}
    G_1(I,a) =  1 - O(a^2)~~\text{as}~~a\to 0^+, \label{eqn3.124}
\end{equation}
uniformly for $\hat{I} \in \mathbb{R}$. The final limiting case is now considered.\\

$(iii)~a \to (1/4)^-$

\begin{figure}
\begin{center}
\includegraphics[width=\textwidth]{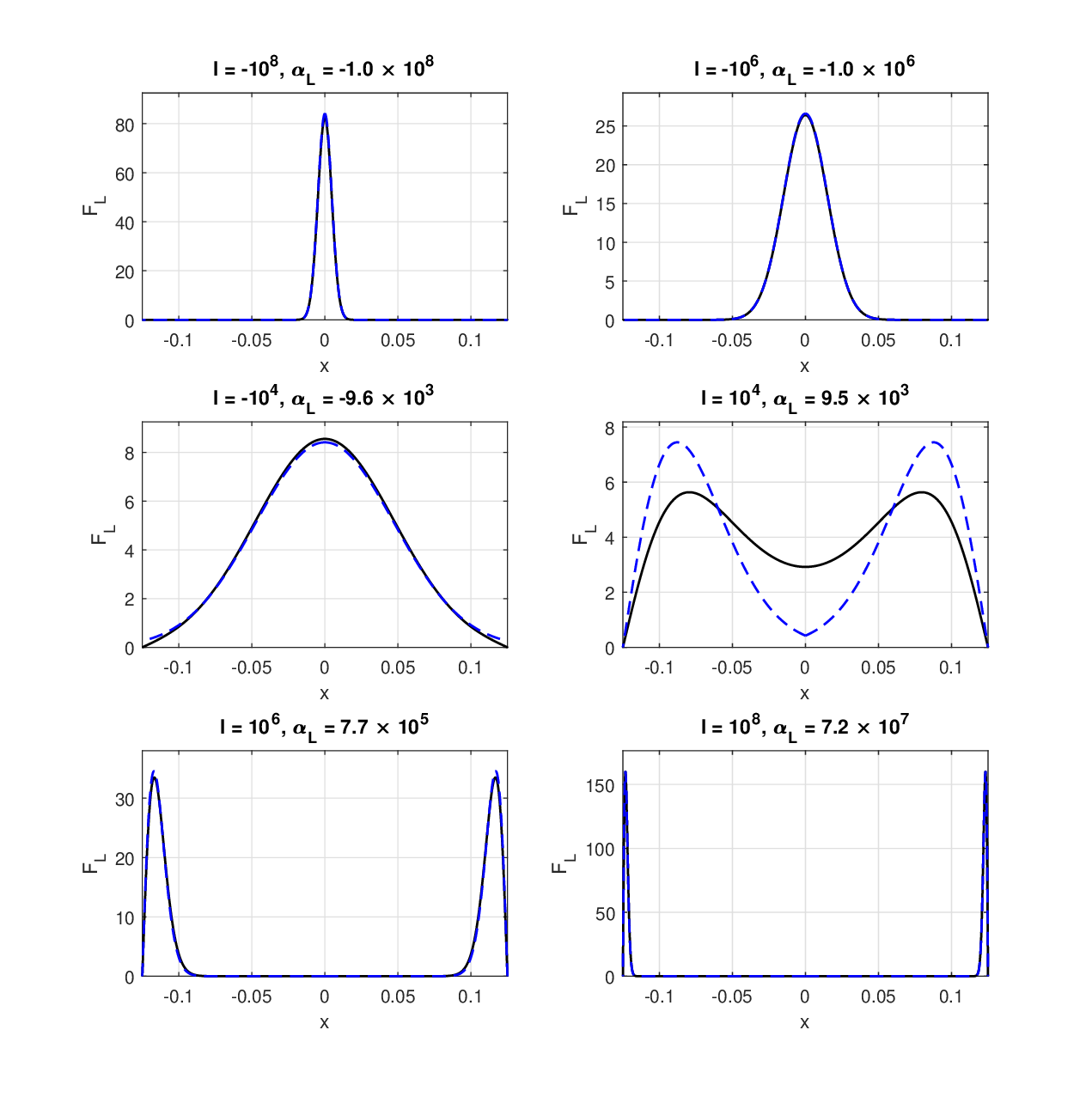}
\caption{The eigenfunction, $\mathcal{F}_L$ and the associated eigenvalue $\alpha_L$ for various values of $I$ when $a = \frac{1}{8}$. The broken blue lines show the leading order asymptotic solution when $|I| \gg 1$.}\label{fig_u}
\end{center}
\end{figure}

\begin{figure}
\begin{center}
\includegraphics[width=\textwidth]{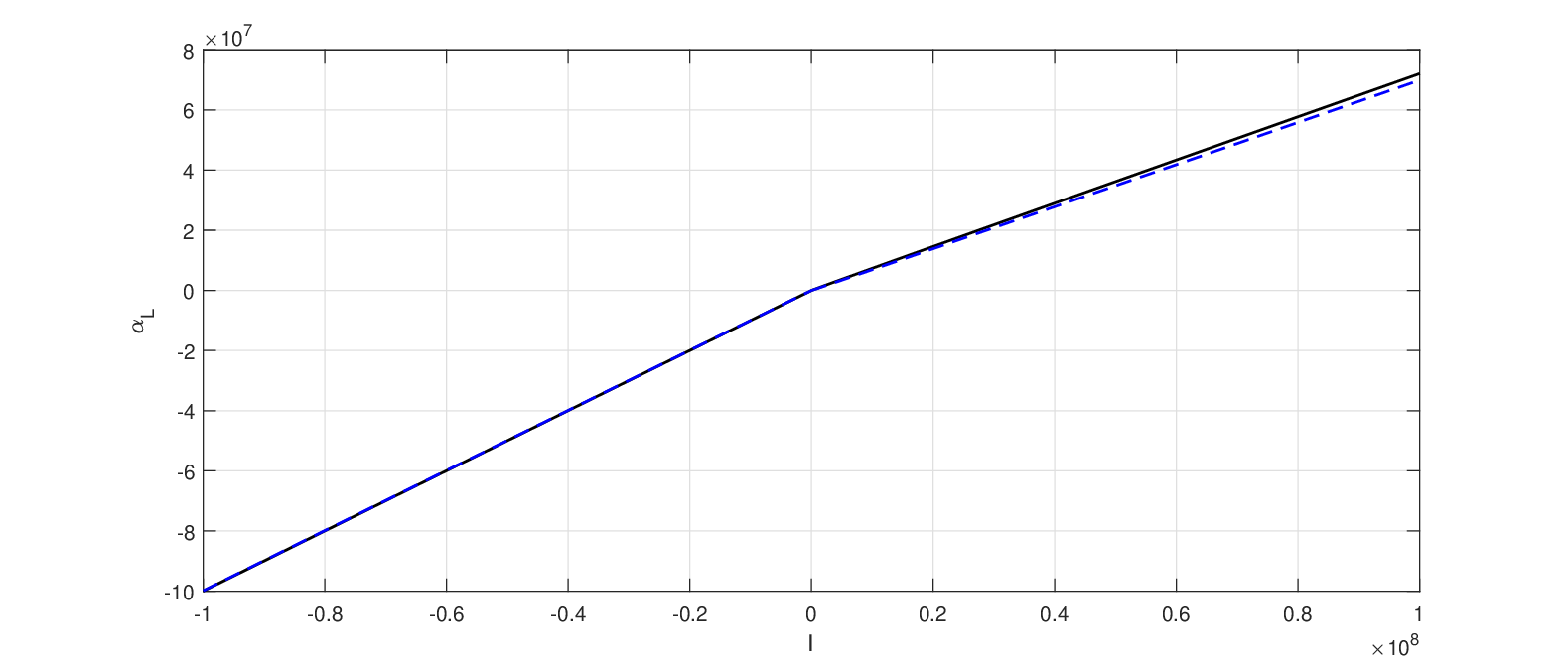}
\caption{A plot of the eigenvalue, $\alpha_L$, as a function of $I$ when $a = \frac{1}{8}$. The broken blue lines show the leading order asymptotic solution when $|I| \gg 1$.}\label{fig_alpha}
\end{center}
\end{figure}

\begin{figure}
\begin{center}
\includegraphics[width=\textwidth]{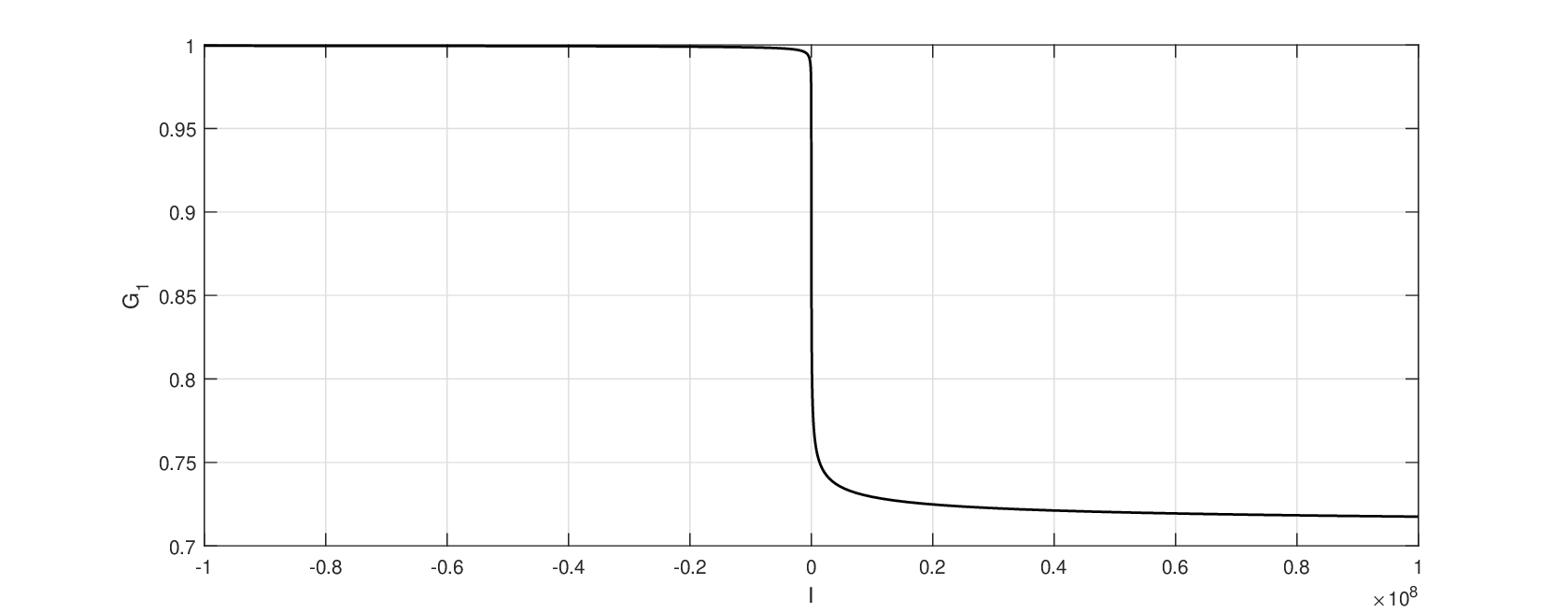}
\caption{A plot of $G_1$ as a function of $I$ when $a = \frac{1}{8}$. Note that $G_1 \to 1$ as $I \to - \infty$ and $G_1 \to \cos(2\pi a) = 1/\sqrt{2}$ as $I \to \infty$, as predicted by the asymptotic solution for $|I|\gg 1$.}\label{fig_G1}
\end{center}
\end{figure}

In this case, for brevity, we omit the detailed asymptotic development, and record only the key results. We write $a = \frac{1}{4} - \overline{a}$ and consider the limit $\overline{a}\to 0^+$. When $I=O(1)$ as $\overline{a}\to 0^+$, we obtain,
\begin{equation}
  \mathcal{F}_L(X,I,a)\sim \overline{f}(X,I),~~X\in [-1,1], \label{eqn3.125'} 
\end{equation}
\begin{equation}
    \alpha_L(I,a)\sim \overline{\alpha}(I), \label{eqnbar}
\end{equation}
\begin{equation}
    G_1(I,a)\sim \int_{-1}^{1}{\overline{f}(Y,I)\cos\left(\frac{1}{2}\pi Y\right)}dY\equiv\overline{g}_1(I).  \label{eqn3.125}
\end{equation}
Here $\overline{f}$ and $\overline{\alpha}$ are the positive principal eigenfunction, normalised in $L^1$, and its associated principal eigenvalue, from the regular Sturm-Liouville eigenvalue problem,
\begin{equation}
\overline{f}_{XX} + \frac{1}{16}\left(\overline{\alpha} - I\cos\left(\frac{1}{2}\pi X\right)\right)\overline{f} = 0,~~X\in (-1,1),
\end{equation}
\begin{equation}
    \overline{f}(\pm1) = 0.
\end{equation}
When $(-I)\gg1$, the approximations in (\ref{eqn3.125'})-(\ref{eqn3.125}) remain uniform, and can be simplified using,
\begin{equation}
    \overline{\alpha}(I)\sim I,
\end{equation}
\begin{equation}
    \overline{f}(X,I)\sim \delta_G(X),~~X\in [-1,1],
\end{equation}
\begin{equation}
    \overline{g}_1(I)\sim 1,
\end{equation}
as $I\to -\infty$. However, when $I\gg1$ the approximations in (\ref{eqn3.125'})-(\ref{eqn3.125}) become nonuniform, and a re-balancing of terms in [SL$(I,a)$] determines that the nonuniformity occurs when $I = O(\overline{a}^{-3})$, with correspondingly, $\overline{\alpha}_L = O(\overline{a}^{-2})$, as $\overline{a} \to 0^+$. We thus introduce $\overline{I} = \overline{a}^3 I = O(1)^+$ as $\overline{a}\to 0^+$ and rescale [SL$(I,a)$] accordingly. Omitting details, we obtain, briefly,
\begin{equation}
 \mathcal{F}_L(X,I,a)\sim \delta_A(X-1) + \delta_A(X+1),~~ X\in [-1,1],   \label{eqn3.126}   
\end{equation}
\begin{equation}
    \alpha_L(I,a)\sim 2\pi \overline{I} \left(1 + c_i\overline{I}^{-\frac{1}{3}}\right)\overline{a}^{-2}, \label{eqm3.127}
\end{equation}
\begin{equation}
   G_1(I,a)\sim 2\pi \left(1 + c_i\overline{I}^{-\frac{1}{3}}\right)\overline{a} .  \label{eqn3.128}
\end{equation}
as $\overline{a}\to 0^+$ with $\overline{I}\ge O(1)^+$. Here the constant $c_i$ is given by,
\begin{equation}
    c_i = (2\pi)^{-\frac{1}{3}} J^{-1} \int_{a_i}^{\infty}{(s - a_i)A_i(s)}ds > 0.
\end{equation}
with $J$ as given in (\ref{eqn3.106}).\\\\
As an illustration, we present a numerical solution to [SL$(I,a)$] with $a=\frac{1}{8}$. In Figure~\ref{fig_u} we graph, at selected values of $I\in \mathbb{R}$, the eigenfunction $\mathcal{F}_L(x,I,a)$ for $x\in [-a,a]$, in Figure~\ref{fig_alpha} we graph the eigenvalue $\alpha_L(I,a)$ against $I$ and finally, in Figure~\ref{fig_G1}, we graph $G_1(I,a)$ against $I$. In each case, there is excellent agreement with the asymptotic solution for $|I|\gg 1$.
The solutions were obtained using a finite difference discretisation, with 500 equally-spaced grid points, solving in $0 \leq x \leq a$, with a symmetry condition applied at $x=0$.

We can now turn to the solution of equation (\ref{eqn3.89}). First (\ref{eqn3.93}), (\ref{eqn3.96}) and (\ref{eqn3.107}) enable a refinement of the earlier bounds on solutions to,
\begin{equation}
 \frac{1}{2}\pi \overline{D}^{-1} \cos(2\pi a) < \pm I(a,\overline{D}) <\frac{1}{2}\pi \overline{D}^{-1}, \label{eqn3.129}
\end{equation}
for the positive and negative kernel perturbations respectively. We now consider each case in turn.

\subsubsection{$\overline{\phi}_N = \overline{\phi}_+$}
In this case, as we have seen earlier, it follows from (\ref{eqn3.92}), (\ref{eqn3.93}) and (\ref{eqn3.96}), together with an application of the Intermediate Value Theorem, that, for each $(a,\overline{D})\in \left(0,\frac{1}{4}\right) \times \mathbb{R}^+$, equation (\ref{eqn3.89}) has a unique solution, which we now label as $I = I^+(a,\overline{D})$, and that this solution is positive, and satisfies the additionally refined inequality,
\begin{equation}
 \frac{1}{2}\pi \overline{D}^{-1} \cos(2\pi a) <  I^+(a,\overline{D})< \frac{1}{2}\pi \overline{D}^{-1}(1-16a^2)^{-1}\cos(2\pi a). \label{eqn3.130}    
\end{equation}
whilst being decreasing in both $a$ and $\overline{D}$, via (\ref{a}) and (\ref{aa}). We conclude that for each $(a,\overline{D})\in \left(0,\frac{1}{4}\right) \times \mathbb{R}^+$, [NBVP] has a unique solution, which is given, via (\ref{eqn3.81}) and (\ref{eqn3.82}), by 
\begin{equation}
    \mathcal{F}(x,\lambda,\overline{D},\overline{\phi}_+) = \mathcal{F}_L(x,I^+(a,\overline{D}),a)~~\forall~~x\in[-a,a],  \label{eqn3.132} 
\end{equation}
\begin{equation}
    \mathcal{I} = \alpha_L(I^+(a,\overline{D}),a).  \label{eqn3.133}
\end{equation}
To determine the structure of $\mathcal{F}$, we next consider the form of $I^+(a,\overline{D})$ for fixed $\overline{D} > 0$, as $a$ increases through the interval $\left(0,\frac{1}{4}\right)$. First, with $\overline{D}=O(1)^+$, it follows from equation (\ref{eqn3.89}) and (\ref{eqn3.124}), that,
\begin{equation}
 I^+(a,\overline{D})= \frac{1}{2}\pi \overline{D}^{-1} - O(a^2)~~\text{as}~~a\to 0^+,  \label{eqn3.134}   
\end{equation}
which we observe is in accord with the bounds in (\ref{eqn3.130}). Conversely, as $a\to \frac{1}{4}^-$, we have, from equation (\ref{eqn3.89}) and (\ref{eqn3.125}), that,
\begin{equation}
  I^+(a,\overline{D}) \to I^*(\overline{D}), \label{eqn3.135}
\end{equation}
where $Y=I^*(\overline{D})$ is the unique positive root of the equation,
\begin{equation}
    \overline{g}_1(Y) - \frac{2}{\pi}\overline{D}Y = 0,
\end{equation}
and it is readily established that $I^*(\overline{D})$ is monotone decreasing with $\overline{D}>0$, and strictly less than $\frac{1}{2}\pi \overline{D}^{-1}$ (in accord with inequalities (\ref{eqn3.129})), whilst,
\begin{equation}
    I^*(\overline{D}) \to 
    \begin{cases}
        0~~\text{as}~~\overline{D}\to \infty,\\
        \infty~~\text{as}~~\overline{D}\to 0.
    \end{cases}
\end{equation}
Between these limiting forms, $I^+(a,\overline{D})$ is monotone decreasing with $a$. We next examine how the structure of $I^+(a,\overline{D})$ against $a\in \left(0,\frac{1}{4}\right)$ develops as $\overline{D}\to \infty$. This follows directly, via (\ref{eqn3.130}), (\ref{eqn3.79'''}), (\ref{eqn3.89}) and (\ref{eqn3.92}), as,
\begin{equation}
    I^+(a,\overline{D}) \sim \frac{1}{2}\pi \overline{D}^{-1}(1 - 16a^2)^{-1} \cos(2\pi a)~~\text{as}~~\overline{D}\to \infty,  \label{eqn3.136}
\end{equation}
uniformly for $a\in \left(0,\frac{1}{4}\right)$. We finally establish the structure of $I^+(a,\overline{D})$ against $a\in \left(0,\frac{1}{4}\right)$ as $\overline{D}\to 0^+$. This is obtained directly from (\ref{eqn3.89}), (\ref{eqn3.107}), (\ref{eqn3.124}) and (\ref{eqn3.128}), as,
\begin{equation}
 I^+(a,\overline{D}) \sim 
 \begin{cases}
     \frac{1}{2}\pi \overline{D}^{-1}\cos(2\pi a),~~a\in \left(0,\frac{1}{4}-O(\overline{D}^{\frac{1}{4}})\right),\\
     \hat{a}^{-3}\psi(\hat{a})\overline{D}^{-\frac{3}{4}},~~a\in \left(\frac{1}{4}-O(\overline{D}^{\frac{1}{4}}), \frac{1}{4}\right),
 \end{cases} \label{eqn3.137}
\end{equation}
as $\overline{D}\to 0^+$, where $\hat{a} = \left(\frac{1}{4} - a\right)\overline{D}^{-\frac{1}{4}} = O(1)^+$ and $\psi=\psi(\hat{a})$ is the unique positive root of the equation,
\begin{equation}
    \psi = \pi^2\hat{a}^4 (1 + c_i\psi^{-\frac{1}{3}}),
\end{equation}
which is monotone increasing with $\hat{a}>0$, and has,
\begin{equation}
    \psi(\hat{a}) \sim
    \begin{cases}
      \pi^{\frac{3}{2}} c_i^{\frac{3}{4}} \hat{a}^3~~\text{as}~~\hat{a}\to 0^+, \\
      \pi^2 \hat{a}^4~~\text{as}~~\hat{a}\to \infty.
    \end{cases}
    \label{eqn3.138}
\end{equation}
We can now use these results to construct $\mathcal{F}$. From (\ref{eqn3.136}), (\ref{eqn3.137}), (\ref{eqn3.132}), and the cases (i)-(iii) earlier, we obtain, after some routine consideration, firstly that, as $\overline{D}\to 0$, with $x\in [0,a]$,
\begin{equation}
    \mathcal{F}(x,\lambda,\overline{D},\overline{\phi}_+) \sim
    \begin{cases}
   \frac{1}{4} a^{-1} \pi \cos\left(\frac{\pi}{2a} x\right),~~a=o(\overline{D}^{\frac{1}{4}})^+,\\
   a^{-1}\hat{f}\left(xa^{-1}, \frac{1}{2}\pi a^4 \overline{D}^{-1}\right),~~a= O(\overline{D}^{\frac{1}{4}})^+,\\
   \delta_A(x+a) + \delta_A(x-a),~~O(\overline{D}^{\frac{1}{4}}) < a \le \frac{1}{4}.
    \end{cases}
\end{equation}
We observe that there exists a structure bifurcation value $a=a_0\overline{D}^{\frac{1}{4}}$ such that $\mathcal{F}$ has single-hump structure for lower values of $a$ and double-hump structure for larger values of $a$. Here $a_0$ is the unique, positive root of the equation,
\begin{equation}
  \hat{\alpha}\left(\frac{1}{2}\pi a_0^4\right) = 0,   
\end{equation}
with $\hat{\alpha}(\hat{I})$ determined via (\ref{eqnA}) - (\ref{eqn3.120'}). When $\hat{\alpha} = 0$ in (\ref{eqnA}), the general solution can be determined analytically in terms of parabolic cylinder functions. This shows that there is no solution that satisfies the boundary condition (\ref{bcA}) when $\hat{I}<0$. When $\hat{I}>0$, a positive solution exists provided that $Y= \left(8 \pi^2 \hat{I}\right)^{\frac{1}{4}}$ is the smallest positive root of
\begin{equation}
-i D_{-\frac{1}{2}}\left(e^{\frac{i\pi}{4}} Y\right) + D_{-\frac{1}{2}}\left(e^{\frac{3i\pi}{4}} Y\right)=0,
\end{equation}
where $D_{-\frac{1}{2}}(\cdot)$ is the Whittaker form of the associated parabolic cylinder function. This gives $a_0 \approx 0.849$.

Conversely, as $\overline{D}\to \infty$, with $x\in [0,a]$, we obtain,
\begin{equation}
 \mathcal{F}(x,\lambda,\overline{D},\overline{\phi}_+) \sim
   \frac{1}{4} a^{-1} \pi \cos\left(\frac{\pi}{2a} x\right), \label{eqnf}
\end{equation}
for each $a\in\left(0,\frac{1}{4}\right)$. We now observe that there is no structure bifurcation, and that $\mathcal{F}$ has a single-hump structure at each $a\in \left(0,\frac{1}{4}\right)$. Finally, when $\overline{D}=O(1)^+$, we have
\begin{equation}
 \mathcal{F}(x,\lambda,\overline{D},\overline{\phi}_+) \sim
   \frac{1}{4} a^{-1} \pi \cos\left(\frac{\pi}{2a} x\right)~~\text{as}~~a\to 0^+,
\end{equation}
for all $x\in [0,a]$. More generally, we determine that there is a value $\overline{D}=\overline{D}^*$, where $Z=\overline{D}^*$ is the unique positive root of the equation,
\begin{equation}
    \overline{\alpha}(I^*(Z))=I^*(Z),  \label{eqn3.140}
\end{equation}
with $I^*(\cdot)$ and $\overline{\alpha}(\cdot)$ as defined in (\ref{eqn3.135}) and (\ref{eqnbar}) respectively, such that:\\

(i) $0<\overline{D}<\overline{D}^*$

There exists $a_c(\overline{D})\in \left(0,\frac{1}{4}\right)$ such that $\mathcal{F}(x,\lambda,\overline{D},\overline{\phi}_+)$ has single-hump structure for $a\in (0,a_c(\overline{D})]$ but double-hump structure for $a\in \left(a_c(\overline{D}),\frac{1}{4}\right)$. Here $a_c(\overline{D})$ is continuous and monotone, with limits
\begin{equation}
     a_c(\overline{D})\to
    \begin{cases}
     0~~\text{as}~~\overline{D} \to 0,\\
     \frac{1}{4}~~\text{as}~~\overline{D} \to \overline{D}^{*},
    \end{cases}    
\end{equation}
and is determined, for each $\overline{D} \in (0,\overline{D}^*]$, as the unique positive root $Z=a_c(\overline{D})$ of the equation,
\begin{equation}
    \alpha_L(I^+(Z,\overline{D})) = I^+(Z,\overline{D}).  \label{eqn3.141}
\end{equation}

\begin{figure}
\begin{center}
\includegraphics[width=\textwidth]{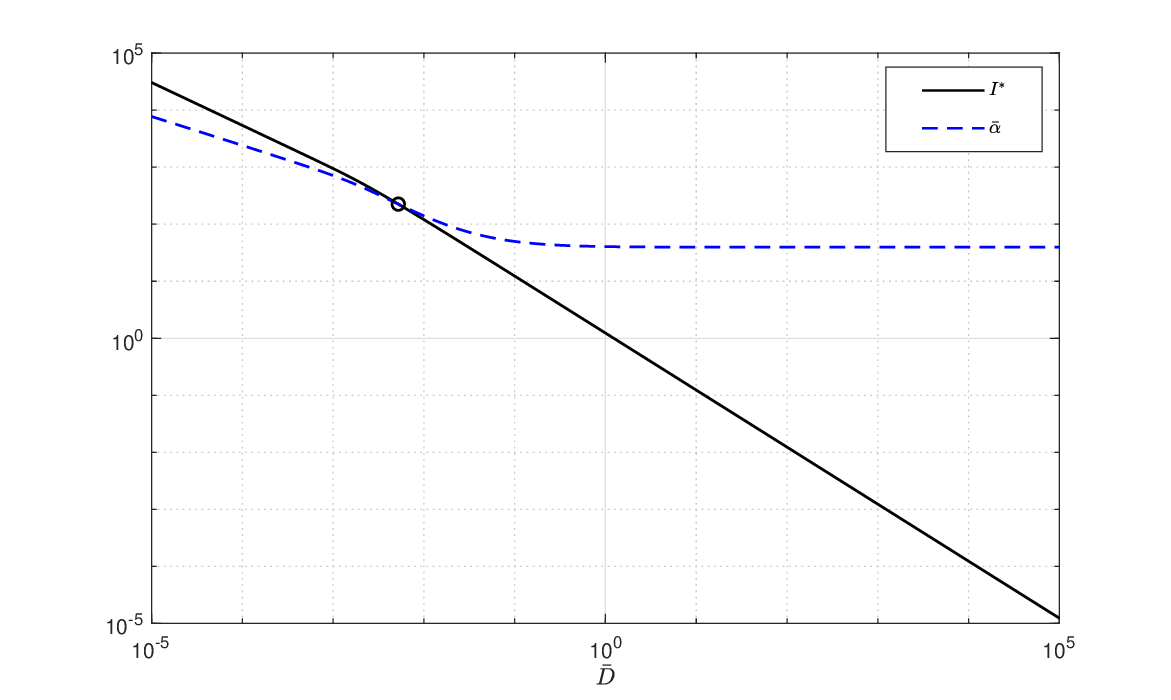}
\caption{A plot of the functions $I^*(\bar{D})$ and $\bar{\alpha}(I^*(\bar{D}))$, whose intersection, indicated by a circle, lies at $\bar{D} = \bar{D}^* \approx 5.22 \times 10^{-3}$.}\label{fig_Dstar}
\end{center}
\end{figure}

\begin{figure}
\begin{center}
\includegraphics[width=0.9\textwidth]{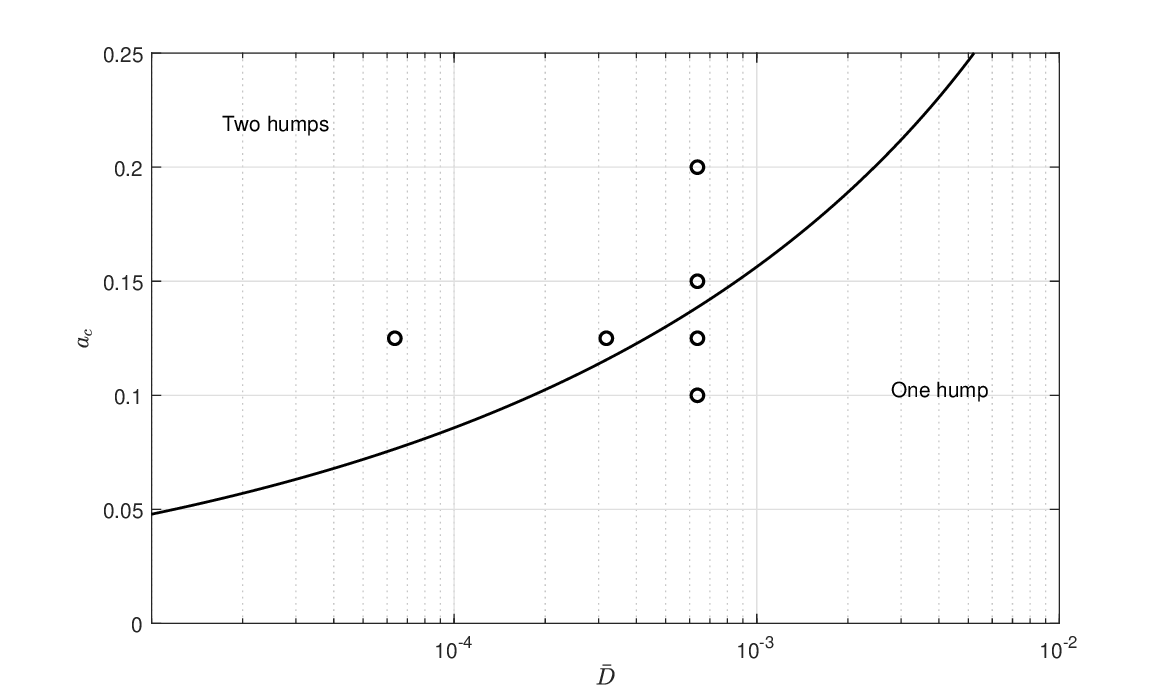}
\caption{A graph of the function $a_c(\bar{D})$, which is defined for $0<\bar{D} \leq \bar{D}^* \approx 5.22 \times 10^{-3}$. The circles show the location of the full numerical solutions shown in Figure~\ref{fig6}.}\label{fig_ac}
\end{center}
\end{figure}

A numerical consideration of equation (\ref{eqn3.140}) gives the value $\overline{D}^* \approx 5.22 \times 10^{-3}$, as shown in Figure~\ref{fig_Dstar},
whilst a numerical solution of equation (\ref{eqn3.141}) determines the graph of $a_c(\overline{D})$ against $\overline{D}$ as shown in Figure~\ref{fig_ac}.
\\

(ii) $\overline{D}\ge\overline{D}^*$

    In this case $\mathcal{F}(x,\lambda,\overline{D},\overline{\phi}_+)$ has single-hump structure for each $a\in \left(0,\frac{1}{4}\right)$.\\\\
This completes the analysis for this kernel perturbation. We have seen that for each $(a,\overline{D}) \in \left(0,\frac{1}{4}\right) \times (0,\infty)$, [NBVP] has a unique solution, and so consequently, via(\ref{eqn3.55}), we conclude that equation (\ref{eqn1.1}) has a unique positive periodic steady state, $\mathcal{F}_p\sim \mathcal{F}(x,\lambda,\overline{D};\overline{\phi}_+)$, for each wavelength $\lambda\in \left(\frac{1}{2},1\right)$ at each $\overline{D}>0$. When $\overline{D}$ is large the structure of $\mathcal{F}$ has the single-hump cosine form at each $a\in\left(0,\frac{1}{4}\right)$, as given in (\ref{eqnf}). This matches  with expansion (\ref{eqn3.4}), to O(1), as $D\to0$, according to the classical Van Dyke asymptotic matching principle. This qualitative single-hump structure is retained at each $\overline{D}\ge \overline{D}^*$. However, for $0<\overline{D}<\overline{D}^*$, there is a bifurcation in structure at $a=a_c(\overline{D})$, with $\mathcal{F}_p\sim \mathcal{F}(x,\lambda,\overline{D};\overline{\phi}_+)$ having single-hump structure for lower values of $a$, but double-hump structure for larger values of $a$. When $\overline{D}$ is small, the structure of $\mathcal{F}_p\sim \mathcal{F}(x,\lambda,\overline{D};\overline{\phi}_+)$ rapidly develops from a cosine single-hump at $a=0^+$, over an interval $a=O(\overline{D}^{\frac{1}{4}})$, to a double-hump structure composed of an approximation with two Airy-type delta functions located at the end points $x=0$ and $x=a$, for $a$ beyond this thin region.

\subsubsection{$\overline{\phi}_N = \overline{\phi}_-$}
In this case, as determined earlier,  there is the possibility of secondary bifurcations occurring in  solutions of equation (\ref{eqn3.89}), and in particular, these will occur across curves on the $(a,\overline{D})$-plane in $(0,1/4) \times \mathbb{R}$ where $\mathcal{J}_-(a,I(a,\overline{D}),\overline{D})$ vanishes. 
However, for each $(a,\overline{D})\in \left(0,\frac{1}{4}\right) \times \mathbb{R}^+$ with $\overline{D}$ sufficiently large, we have seen earlier that equation (\ref{eqn3.89}) has a unique solution, which we label as $I = I^-(a,\overline{D})$, and that this solution is negative, and satisfies the inequality,
\begin{equation}
 -\frac{1}{2}\pi \overline{D}^{-1} <  I^-(a,\overline{D})< -\frac{1}{2}\pi \overline{D}^{-1}(1-16a^2)^{-1}\cos(2\pi a). \label{eqn3.142}    
\end{equation}
Numerical continuation of this branch of solutions determines that the branch continues onto all of $(0,1/4)\times \mathbb{R}^+$, with corresponding values of $\mathcal{J}_-$ being strictly positive throughout, and so, in fact, no secondary bifurcations occur from this primary branch. In addition, it then, again, follows from the Implicit Function Theorem, via (\ref{eqn3.92}) and (\ref{eqn3.96}), that
\begin{equation}
    I^-\in C^{1,1}\left(\left(0,\frac{1}{4}\right) \times \mathbb{R}^+\right),\label{eqn3.143}
\end{equation}
and now
\begin{equation}
    I^-_a(a,\overline{D}),~I^-_{\overline{D}}(a,\overline{D}) > 0   \label{eqn3.144} 
\end{equation}
for each $(a,\overline{D})\in \left(0,\frac{1}{4}\right) \times \mathbb{R}^+$. We conclude that for each $(a,\overline{D})\in \left(0,\frac{1}{4}\right) \times \mathbb{R}^+$, [NBVP] again has a unique solution, which is given, via (\ref{eqn3.81}) and (\ref{eqn3.82}), by 
\begin{equation}
    \mathcal{F}(x,\lambda,\overline{D},\overline{\phi}_-) = \mathcal{F}_L(x,I^-(a,\overline{D}),a)~~\forall~~x\in[-a,a],  \label{eqn3.145} 
\end{equation}
\begin{equation}
    \mathcal{I} = \alpha_L(I^-(a,\overline{D}),a).  \label{eqn3.146}
\end{equation}
To determine the structure of $\mathcal{F}$, we next consider the form of $I^-(a,\overline{D})$ for fixed $\overline{D} > 0$, as $a$ increases through the interval $\left(0,\frac{1}{4}\right)$. First, with $\overline{D}=O(1)^+$, it follows from equation (\ref{eqn3.89}) and (\ref{eqn3.124}), that,
\begin{equation}
 I^-(a,\overline{D})= -\frac{1}{2}\pi \overline{D}^{-1} + O(a^2)~~\text{as}~~a\to 0^+,  \label{eqn3.147}   
\end{equation}
which we observe is in accord with the bounds in (\ref{eqn3.130}). Conversely, as $a\to \frac{1}{4}^-$, we have, from equation (\ref{eqn3.89}) and (\ref{eqn3.125}), that,
\begin{equation}
  I^-(a,\overline{D}) \to I^*(\overline{D}), \label{eqn3.148}
\end{equation}
where now $Y=I^*(\overline{D})$ is the unique negative root of the equation,
\begin{equation}
    \overline{g}_1(Y) + \frac{2}{\pi}\overline{D}Y = 0.
\end{equation}
and it is readily established that $I^*(\overline{D})$ is monotone increasing with $\overline{D}>0$, and strictly larger than $\frac{1}{2}\pi \overline{D}^{-1}$ (in accord with inequalities (\ref{eqn3.142})), whilst,
\begin{equation}
    I^*(\overline{D}) \to 
    \begin{cases}
        0~~\text{as}~~\overline{D}\to \infty,\\
        -\infty~~\text{as}~~\overline{D}\to 0.
    \end{cases}
\end{equation}
Between these limiting forms, $I^-(a,\overline{D})$ is monotone increasing with $a$.  We next examine how the structure of $I^-(a,\overline{D})$ against $a\in \left(0,\frac{1}{4}\right)$ develops as $\overline{D}\to \infty$. This follows directly, via (\ref{eqn3.130}), (\ref{eqn3.79'''}), (\ref{eqn3.89}) and (\ref{eqn3.92}), as,
\begin{equation}
    I^-(a,\overline{D}) \sim -\frac{1}{2}\pi \overline{D}^{-1}(1 - 16a^2)^{-1} \cos(2\pi a)~~\text{as}~~\overline{D}\to \infty,  \label{eqn3.149}
\end{equation}
uniformly for $a\in \left(0,\frac{1}{4}\right)$. We finally establish the structure of $I^-(a,\overline{D})$ against $a\in \left(0,\frac{1}{4}\right)$ as $\overline{D}\to 0^+$. This is obtained directly from (\ref{eqn3.89}), (\ref{eqn3.107}), (\ref{eqn3.124}) and (\ref{eqn3.128}), as,
\begin{equation}
 I^-(a,\overline{D}) \sim -\frac{1}{2}\pi \overline{D}^{-1}~~ \text{as}~~ \overline{D}\to 0^+ 
 \label{eqn3.150}
\end{equation}
uniformly for $a\in \left(0,\frac{1}{4}\right)$. We now use these results to construct $\mathcal{F}$. From (\ref{eqn3.136}), (\ref{eqn3.137}), (\ref{eqn3.132}), and the cases (i)-(iii) earlier, we obtain, after some routine consideration, firstly that, as $\overline{D}\to 0$, with $x\in [0,a]$,
\begin{equation}
    \mathcal{F}(x,\lambda,\overline{D},\overline{\phi}_-) \sim
    \begin{cases}
   \frac{1}{4} a^{-1} \pi \cos\left(\frac{\pi}{2a} x\right),~~a=o(\overline{D}^{\frac{1}{4}})^+,\\
   a^{-1}\hat{f}\left(xa^{-1}, -\frac{1}{2}\pi a^4 \overline{D}^{-1}\right),~~a= O(\overline{D}^{\frac{1}{4}})^+,\\
   \delta_G(x),~~O(\overline{D}^{\frac{1}{4}}) < a \le \frac{1}{4}.
    \end{cases}
\end{equation}
We observe that $\mathcal{F}$ has a single-hump structure at each $a\in \left(0,\frac{1}{4}\right)$. Conversely, as $\overline{D}\to \infty$, with $x\in [0,a]$, we obtain,
\begin{equation}
 \mathcal{F}(x,\lambda,\overline{D},\overline{\phi}_-) \sim
   \frac{1}{4} a^{-1} \pi \cos\left(\frac{\pi}{2a} x\right), \label{eqnf'}
\end{equation}
for each $a\in\left(0,\frac{1}{4}\right)$. Again there is no structure bifurcation and $\mathcal{F}$ has a single-hump structure at each $a\in \left(0,\frac{1}{4}\right)$. Finally, when $\overline{D}=O(1)^+$, we again have
\begin{equation}
 \mathcal{F}(x,\lambda,\overline{D},\overline{\phi}_-) \sim
   \frac{1}{4} a^{-1} \pi \cos\left(\frac{\pi}{2a} x\right)~~\text{as}~~a\to 0^+,
\end{equation}
for all $x\in [0,a]$, with $\mathcal{F}(x,\lambda,\overline{D},\overline{\phi}_-)$ retaining a single-hump structure for each $a\in [0,\frac{1}{4}]$.\\\\
This completes the analysis for this kernel perturbation. We have again seen that for each $(a,\overline{D})\in \left(0,\frac{1}{4}\right) \times (0,\infty)$, [NBVP] has a unique solution, and so consequently, via (\ref{eqn3.55}), we conclude that equation (\ref{eqn1.1}) has a unique positive periodic steady state, $\mathcal{F}_p\sim \mathcal{F}(x,\lambda,\overline{D};\overline{\phi}_-)$, for each wavelength $\lambda\in \left(\frac{1}{2},1\right)$ at each $\overline{D}>0$. When $\overline{D}$ is large the structure of $\mathcal{F}$ has the single-hump cosine form at each $a\in\left(0,\frac{1}{4}\right)$, as given in (\ref{eqnf}). This again matches  with expansion (\ref{eqn3.4}), to O(1), as $D\to0$, according to the classical Van Dyke asymptotic matching principle. This qualitative single-hump structure is, in this case, retained at each $\overline{D}>0$. When $\overline{D}$ is small, the structure of $\mathcal{F}_p\sim \mathcal{F}(x,\lambda,\overline{D};\overline{\phi}_-)$ rapidly develops from a cosine single-hump at $a=0^+$, over an interval $a=O(\overline{D}^{\frac{1}{4}})$, to a single-hump structure composed of a single Gaussian-type delta function located at the origin, for $a$ beyond this thin region. \\

\begin{figure}
\begin{center}
\includegraphics[width=0.9\textwidth]{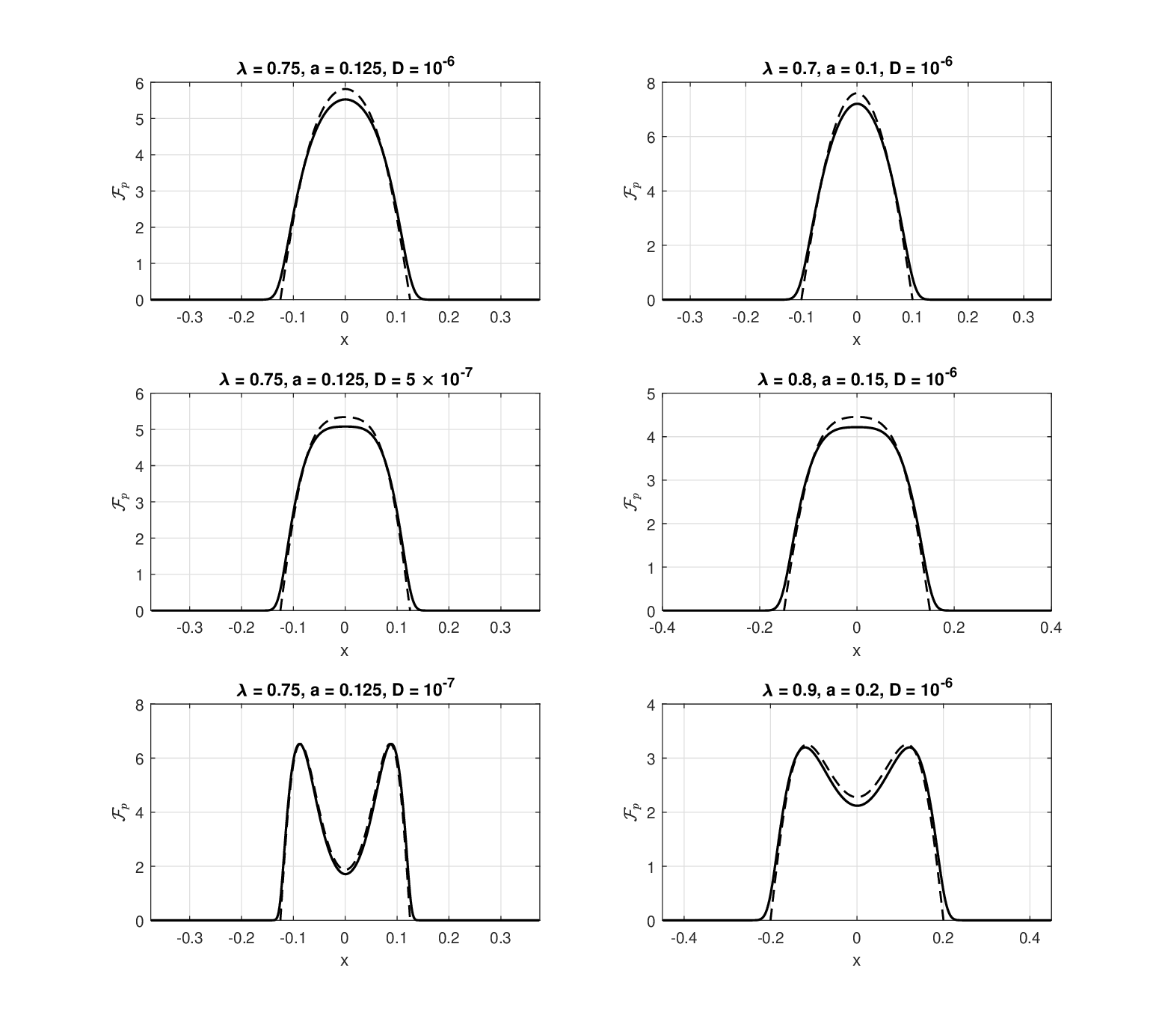}
\caption{Numerical solutions of [FPP] for various values of $D$ and $a$, with $\epsilon =  10^{-3}$ are shown as solid lines. The broken lines are the corresponding solutions of [NBVP] at $\overline{D} = \epsilon^{-1}D$ and the same value of $a$, as in subsections 3.2.2 and 3.2.3. The location of these solutions in the $(\overline{D}, a_c)$-plane is shown in Figure~\ref{fig_ac}.}\label{fig6}
\end{center}
\end{figure}

 To complete this section we determine the periodic steady state $\mathcal{F}_p(x,\lambda,D;\overline{\phi})$, at a number of decreasing values of $D$, and with $(\lambda,D)\in \Omega_1$, via numerical solution of the exact periodic boundary value problem, from (\ref{eqn1.1}),
 \begin{equation}
 D\mathcal{F}_p'' + \mathcal{F}_p\left(1 - \int_{x-\frac{1}{2}}^{x+\frac{1}{2}}{\mathcal{F}_p(y)}dy -\int_{-\infty}^{\infty}{\mathcal{F}_p(y)\overline{\phi}(x-y)}dy \right) = 0,~~x\in \mathbb{R},
 \end{equation}
 subject to,
 \begin{equation}
\mathcal{F}_p \in P_{\lambda}(\mathbb{R})\cap C^2(\mathbb{R}),~~\mathcal{F}_p~~\text{is even}, 
\end{equation}
which we refer to as [FPP]. For comparison with the above theory, we restrict attention to kernel perturbations $\overline{\phi} = \epsilon \overline{\phi}_{\pm}$ (so that $||\overline{\phi}||_1^m = \epsilon$), with $\epsilon = 10^{-3}$ (and so $D= 10^{-3}\overline{D}$). For the $'+'$ kernel perturbation Figure~\ref{fig6} compares numerical solutions of [NBVP] and [FPP] at the six points in parameter space indicated by circles in Figure~\ref{fig_ac}. The transition from one hump to two hump periodic steady states is clearly visible, and consistent, at least to leading order in $\epsilon$, with the leading order location of the boundary, $a_c(\bar{D})$, between these two qualitatively different types of structure. For the $'-'$ kernel perturbation, a similar level of agreement between the the numerical solutions of [NBVP] and [FPP] is observed.

Finally, we return to the comments made at the beginning of this section in Remark 3.1. In particular we now address the region in the structure adopted in equation (\ref{eqn3.53}) where $\mathcal{F}_p(x,\lambda,D;\overline{\phi})$ is exponentially small in $||\overline{\phi}||_1^m$ as $||\overline{\phi}||_1^m \to 0$ (that is $a(\lambda)< |x| \le \frac{1}{2}\lambda$). We can analyse this region via the same approach as that detailed in (NB) (see subsection 4.1), on taking $D=O(||\overline{\phi}||_1^m)$ as $||\overline{\phi}||_1^m \to0$. In doing so (without giving details, which follow directly those in (NM)) we find firstly that a sufficient condition for such a structure to be constructed, and asymptotically matched, so that the form in (\ref{eqn3.53}) continues to hold, is that the corresponding solution to [NBVP] has a single hump structure for $x\in [-a,a]$. Conversely, when a double hump structure becomes sufficiently well developed in the solution to [NBVP], evolving into two separating humps which approach the ends of the core asymptotic region as $\overline{D}$ decreases, the form in (\ref{eqn3.53}) may fail, when it becomes possible for a secondary bifurcation to occur via the emergence, and subsequent growth, of a new hump, localised at the centre of the exponentially small region. Immediate consequences of this analysis in relation to the two representative kernel perturbations studied in this subsection are:

\begin{itemize}
    \item When the kernel perturbation has the negative sign no secondary bifurcations take place from the principal branch of periodic steady states on $(a,\overline{D}) \in (0,1/4) \times \mathbb{R}^+$.
    \item When the kernel perturbation has the positive sign no secondary bifurcations take place from the principal branch of periodic steady states on $(a,\overline{D}) \in ((0,1/4) \times \mathbb{R}^+) \setminus \{(a,\overline{D}):~0<\overline{D}<\overline{D}^*, ~a_c(\overline{D})<a<1/4~\}$. However on $\{(a,\overline{D}):~0<\overline{D}<\overline{D}^*, ~a_c(\overline{D})<a<1/4~\}$ secondary bifurcations from the principal branch become possible, and these are predicated by the asymptotic approximation (\ref{eqn3.53}) to $\mathcal{F}_p$ developing a double hump structure in the core asymptotic region, and are identified by a failure of the asymptotic form for $\mathcal{F}_p$ in (\ref{eqn3.53}) at the centre of the exponentiaaly small asymptotic region, where there is incipient secondary hump formation.
\end{itemize}
The detailed analysis, in the last three subsections, for the two complementary cases of the normalised kernel perturbations $\overline{\phi}_{\pm}$, lead us to suggest and formulate the following more general conjectures (in effect, as extensions of Theorem 3.2), the detailed confirmation of which is beyond the scope of the current paper, but which, nonetheless are worthy of consideration here (it is convenient to adopt the notation of the last three subsections):

\begin{cj}
Let $\overline{\phi}\in \overline{K}(\mathbb{R})$ be an admissible kernel perturbation which has support contained in $[-1/2,1/2]$ and with $||\overline{\phi}||_1^m$ small. Suppose first that the kernel perturbation is of single-hump form, with a positive maximum at the centre and negative towards the ends of the interval. Then we propose that, when $(a,\overline{D}) \in (0,1/4)\times\mathbb{R}^+$:
\begin{itemize}
    \item[C1] The continuation of the positive periodic steady state $\mathcal{F}_p(x,\lambda,D;\overline{\phi})$ from $\Omega^+(\overline{\phi})$ into $\Omega^-(\overline{\phi})$ (which we have denoted by $\mathcal{F}(x,\lambda,\overline{D}; \overline{\phi}_N)$ on the core asymptotic region, at leading order as $||\overline{\phi}||_1^m \to 0$, with $\lambda = 2a +1/2$, and the parameter domain corresponding to $(\lambda,D)\in\Omega^-(\overline{\phi})$ being $(a,\overline{D}) \in (0,1/4)\times\mathbb{R}^+$) provides a unique positive periodic steady state at each $(a,\overline{D})$ with $\overline{D}$ sufficiently large and $a\in (0,1/4)$. This periodic steady state has the structure given in (\ref{eqn3.53}) and (\ref{eqn3.55}), with a symmetric single hump form on the core asymptotic region $[-a,a]$. However it is anticipated that, with decreasing $\overline{D}$, and in particular for $a$ at the upper end on the interval $(0,1/4)$, the concentration of the perturbed kernel towards the centre of its support interval and away from the edges, has the effect of driving the structure of the periodic steady state in the core asymptotic region to develop into a symmetric two-hump form. When these two symmetric humps begin to separate and approach the edges of the core region (as they do with decreasing $\overline{D}$), their presence induces the breakdown of the structure in the exponentially small region, leading to the emergence of a new localised hump. This appears at the centre of the exponentially small region, and triggers a local secondary bifurcation as a fold in the primary branch.This type of process in studied further in the next section via a detailed numerical investigation of [FPP].

    Next suppose that the kernel perturbation is of an inverted single hump form, with a negative minimum at the centre and positive towards the ends of the interval. Then we next propose that, when $(a,\overline{D}) \in (0,1/4)\times\mathbb{R}^+$:
    \item[C2] Again the continuation of the positive periodic steady state $\mathcal{F}_p(x,\lambda,D;\overline{\phi})$ from $\Omega^+(\overline{\phi})$ into $\Omega^-(\overline{\phi})$  provides a unique positive periodic steady state at each $(a,\overline{D})$ with $\overline{D}$ sufficiently large and $a\in (0,1/4)$. This periodic steady state has the structure given in (\ref{eqn3.53}) and (\ref{eqn3.55}), with a symmetric single hump form on the core asymptotic region $[-a,a]$. However, in this case, the concentration of the perturbed kernel towards the edges of its support interval and away from the centre, now has the effect of driving the structure of the periodic steady state in the core asymptotic region to retain its single hump form, and focus towards a single spike at the centre of the core asymptotic region as $\overline{D}$ decreases. This process simply reinforces the structure in the exponentially small region, and as such inhibits the occurrence of secondary bifurcations from the principal branch.
\end{itemize}
\end{cj}

We now move on to consider these conjectures via careful numerical solution of the full boundary value problem [FPP].

\section{Secondary bifurcation structure when $\overline{D}=O(1)$ for periodic steady states with close to unit wavelength}
The key point to recall (as discussed in some detail in the last section) at the beginning of this section is that \emph{all of the results presented in subsection~\ref{subsec_smallD} are based upon the asymptotic approximation (\ref{eqn3.53}) for $F_p$ with $D = O(||\overline{\phi}||_1^m)$ as $||\overline{\phi}||_1^m\to 0$} ($\overline{D}=O(1)$), and this requires, following the analysis in (NB) (section 4), that the local bifurcation (from the equilibrium state $u=1$), which generates these periodic steady states, takes place, formally, when $D=O(1)$. This immediately places the restriction
\begin{equation}
\lambda\in \left( \frac{1}{2}+ O(\sqrt{D}), 1-O(D) \right)
\end{equation}
on the wavelength when $D$ becomes small. Consequently the results of subsection~\ref{subsec_smallD} require the formal restriction,
 \begin{equation}
 \lambda \in \left(\frac{1}{2} + O(\sqrt{||\overline{\phi}||_1^m} ), 1-O(||\overline{\phi}||_1^m)\right), \label{eqn''}
 \end{equation}
 throughout, and equivalently,
 \begin{equation}
     a(\lambda) \in 
 \left( O(\sqrt{||\overline{\phi}||_1^m }), \frac{1}{4}-O(||\overline{\phi}||_1^m)\right),
 \end{equation}
with $D = O(||\overline{\phi}||_1^m)$ as $||\overline{\phi}||_1^m\to 0$. In this section our objective is to complement the analysis of subsection~\ref{subsec_smallD} by considering periodic steady states, via careful direct numerical investigation of [FPP], in the remaining regions of $\Omega_1$, were $\lambda \sim \frac{1}{2} + O(\sqrt{||\overline{\phi}||_1^m})$ and $\lambda \sim 1-O(||\overline{\phi}||_1^m)$ when $D= O(||\overline{\phi}||_1^m)$, which are outside the range of asymptotic uniformity of the theory developed in the previous section, and similarly consider the possibility of secondary bifurcations when $D=o(||\overline{\phi}||_1^m)$ with  $\lambda \in \left(\frac{1}{2} + O(\sqrt{||\overline{\phi}||_1^m} ), 1-O(||\overline{\phi}||_1^m)\right)$. We focus on the kernel perturbation (\ref{eqn3.85}), taking the positive sign, which is the case for which two hump structures emerge in the family of periodic steady states identified and analysed in subsection~\ref{subsec_smallD} (under the restriction (\ref{eqn''}), and shown in Figure~\ref{fig6}. Our objective is to uncover support for Conjecture 3.1. First we find that when $\lambda \sim \frac{1}{2} + O(\sqrt{||\overline{\phi}||_1^m})$, the situation remains very close to that of subsection~\ref{subsec_smallD} and needs no further discussion. However the situation is quite different to that determined in subsection~\ref{subsec_smallD} when $\lambda \sim 1-O(||\overline{\phi}||_1^m)$, with an interesting secondary bifurcation structure emerging when $D<<||\overline{\phi}||_1^m)$, in accord with Conjecture 3.1.

Following the above discussion we now analyse periodic steady states of [FPP] with the magnitude of the kernel perturbation $||\overline{\phi}||_1^m = \epsilon$ fixed as positive and small, the wavelength $\lambda$ fixed according to our intended investigation, whilst $\overline{D} = \epsilon^{-1}D$ is decreased from its local bifurcation value on the boundary of $\Omega_1$ down towards zero, and is regarded as a bifurcation parameter. The primary branch of periodic steady states emerging from this local bifurcation, as identified in the preliminary part of section 3, is then numerically path followed as $\overline{D}$ is decreased. The periodic steady state on this branch, as $\overline{D}$ decreases initially (after the very early weakly nonlinear phase), develops the 'single humps separated by exponentially small regions' structure and then at a lower value of $\overline{D}$, within the hump region, the single hump develops two symmetric peaks from the single peak previously situated at the origin (this is a smooth development without formal bifurcation). For simpliciy of notation, we will henceforth refer to this primary branch of periodic steady states as the \emph{one-peak branch} (bearing in mind that, as noted above, two peaks do develop whilst on this branch). As we shall see a secondary bifurcation occurs from this branch, and this secondary bifurcation structure, as $\overline{D}$ is decreased further, is then strongly influenced by the behaviour in the exponentially small part of the periodic steady state. In consequence, to retain careful numerical accuracy, we define $W = \log u$ and write the steady version of (\ref{eqn1.1}) as
\begin{equation}
    D\left(W^{\prime \prime} + W^{{\prime}^2}\right) + 1 - \phi * e^{W}=0,\label{bif_eqn}
\end{equation}
where a prime denotes $d/dx$ and 
\begin{equation}
    \phi(y) = \left(1+ \overline{\epsilon}\cos 2 \pi y\right) H\left(\frac{1}{4} - y^2\right).\label{bif_kernel}
\end{equation}
where, for convenience, we have written $\overline{\epsilon} = \frac{1}{2} \pi \epsilon$. This allows us to compute solutions with $W$ large and negative, for which a numerical method based on the use of (\ref{eqn1.1}) to find $u$ would fail to accurately resolve the key features.

We compute numerical solutions of (\ref{bif_eqn}) on a domain of wavelength $\lambda$ with periodic boundary conditions. We discretise at equally spaced points, evaluate derivatives using five point central differences and calculate the convolution using fast Fourier transforms. The resulting algebraic equations are solved in MATLAB using {\em fsolve}, which is an implementation of the trust region dogleg method, with an analytical Jacobian provided to speed up the calculation. We also use deflation to avoid the numerical solution converging to the steady state $W=0$, by dividing the whole system by $\tanh(|{\bf W}|^2/N)$, where ${\bf W}$ is the vector of unknowns and $N$ the number of unknowns. 

\begin{figure}
\begin{center}
\includegraphics[width=0.65\textwidth]{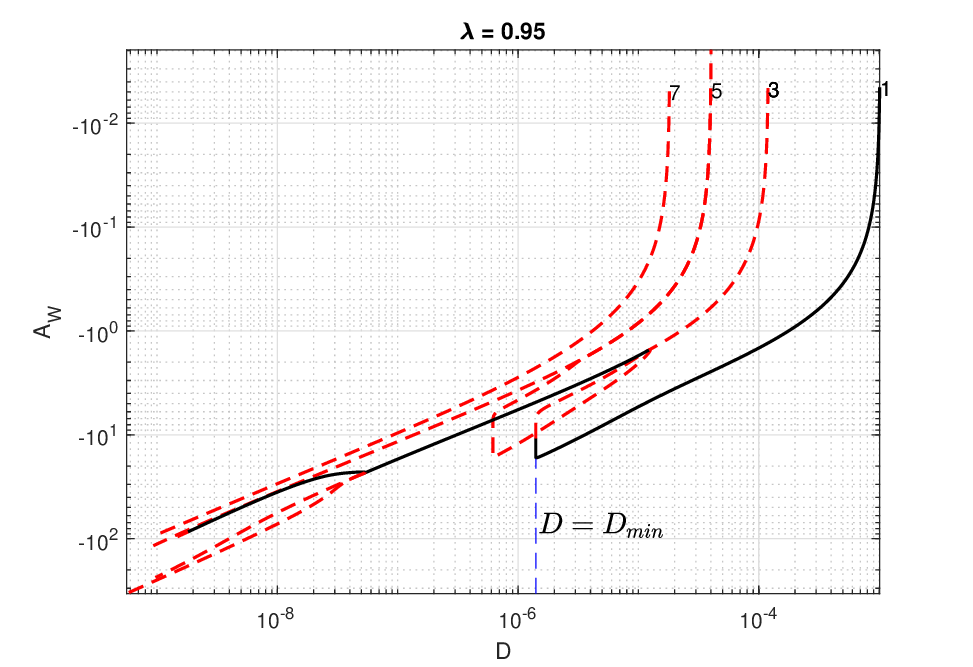}
\includegraphics[width=0.65\textwidth]{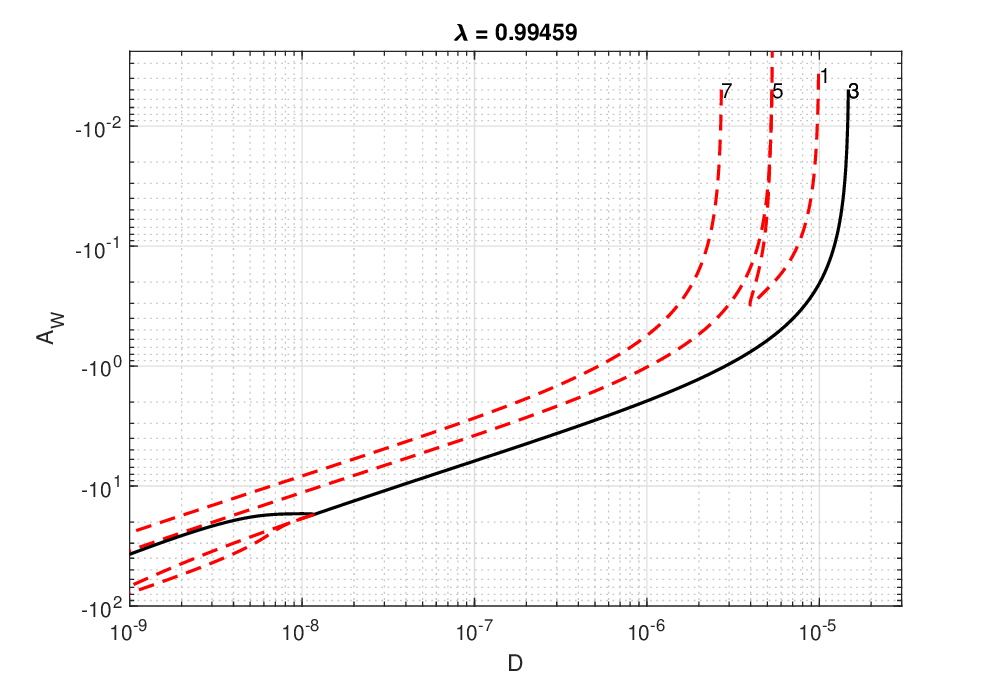}
\caption{The bifurcation diagram for [FPP] with $\overline{\epsilon} = 0.01$ and $\lambda = 0.95$ (upper panel) and $\lambda = 0.99459$ (lower panel). Note that for this value of $\overline{\epsilon}$ the largest wavelength on $\Omega_1$ is given by $\lambda_{\rm max}(\overline{\epsilon}) = \sqrt{1 - \overline{\epsilon}} \approx 0.995$. The red broken lines denote unstable periodic steady states and the solid black lines stable periodic steady states. The numbers indicate the number of peaks over one wavelength $\lambda$ in the periodic steady state in the neighbourhood of the local bifurcation curve at the boundary of $\Omega_1$. The value of $D_{\rm min}$, at which the one-peak bifurcation curve first loops around, is also indicated (see Figure~\ref{fig_Dmin}).}\label{fig_bif}
\end{center}
\end{figure}

For reasons that will become clear below, we compute the resulting bifurcation diagram using adaptive arc-length continuation in the $(\log_{10} D, \log_{10}(-A_W))$-plane, where
\begin{equation}
A_W = \int_0^{\lambda} W(x)\,dx.
\end{equation}
We use the small-amplitude periodic steady state, given by a weakly nonlinear analysis of the local steady state supercritical pitchfork bifurcation (from the equilibrium state $u=1$) as $\overline{D}$ passes into $\Omega_1$ through its boundary (see Appendix~\ref{appendix}), to provide an accurate initial guess to start the computation of the bifurcation curve. Figure~\ref{fig_bif} shows the bifurcation diagram for two values of the wavelength, $\lambda$, with fixed $\overline{\epsilon} = 0.01$ (for convenience we use $D$ rather than $\overline{D}$ in the Figures of this section).
Much more could be said about how this structure changes as $\overline{\epsilon}$ increases, but the focus of this paper is the robustness of the top hat kernel to small perturbations. 

\begin{figure}
    \centering
    \includegraphics[width = 0.8\textwidth]{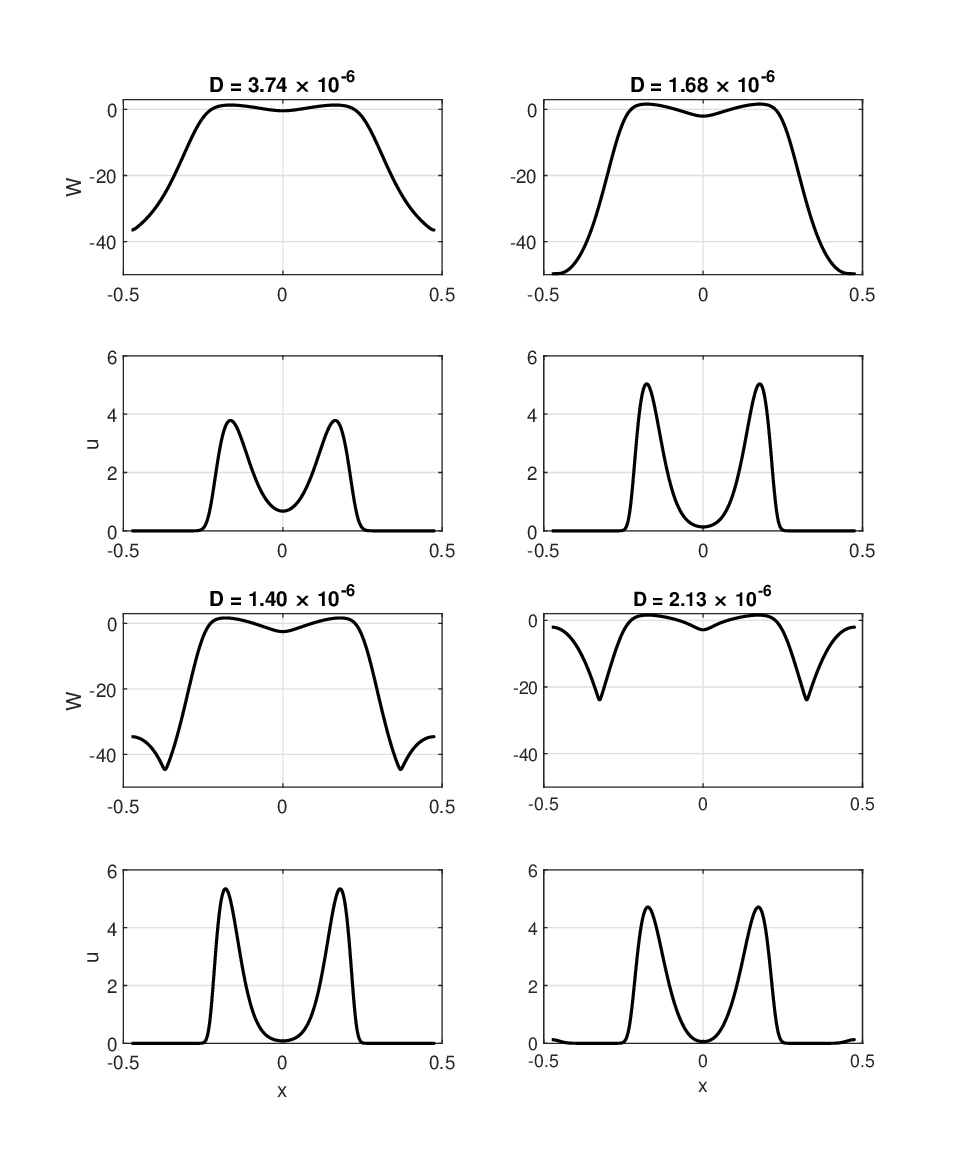}
    \caption{Four successive periodic steady states on the one-peak bifurcation curve when $\overline{\epsilon} = 0.01$ and $\lambda = 0.95$ (see Figure~\ref{fig_bif}, upper panel). These plots indicate the qualitative changes in the periodic steady state as the first loop of the one-peak bifurcation curve is traversed.}
    \label{fig_bif3}
\end{figure}

Each panel of Figure~\ref{fig_bif} shows branches initiated from the weakly nonlinear theory of Appendix~\ref{appendix} with small amplitude (relative to $u=1$) periodic steady states of wavelength $\lambda$, but with $1$, $3$, $5$ or $7$ peaks, or equivalently, the weakly nonlinear steady states with wavelengths $\lambda$, $\lambda/3$, $\lambda/5$ and $\lambda/7$ on the periodic domain of length $\lambda$. Note that periodic steady states do not exist for wavelengths that are even fractions of $\lambda$, since they lie outside the boundaries of the tongues $\Omega_i$ in the $(\lambda, D)$-plane, as discussed at the start of section 3. We computed these branches because, as can be seen in the upper panel of Figure~\ref{fig_bif}, the one-peak branch turns close to vertically at $D \approx 1.4 \times 10^{-6}$ and then touches the three-peak branch at $D \approx 1.5 \times 10^{-5}$ before turning around again and terminating on the five-peak branch. As discussed towards the end of section 3, the mechanism by which these peaks are created originates in the behaviour of the periodic steady state in the region separating each core region double hump form, where $u$ is exponentially small. Briefly, to understand this mechanism, we see that in the regions between the humps, where $u$ is exponentially small, $W$ is large and negative. If we then rescale using $W = D^{-1/2} \bar{W}$, with $D=O(\overline{\epsilon})$, we obtain from (\ref{bif_eqn}) the leading order equation
\begin{equation}
    \bar{W}^{{\prime}^2} = \phi*e^W - 1, \label{eqn_expsmall}
\end{equation}
which dictates the structure in the region where $u$ is exponentially small. The nonlocal nature of the equation is crucial here, as the convolution on the right hand side of (\ref{eqn_expsmall}) is determined from the periodic steady state in the hump region, where $e^W = u = O(1)$ when $D=O(\overline{\epsilon})$ is small. For a solution to (\ref{eqn_expsmall}) to exist in the exponentially small region, we clearly require $C(W) = \phi*e^W-1> 0$ when evaluated in that region. However, investigating our numerical solution of [FPP] along the lower part of the one peak branch in the upper panel of Figure~\ref{fig_bif}, the profile of $C(W)$ across the exponentially small region is initially everywhere positive. However, it then decreases with decreasing $D$, and eventually reaches zero at the  symmetry point in the centre of the exponentially small region, when $D \approx 1.4 \times 10^{-6}$, the value where the bifurcation curve turns vertically, and thereafter initiates a neighbourhood of the symmetry point where it is positive - at this stage the approximation (\ref{eqn_expsmall}) associated with the exponentially small region fails, and we see this in the solution to [FPP] as the growth of a new incipient peak at the centre of the exponentially small region. The sequence of periodic steady states shown in Figure~\ref{fig_bif3} illustrates this phenomenon (it should be noted that in these Figures the exponentially small region is approximately$(-0.475,-0.225)\cup (0.225,0.475)$, with the 'centre points' of the exponentially small region at locations $x = \pm 0.475$). It can also be seen in the Figure how this phenomenon is driven by the separation of the two humps and their drift towards the edges of the core region. 
The sequence of changes along the bifurcation curve is most easily understood by watching the video available \href{https://drive.google.com/file/d/1xfcDbtDMlL9VEb7jke-m3nuHfid1DLZw/view?usp=sharing}{here}, but we will also describe these changes.

As $D$ initially decreases, the profile values of $W$ in the exponentially small region between the humps decrease, with the slope at the centre points of the exponentially small region ($x=\pm 0.475$ in the Figures) approaching zero (the upper four panels of Figure~\ref{fig_bif3}). As we continue along the bifurcation curve, the region between the humps where $u$ is exponentially small develops a local maximum in $W$ (at $x=\pm 0.475$ in the Figures), which increases as the bifurcation curve turns around (bottom four panels). Eventually, $W$ approaches zero at this maximum ($x=\pm 0.475$), and a new hump in $u$ starts to form, which can be seen in the bottom right panel, which is the genesis of a three hump periodic steady state. We note that, regarding the $W$-profiles, whilst the exponentially-small part of the steady state changes dramatically, the plot of $u$ shows that the $O(1)$ parts of the periodic steady state change very little. We also observe that the periodic steady state loses stability roughly at the profile given at the bottom left of Figure~\ref{fig_bif3}. It should now be clear why we use $A_W$ as the bifurcation parameter. We need a quantity that measures the difference between the periodic steady states shown in Figure~\ref{fig_bif3}, which is a consequence of $O(1)$ changes in $W = \log u$, not $u$.

\begin{figure}
    \centering
    \includegraphics[width = 0.8\textwidth]{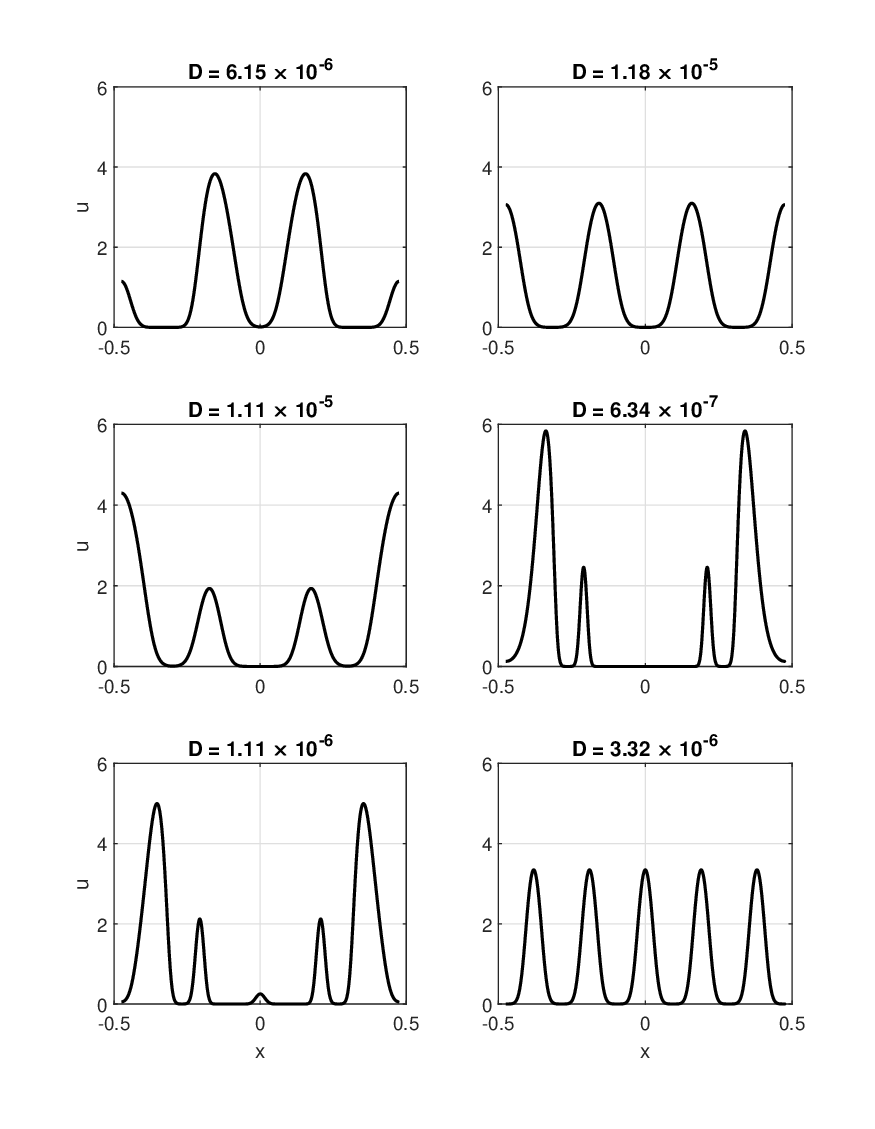}
    \caption{Six successive periodic steady states on the one-peak bifurcation curve when $\overline{\epsilon} = 0.01$ and $\lambda = 0.95$ as it curves around to first touch the three-peak curve and then terminate on the five-peak curve (see Figure~\ref{fig_bif}, upper panel). This sequence follows on from the sequence shown in Figure~\ref{fig_bif3}, and indicates the qualitative changes in the steady state as the second and third loops of the bifurcation curve are traversed.}
    \label{fig_bif4}
\end{figure}

The three-peak steady state formed in this process touches the three-peak branch in a further secondary bifurcation at $D\approx 1.5 \times 10^{-5}$, at which point the three humps are of equal height. The branch then doubles back with one hump growing and two shrinking, as shown in Figure~\ref{fig_bif4}, and the process shown in Figure~\ref{fig_bif3} repeats itself at the large hump, which splits in two, leading to the formation of a five-peak periodic steady state. The bifurcation curve then terminates at the five-peak branch.
This process transforms a group of three humps into five humps. There does not appear to be an equivalent mechanism to transform other combinations of humps, so there is no connection between the five- and seven-peak branches.

\begin{figure}
    \centering
    \includegraphics[width = 0.9\textwidth]{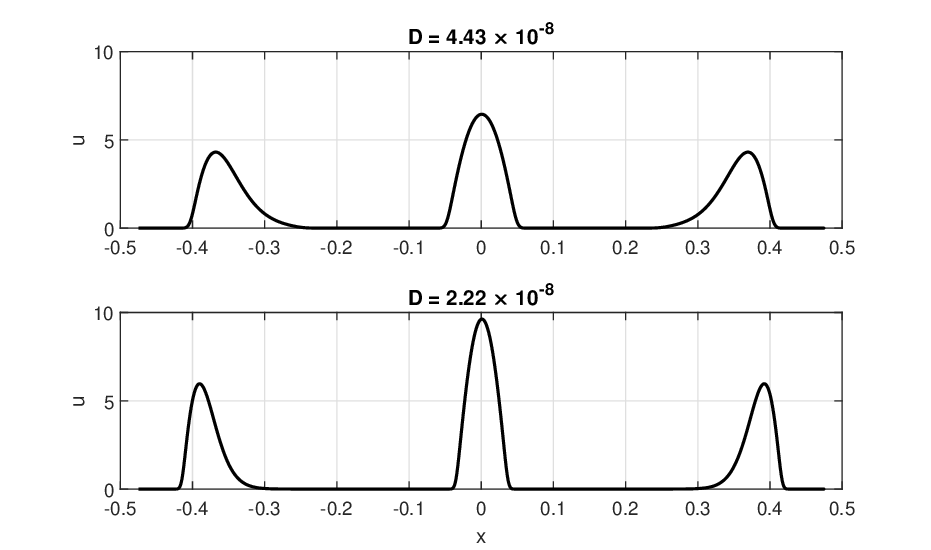}
    \caption{Two typical, asymmetric, three-peak periodic steady states on the branch that bifurcates from the three-peak branch at $D \approx 6 \times 10^{-8}$ when $\overline{\epsilon}= 0.01$ and $\lambda = 0.95$.}
    \label{fig_bif5}
\end{figure}

The other bifurcations of interest lie on the three-peak branch which, whilst initially unstable, stabilises when touched by the originally one-peak branch, then loses stability to an asymmetric bifurcation at $D \approx 6 \times 10^{-8}$. Two typical periodic steady states on this asymmetric branch are shown in Figure~\ref{fig_bif5}.
This bifurcation structure leads to a range of values of $D$ for which two stable steady states coexist, one of wavelength $\lambda$ and one of wavelength $\lambda/3$. For the case shown in the upper panel this range is approximately $1.5 \times 10^{-6} < D < 1.5 \times 10^{-5}$. There is a similar, but very small, region where there are two stable steady states close to the bifurcation point at $D \approx 6 \times 10^{-8}$.

\begin{figure}
    \centering
    \includegraphics[width = 0.9 \textwidth]{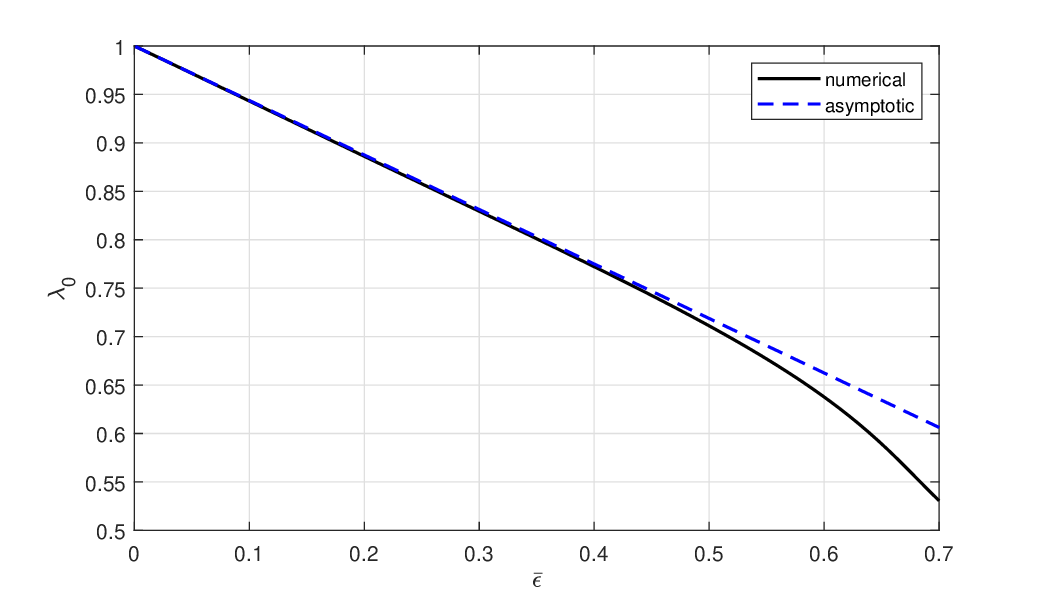}
    \caption{The value of $\lambda$ at which the qualitative nature of the bifurcation diagram changes. This is the value of $\lambda$ at which a periodic steady state of wavelength $\lambda$ and a periodic steady state of wavelength $\lambda/3$ emerge at the primary steady state pitchfork bifurcation at the same value of $D$ on the boundary of $\Omega$ (on entering $\Omega_1$ and $\Omega_2$ respectively).}
    \label{fig_lambda0}
\end{figure}

For values of $\lambda$ close enough to $\lambda_{\rm max}(\overline{\epsilon})$ (the maximum value of $\lambda$ in $\Omega_1$) so that the one-peak periodic steady state branch meets the boundary of $\Omega_1$ at a smaller value of $D$ than the three-peak periodic steady state branch meets the boundary of $\Omega_2$, an example of which is shown in the lower panel of Figure~\ref{fig_bif}, only the three-peak periodic steady state is stable, until an asymmetric bifurcation, at $D \approx 1.5 \times 10^{-8}$ for this value of $\lambda$. Recall that we use 'one-peak periodic steady state' to refer to a periodic steady state on the branch that originates at the boundary of $\Omega_1$ with one peak. As we have seen, these periodic steady states actually develop two peaks for sufficiently small $D$. The (unstable) one-peak branch connects to a different branch. In the lower panel of Figure~\ref{fig_bif} this is the five-peak branch. Numerical exploration provides strong evidence that there is a value $\lambda_0(\overline{\epsilon})$ such that this is the generic picture for $\lambda_0(\overline{\epsilon})<\lambda<\lambda_{\rm max}(\overline{\epsilon})\sim 1-\frac{1}{2}\overline{\epsilon}$. It is straightforward to show that $\lambda_0(\overline{\epsilon}) \sim 1 - \frac{9}{16} \overline{\epsilon}$ as $\overline{\epsilon} \to 0$, and this gives an excellent approximation for $0 \leq \overline{\epsilon}< 0.5$, as shown in Figure~\ref{fig_lambda0}.
In other words, the one-peak branch is always unstable for $\lambda_0(\overline{\epsilon})< \lambda < \lambda_{\rm max}(\overline{\epsilon})$. 

\begin{figure}
    \centering
    \includegraphics[width = 0.9\textwidth]{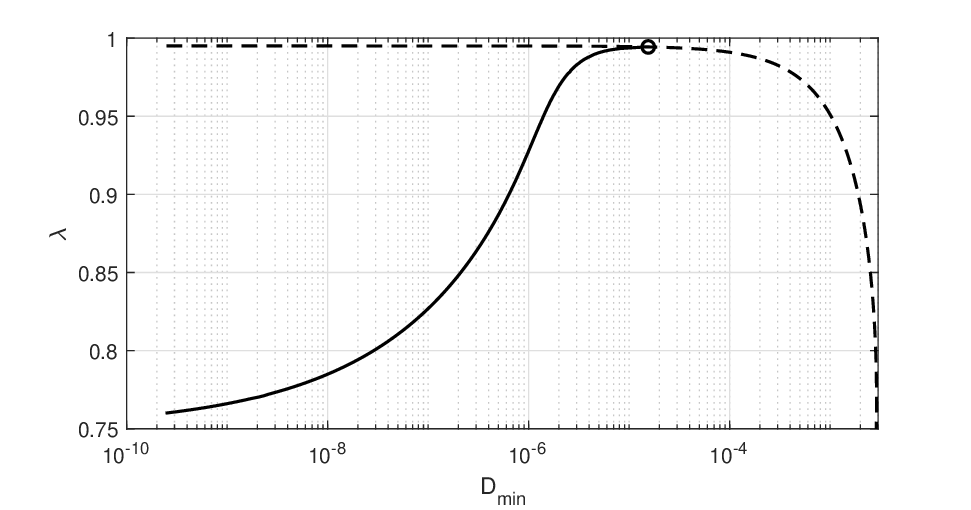}
    \caption{The numerically-calculated value of $D_{\rm min}$, the smallest value of $D$ for which a one-peak periodic steady state exists for $\overline{\epsilon}= 0.01$ (see the upper panel of Figure~\ref{fig_bif} for an example of where this lies on the bifurcation curve when $\lambda = 0.95$). Note that $\lambda_0 \approx 0.9944$, which meets the boundary (neutral) curve of $\Omega_1$ (shown as a broken line) at $D \approx 1.55 \times 10^{-5}$, indicated on the graph with a circle. Note that $D_{\rm min}(\lambda,\overline{\epsilon})\to 0$ as $\overline{\epsilon} \to 0$ uniformly for $\lambda \in [\frac{3}{4},\lambda_0(\overline{\epsilon})]$.}
    \label{fig_Dmin}
\end{figure}

We now return to our discussion of the exponentially small part of the principal one peak-periodic steady state, in particular (\ref{eqn_expsmall}). As we have seen, when $D$ becomes $O(\overline{\epsilon})$ the one-peak periodic steady state develops two peaks on its effective compactly supported hump region, via local hump splitting at the origin, which, when $D\ll\overline{\epsilon}$, has developed into the form of two delta functions, each of weight $\frac{1}{2}$, located at $x = \pm \frac{1}{2}\left(\lambda-\frac{1}{2}\right)$. Now, for this structure to then admit development into three or more peaks, its wavelength must be short enough for the unit width of the kernel to span at least three of these delta functions at all points where $u$ is exponentially small, so that the right hand side of (\ref{eqn_expsmall}) is positive, as is required for a solution of (\ref{eqn_expsmall}) to exist. The marginal case for this condition to be satisfied is when the wavelength $\lambda$ is such that a unit span from the midpoint of the exponentially small region encompasses precisely four of the spikes, so that $2\left(\lambda-\frac{1}{2}\right) + \frac{1}{2} = 1$, and hence $\lambda = \frac{3}{4}$. Thus, for $\frac{3}{4} < \lambda < \lambda_0(\overline{\epsilon})$ it is possible for the one-peak periodic steady state  to ultimately undergo the transition to a three-peak periodic steady state, with the one-peak periodic steady state existing only until $D$ reaches a critical value, which has $D\le O(\overline{\epsilon})$. However, for $\frac{1}{2} < \lambda < \frac{3}{4}$, this transition is ruled out, and the (stable) one-peak periodic steady state persists with decreasing $D= o(\overline{\epsilon})$, and although it does eventually undergo the hump splitting process in the core region, it does not undergo any secondary bifurcations. It is worth noting, in additional support of this argument, that for $\frac{1}{2} < \lambda < \frac{3}{4}$, the value of $\lambda/3$ does not lie within the second tongue $\Omega_2$, and so no three-peak periodic steady state is available. To numerically verify this condition, we computed the smallest value of $D = D_{\rm min}(\lambda, \overline{\epsilon})$ for which the one-peak branch exists as a function of wavelength $\lambda$ when $\overline{\epsilon} = 0.01$, i.e. the point at which the one-peak curve turns vertically in Figure~\ref{fig_bif}, upper panel, and this is shown in Figure~\ref{fig_Dmin}.
We can see that this curve terminates at $\lambda = \lambda_0(\overline{\epsilon})$ on the neutral curve and is consistent with $D_{\rm min}(\lambda,\overline{\epsilon}) \to 0 $ as $\lambda \to \frac{3}{4}$ when $\overline{\epsilon}>0$. It is difficult to calculate accurate numerical solutions in a reasonable time for $D$ smaller than about $10^{-10}$. Note that $D_{\rm min}(\lambda,\overline{\epsilon}) = O(\overline{\epsilon})$ as $\overline{\epsilon} \to 0$ uniformly for $\lambda\in [\frac{3}{4},\lambda_0(\overline{\epsilon})]$. 

We emphasise that when $\overline{\epsilon} = 0$, the unperturbed top hat kernel, the one-peak branch is stable, there is no hump splitting in the core region of the periodic steady state, all the other branches are unstable and there are no secondary bifurcations from the primary one-peak periodic steady state branch. Indeed we recall that all of the discussion above is of features that arise because of the \emph{small but nontrivial perturbation} of the kernel.

A similar investigation has been undertaken for the kernel perturbation (\ref{eqn3.85}) with now the negative sign taken. We will not give the details here, but a similar argument to that given above indicates that the generation of a new peak in the exponentially small region can only occur for wavelengths that satisfy $1 < \lambda <\sqrt{1+\overline{\epsilon}}$. Numerical solutions then indicate that the one-peak branch connects with the two-peak branch, which lies within the first tongue, $\Omega_1$, for this range of wavelengths, which is small when $\overline{\epsilon}$ is small. However, for $1/2<\lambda<1$ the one-peak branch does not undergo any secondary bifurcations or hump splitting process.\\

To finish this section we observe that all of the above results and observations are in agreement with the theory of section 3, provide support for Conjectures 3.1, and give additional specific information, regarding the nature and location of secondary bifurcations from the singularly perturbed principal branch of periodic steady states in the parameter region $\Omega^-(\overline{\phi})\subset \Omega_1$. It is worth reviewing what we have uncovered in this section and in doing so we will rephrase the results in terms of the parameters $(a,\overline{D})$ rather than $(\lambda, D)$ - this enables us to integrate these discoveries with the conclusions of section 3 in a consistent way. When $\overline{\epsilon}$ is fixed and small, we recall that, at leading order in $\overline{\epsilon}$, the parameter region we are confined to is again $(a,\overline{D}) \in (0,1/4)\times\mathbb{R}^+$. On including the correction at $O(\overline{\epsilon})$ the upper boundary of this parameter region is adjusted to 
\begin{equation}
    a = a_R(\overline{D},\overline{\epsilon})\equiv \frac{1}{4} - \frac{1}{4}(16 \pi \overline{D} \pm 1)\overline{\epsilon} + O(\overline{\epsilon}^2)~~\text{as}~~\overline{\epsilon}\to 0,
\end{equation}
with $0<\overline{D}\le O(1)$, with the sign chosen according to the plus or minus kernel perturbation respectively.
This adjusted parameter region is now labeled as $\Gamma(\overline{\epsilon})$. It is also convenient to represent the hump splitting curve $a=a_c(\overline{D})$ (determined is subsection 2.2.3) in inverse form as $\overline{D} = \overline{D}_c(a)$ for $a\in (0,1/4)$, and the location of the secondary folding bifurcation from the primary one peak branch as $\overline{D}_{min}(a,\overline{\epsilon}) \equiv \overline{\epsilon}^{-1} D_{min}(\lambda(a),\overline{\epsilon}) = O(1)$ as $\overline{\epsilon}\to 0$ for $a\in (1/8,(1/4)-(9/32)\overline{\epsilon})$, with $\overline{D}(1/8,\overline{\epsilon})=0$ and $\overline{D}((1/4)-(9/32)\overline{\epsilon},\overline{\epsilon}) =1/128\pi$. We can now reinterpret the above conclusions in this parameter region which, after a little calculation, leads to the following:
\begin{summary}
    Fix $\epsilon>0$ small. To best illustrate the bifurcation structure, we then fix $a \in (0, a_R(0,\overline{\epsilon}))$ and allow $\overline{D}$ to decrease from a sufficiently large value (where only the principal 1-peak branch is present) down to zero. Then for the positive sign kernel perturbation $\overline{\phi}_+$ we have the following summary of the secondary bifurcation structure in the parameter region $(a,\overline{D}) \in \Gamma(\overline{\epsilon})$:
    \begin{itemize}
    \item[S1]First when $a\in (0,1/8]$ there are no secondary bifurcations from the principle branch, although it passes through the hump splitting curve when $\overline{D} = \overline{D}_c(a)$, at which stage the single peak splits in two at the centre of the periodic steady state profile and its profile eventually approaches a form with simply two Airy-type delta functions located at the ends of the core region as $\overline{D} \to 0$. However, when $a\in (1/8,(1/4)-(9/32)\overline{\epsilon}]$ the principal branch is terminated at a secondary fold bifurcation when $\overline{D} = \overline{D}_{min}(a,\overline{\epsilon})~(<\overline{D}_c(a))$, after which multiple  subsequent bifurcations occur from the folded branch, with the bifurcation diagram now having the qualitative form of that represented in the upper panel of Figure~\ref{fig_bif}. Last,in a thin upper edge region when $a\in ((1/4)-(9/32)\overline{\epsilon},a_R(0,\overline{\epsilon}))$, the principal branch is again terminated at a secondary fold bifurcation when $\overline{D} = \overline{D}_{min}(a,\overline{\epsilon})~(<\overline{D}_c(a))$, after which multiple  subsequent bifurcations occur from the folded branch, but the bifurcation diagram now has the qualitative form of that represented in the lower panel of Figure~\ref{fig_bif}.

    Next for the negative sign kernel perturbation $\overline{\phi}_-$ we have the following summary of the secondary bifurcation structure in the parameter region $(a,\overline{D}) \in \Gamma(\overline{\epsilon})$:
    
    \item[S2] When $a\in(0,1/4]$ there are no secondary bifurcations from the principle branch, and in the core region the profile of the periodic steady state remains as a single hump, and approaches a single, central Gaussian-type delta function as $\overline{D} \to 0$. In a thin upper edge region when $a\in (1/4,a_R(0,\overline{\epsilon}))$ the principal branch undergoes a folding bifurcation again, much as in the similar case in the thin edge region above.
    \end{itemize}
\end{summary}
This completes our numerical investigation of positive periodic steady states. We are now in a position to use the theory in sections 1-4 to draw further conclusions concerning the primary conjecture (P2) relating to the evolution problem (IBVP)$_p$. Before doing this, it is instructive to consider some numerical direct numerical solutions to (IBVP)$_p$.
 
\section{Numerical solution of (IBVP)$_p$}
In this section we discuss some numerical solutions of the full evolution problem(IBVP)$_p$ with kernel perturbations again given given by $\overline{\phi} = \epsilon\overline{\phi}_{\pm}$, for a range of values of $\overline{\epsilon} = \frac{1}{2} \pi \epsilon$. In particular, we investigate the principal conjectures (P1) and (P2) in relation to the current detailed  theory developed in section 2 and section 3, concerning equilibrium state temporal stability, and positive periodic steady states. Numerical simulations for the top hat kernel ($\overline{\epsilon} = 0$), as discussed in \cite{NBLM} (section 2), show that localised initial conditions generate a pair of diverging travelling wavefronts that leave behind a stable and stationary periodic steady state with wavelength that approaches $\frac{1}{2}$ as $D \to 0$. This occurs through a mechanism that initiates incipient spikes ahead of the wavefront, where $u$ is exponentially small, and which mature on passing through the wavefront and organise into the periodic steady state with wavelength close to $\frac{1}{2}$ at the rear of the wavefront. For the perturbed top hat kernel, when $|\overline{\epsilon}|$ is sufficiently small, and $D=O(1)$, numerical solutions confirm that the evolution is simply a regular perturbation on that when $\overline{\epsilon} =0$, which is anticipated from the theory in the earlier sections. 

\begin{figure}
\begin{center}
\includegraphics[width=\textwidth]{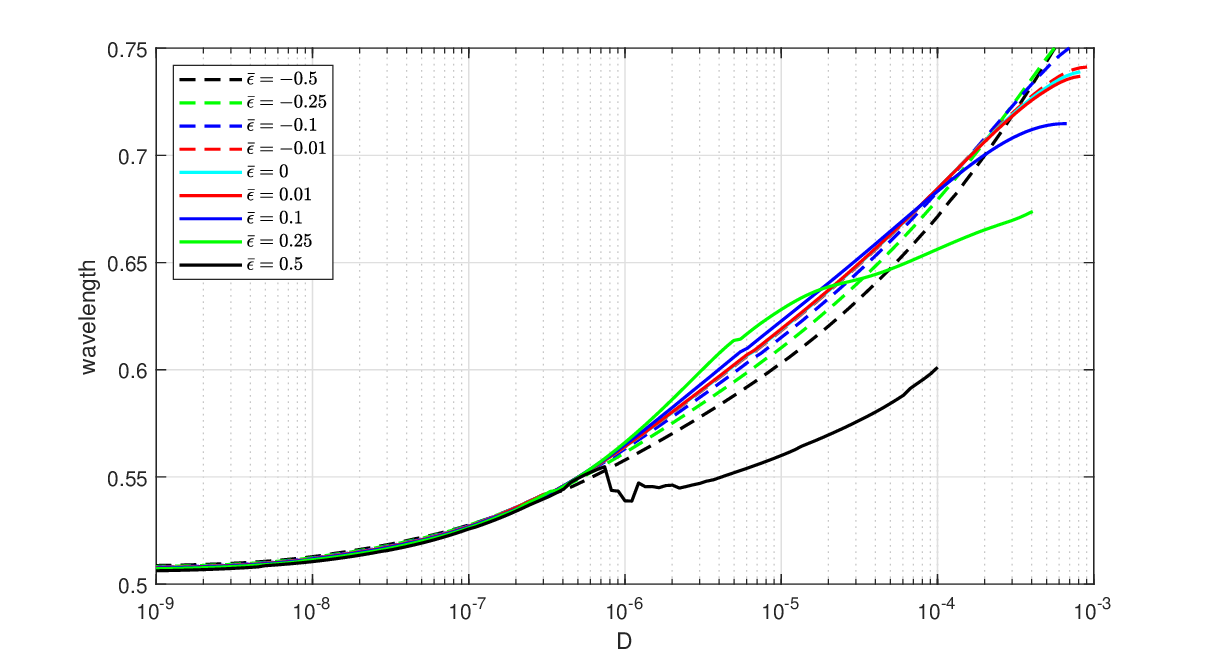}
\caption{The wavelength of the spatially-periodic steady state left behind the wavefront, calculated numerically as a function of $D$, for various values of $\overline{\epsilon}$. The broken lines are the results for negative values of $\overline{\epsilon}$.}\label{fig_wavelength}
\end{center}
\end{figure}

\begin{figure}
\begin{center}
\includegraphics[width=\textwidth]{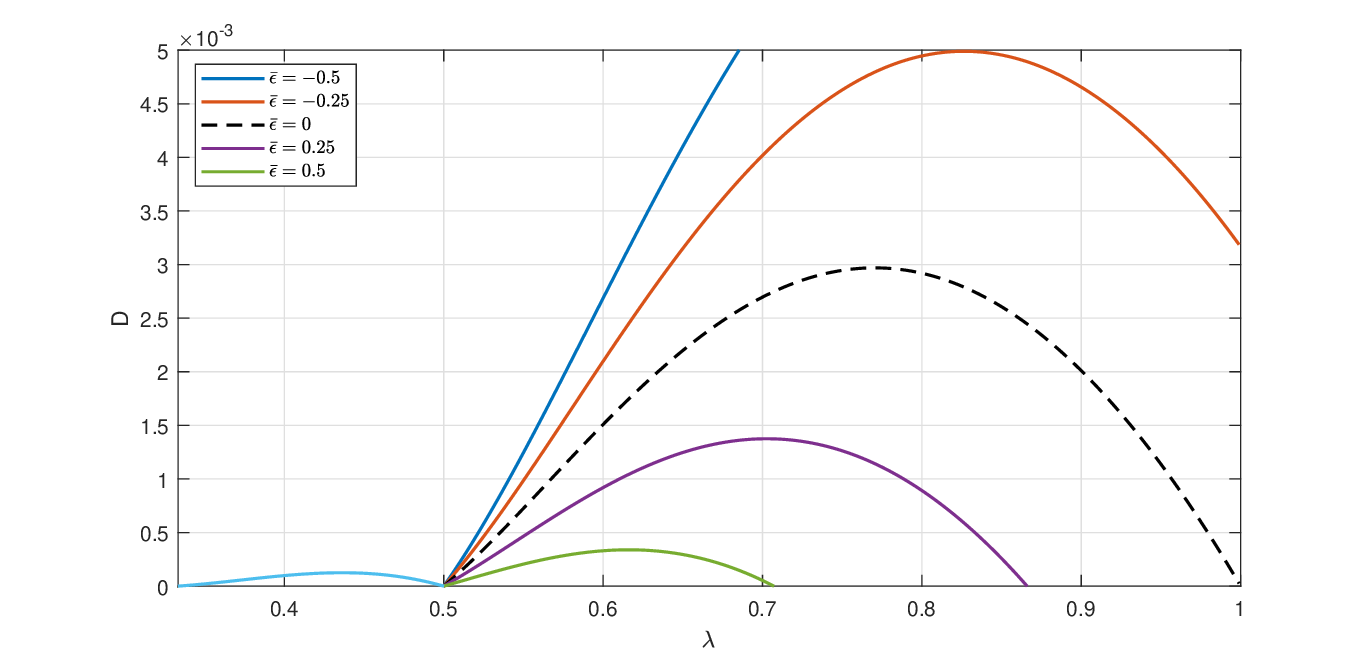}
\caption{The neutral curve given by (\ref{eqn_neutral}) for various values of $\overline{\epsilon}$. The broken line is the neutral curve for the top hat kernel, $\overline{\epsilon} = 0$.}\label{fig_neutral}
\end{center}
\end{figure}

However, of particular interest is how the evolution is affected when $D$ is decreased so that $D\le O(\overline{\epsilon})$, when the above sections suggest that the evolution in (IBVP)$_p$ will now be a singular perturbation from that in (IBVP). In fact, when this is so, the numerical solutions of (IBVP)$_p$ reveal that the evolution mechanism remains close to that when $\overline{\epsilon}=0$, again creating a periodic steady state where the wavelength generated remains uniformly close to that generated when $\overline{\epsilon}= 0$, as shown in Figure~\ref{fig_wavelength}. This means that the complicated bifurcation structure which occurs in the periodic steady states $(\lambda,\overline{D})$ diagram is not involved, since this is located in the region with wavelengths $\lambda\in \left(\frac{3}{4},1\right)$.
It is also clear from Figure~\ref{fig_wavelength} that as $\epsilon$ increases up towards unity, the wavelength generated starts to change for moderate values of $D$, although this effect is not seen for negative values of $\overline{\epsilon}$. As shown in Figure~\ref{fig_neutral}, the first tongue of periodic steady states, given by (\ref{eqn_neutral}), grows as $\overline{\epsilon}$ becomes more negative, shrinks as it becomes more positive and flips into the region $\frac{1}{3} < \lambda < \frac{1}{2}$ for $\overline{\epsilon} > \frac{3}{4}$.
For sufficiently large values of $\overline{\epsilon}$, the dynamics of the solution of (IBVP)$_p$ therefore become qualitatively different, and beyond the scope of the present paper. It is, however, worth noting that for $D$ sufficiently small, even for $|\overline{\epsilon}|$ as large as $0.5$, the same spike creation mechanism acts to generate a periodic steady state of wavelength close to $\frac{1}{2}$, and has the form of those periodic steady states near $\lambda=\frac{1}{2}$ when $\overline{D} =O(1)$ identified in section 3. Note that the slight irregularity of the wavelength generated when $\overline{\epsilon} = 0.5$, shown in Figure~\ref{fig_wavelength}, arises from the disordered and long-lived transient behaviour behind the wavefront for this large perturbation of the top hat kernel. The spike initiation mechanism moves from the region behind the wavefront where $u = O(1)$ to the region ahead of the wavefront where $u$ is exponentially small as $D$ decreases past about $5\times 10^{-7}$ for $\overline{\epsilon} = 0.5$. Typical transient behaviour is shown \href{https://drive.google.com/file/d/1PxvNLm_9ftDHhW901CJDEFqSNl_1KEhD/view?usp=sharing}{in this video} for $D = 10^{-6}$ and $\overline{\epsilon} = 0.5$, but a detailed analysis is beyond the scope of the present paper, whose focus is on universal features of localised kernels, not the specifics of particular large perturbations. Our claim is that the generation of stationary periodic steady states with wavelength $\frac{1}{2}$ (the half-width of the kernel in physical variables) for $D$ sufficiently small is one such universal feature. 

It is important to note here that should we consider the evolution problem now with non-negative \emph{periodic} initial data, the wavelength of the terminal periodic steady state is selected by the chosen wavelength of the initial data. Therefore, when $D\le O(\overline{\epsilon})$, we can select the initial data to bring into contention the complex structure of the steady periodic states created through the multiple steady state bifurcations, which emerge and develop when $D=O(\overline{\epsilon})$ and $\lambda\in \left(\frac{3}{4},1\right)$. In this situation, the evolution is a singular perturbation of that at the same value of $D$, but with $\overline{\epsilon}=0$.

 \section{Conclusions}
 All of our findings have been reviewed in the Introduction, and drawn together at the end of each subsequent section. As such it is unnecessary to repeat those outcomes in detail again. However we note that each of the outcomes throughout the paper support the following very general conclusions:
\begin{itemize}
    \item For admissible kernel perturbations $\overline{\phi} \in \overline{K}(\mathbb{R})$, which also satisfy condition (\ref{eqn3.62}), and have fixed $||\overline{\phi}||_1^m$ small, the dynamics and associated coherent structures of the perturbed evolution problem (IBVP)$_p$, are  a uniform regular perturbation of those associated with the unperturbed evolution problem (IBVP) when the dimensionless diffusivity (being the ration of the diffusion length scale to the nonlocal length scale) $D=O(1)$ as $||\overline{\phi}||_1^m \to 0$. However, this becomes a singular perturbation when $D=O(||\overline{\phi}||_1^m)$ as $||\overline{\phi}||_1^m \to 0$, leading to the possibility of significant structural changes in this limit, however small $||\overline{\phi}||_1^m$ is chosen.

    \item When the kernel perturbation also retains the support [-1/2,1/2], the singular perturbation in the limit $D=O(||\overline{\phi}||_1^m)$ as $||\overline{\phi}||_1^m \to 0$ can lead to complex bifurcation structures regarding periodic steady states, involving the formation, motion and division of humps and spikes. However, a number of key evolutionary mechanisms, and underlying temporal stability properties, are preserved from (IBVP) into (IBVP)$_p$. 
    \end{itemize}

\section*{Funding Statement}
This research did not receive any specific grant from funding agencies in the public, commercial, or not-for-profit sectors.

\section*{Declaration of Interests}
The authors have nothing to declare.

\printcredits

\bibliographystyle{plain}

\bibliography{bibliography}



\begin{appendix}
\section{Weakly nonlinear analysis at the boundary (neutral curve) of $\Omega_1$}\label{appendix}
We seek a solution of (\ref{bif_eqn}) close to the boundary (neutral curve) of the first tongue, $\Omega_1$, which, for the kernel (\ref{bif_kernel}) is given by
\begin{equation}
    D = D_0(\lambda,\overline{\epsilon}) = - \frac{\lambda^3}{4 \pi^3}\left(1-\frac{\overline{\epsilon}}{1-\lambda^2}\right) \sin \left(\frac{\pi}{\lambda}\right).\label{eqn_neutral}
\end{equation}
Note that, since $\lambda<1$, this means that the size of the first tongue, and indeed all the tongues, decreases as $\overline{\epsilon}$ increases, and the small epsilon tongue structure is preserved for all $0\leq \overline{\epsilon} < \frac{3}{4}$. In the first tongue, the longest possible wavelength is $\lambda_{\rm max}(\overline{\epsilon})= \sqrt{1-\overline{\epsilon}}$. In this weakly nonlinear analysis, we will assume that $0 \leq \overline{\epsilon} < \frac{3}{4}$, but not that $\overline{\epsilon}$ is necessarily small.

The weakly nonlinear analysis proceeds in the usual manner \cite{KBO}, by writing $D = D_0(\lambda,\overline{\epsilon}) + \delta^2 D_2$, with $\delta \ll 1$, and expanding $u = 1 + \delta u_0 + \delta^2 u_1 + \delta^3 u_2 + O(\delta^4)$. At $O(\delta)$ we recover the neutral curve and seek a solution of the form 
\begin{equation}
    u_0 = A_0 \cos \left(\frac{2 \pi x}{\lambda}\right).
\end{equation}
At $O(\delta^2)$, there are no secular terms on the right hand side of the equation, and a constant term and a term $\cos \left(4 \pi x/\lambda\right)$ appear in $u_1$. At $O(\delta^3)$, a secular term appears on the right hand side of the equation, and suppressing this leads to the condition
\begin{equation}
    D_2 = -\frac{1}{4} D_0 A_0^2 f(k_1,k_2), 
\end{equation}
where
\begin{equation}
   f(k_1,k_2) = \frac{k_1+k_2}{k_2-4k_1} + 2k_1+2 >0,~~k_n = \frac{\lambda}{n\pi}\left(1-\frac{n^2\overline{\epsilon}}{n^2-\lambda^2}\right) \sin\left(\frac{n\pi}{\lambda}\right).
\end{equation}
Rearranging and eliminating $\delta$ gives the weakly nonlinear periodic steady state, of principal wavelength $\lambda$, as
\begin{equation}
    u \sim 1 + 2 \sqrt{\frac{D_0-D}{D_0 f(k_1,k_2)}}\cos \left(\frac{2 \pi x}{\lambda}\right)\label{WNL}
\end{equation}
for $D_0-D \ll 1$, and confirms that the bifurcation is supercritical ($D < D_0$ is required for solutions other than $u=0$ and $u=1$ to exist locally.). Numerical solutions of (\ref{bif_eqn}) close to the boundary of $\Omega_1$ are in excellent agreement with (\ref{WNL}).
\end{appendix}
\end{document}